\documentclass[a4]{article}

\usepackage[dvips]{graphicx}
\usepackage{amssymb, amsmath, amsthm}

\newtheorem{theo}{Theorem}[section]
\newtheorem{coro}{Corollary}[section]
\newtheorem{lemm}{Lemma}[section]
\newtheorem{prop}{Proposition}[section]
\newtheorem{rema}{Remark}[section]
\newtheorem{defi}{Definition}[section]
\newtheorem{clai}{Claim}[section]

   \def\DD{{\mathbb D}}

 \def\NN{{\mathbb N}}  
 \def\RR{{\mathbb R}}  \def\TT{{\mathbb T}}

\def\De{\Delta}

\def\orb{\mathrm{orb}}
\def\loc{\mathrm{loc}}
\def\mka{\mathfrak{a}}
\def\mkt{\mathfrak{t}}

\def\mkd{\mathfrak{d}}
\def\bfR{\mathbf{R}}
\def\supp{\mathrm{supp}}

\def\cA{{\cal A}}   \def\cM{{\cal M}} 
    
\def\cC{{\cal C}}  \def\cI{{\cal I}} \def\cO{{\cal O}} \def\cU{{\cal U}}
\def\cD{{\cal D}}    
   \def\cQ{{\cal Q}}

\title{A mechanism for ejecting a horseshoe 
from a partially hyperbolic chain recurrence class}

\author{Christian Bonatti and Katsutoshi Shinohara}

\begin{document}

\maketitle

\begin{abstract}
We give a $C^1$-perturbation technique for ejecting 
an \emph{a priori} given finite set of 
periodic points 
preserving a given finite set of 
homo/hetero-clinic intersections
from a chain recurrence class of a periodic point.
The technique is first stated under a simpler setting called 
Markov iterated function system, a
two dimensional iterated function system 
in which the compositions are chosen in Markovian way. 
Then we apply the result to the setting of 
three dimensional partially hyperbolic diffeomorphisms.

{ \medskip

}


\end{abstract}


\section{Introduction}

\subsection{Backgrounds}
When we describe the global structure of dynamical systems,
periodic points play important roles. 
For instance, if we have a hyperbolic periodic point whose 
all eigenvalues have absolute values smaller than one, then 
we immediately know that there is an non-empty open region 
of the phase space which are attracted to the orbit of the 
periodic point.  

For $C^1$-generic chaotic dynamical systems, it is known that
several objects which recapitulate the properties of the system 
are well approximated by periodic orbits. For instance,
for $C^1$-generic diffeomorphisms we have the following (see
\cite{Bo}, \cite{BDV} for more comprehensive accounts on 
backgrounds and references): 
\begin{itemize}
\item Pugh's closing lemma implies that 
the non-wandering set is equal to 
the closure of the set of periodic points. 
\item Ma\~n\'e's ergodic closing lemma implies 
that every ergodic probability measure is the
weak limit of Dirac measures supported on periodic orbits 
which converges to the support of the measure 
in the Hausdorff distance. 
\item More recently, as a consequence of Hayashi's 
connecting lemma, in \cite{BC} it was proved that 
the chain recurrent set is the closure of the set of 
periodic points. More precisely, according to \cite{C}, 
every chain transitive compact set is the Hausdorff limit of 
a sequence of periodic orbits. 
Furthermore, \cite{BC} shows that every chain recurrence class 
of a periodic orbit is indeed its homoclinic class 
(closure of the set of the transverse intersections 
points of its invariant manifolds).
\end{itemize}

However, in \cite{BD1, BD2} it was proved that 
there are open sets of $\mathrm{Diff}^1(M)$ 
in which every $C^1$-generic diffeomorphism has 
uncountably many chain recurrence classes that do not contain periodic points. This phenomenon is not exceptional in the sense 
that it occurs for every $C^1$-generic diffeomorphism having 
a homoclinic class that robustly fails to carry any kind of 
dominated splittings. 
This result naturally leads to the notion 
of \emph{aperiodic classes}, chain 
recurrence classes which do not contain any periodic points.

\bigskip

The previously known constructions 
of $C^1$-locally generic diffeomorphisms 
with aperiodic classes follow essentially
the same process (see for instance \cite{BCDG}): By performing 
successive perturbations of a given diffeomorphisms, we 
first build a nested family of periodic attracting/repelling regions whose components have diameters tending to $0$ and whose
periods tend to infinity. 
The aperiodic class is the intersection of these periodic regions and therefore the dynamics on it conjugates to
an adding machine. 
In particular, in all the known examples, the aperiodic classes  of $C^1$-generic diffeomorphisms are minimal and 
uniquely ergodic. 

The lack of examples of aperiodic classes 
is a huge hindrance 
for the understanding the general behavior 
of $C^1$-generic diffeomorphisms, 
in particular in a neighborhood of aperiodic classes. 
This paper is part of a research, in the sequel to \cite{BS1,BS2},
for building aperiodic classes with totally different behaviors:  
non-unique ergodicity or even non-transitivity. 

Let us briefly see what were done in the previous works.
In \cite{BS1} we defined the notion of 
\emph{$\varepsilon$-flexible 
periodic points} 
and discussed their principal property: 
Its stable manifold in a fixed fundamental domain
can be deformed into an arbitrarily 
prescribed shape by performing an 
$\varepsilon$-perturbation of the diffeomorphism. 
We also showed their $C^1$-generic existence among 
certain kind of partially hyperbolic homoclinic classes for 
arbitrarily small $\varepsilon >0$.
In \cite{BS2}, we introduced the notion 
of \emph{partially hyperbolic filtrating Markov partitions}
which is an assembly of the information about 
partial hyperbolicity and the chain recurrence in a region.
In this setting we showed that, 
if it contains an $\varepsilon$-flexible 
points with a large stable manifold, then it can be ejected 
from the chain recurrence class by performing 
a $C^1$-$\varepsilon$-small perturbation. 
As a consequence, assuming additional information
about the partial hyperbolicity which guarantees 
the abundance of flexible points, we proved that 
the $C^1$-generic diffeomorphisms in the neighborhood 
of a diffeomorphisms having partially hyperbolic 
filtrating Markov partition are 
wild: They admit infinitely many periodic points 
with trivial homoclinic classes (saddles). 
In regard to this type of construction, see also the 
recent progress of Wang \cite{W}, in which the creation of 
weak periodic orbits keeping the connection with 
the initial homoclinic class was discussed.

In this paper, based on these preparations,
we discuss the main technical issues of the project. 
Let us explain it.
Consider a diffeomorphism of a $3$-manifold
admitting a partially hyperbolic filtrating Markov partition.
We assume that it contains 
a finite family of $\varepsilon$-flexible points $\{q_i\}$ 
with large stable manifolds and 
a finite family of homoclinic/heteroclinic points $\{Q_j\}$ 
among $\{q_i\}$. We take another periodic point $p$.
Then we want to 
find an $\varepsilon$-perturbation which
ejects a transitive hyperbolic basic set containing 
the periodic points $\{q_i\}$ and the chosen homo/heteroclinic 
orbits $\{Q_j\}$ in such a way that the chain recurrence
class of the ejected hyperbolic set 
does not contain $p$. The aim of this paper is to 
describe such a perturbation technique. 

We want to find such a perturbation because 
it leads us to construct new kinds of aperiodic classes:
The ejected hyperbolic set admits 
a filtrating Markov partition and admits flexible points. 
Hence, we can inductively proceed this construction.  
We eject nested sequence of hyperbolic sets, 
by smaller and smaller perturbations. 
The aperiodic classes will be obtained as the 
limit of the successively ejected hyperbolic sets. 
Its dynamics depends on 
the choice of the ejected hyperbolic sets. 
By controlling the choice of intermediate dynamics, 
we expect that we can produce aperiodic classes
having a great variety of different dynamical behaviors. 
The confirmation of such properties will be the topic of  
the next paper \cite{BS3}. 

\subsection{Main results}
Let us give the precise statement of our results.
Our main result is about the bifurcation of chain recurrence 
classes appearing near a one having specific conditions. 
We freely use the basic notions of 
topological dynamical systems such as 
attracting/repelling sets (we also use the phrases 
attracting/repelling regions), chain recurrence classes and 
filtrating sets based on the convention \cite{BS2}
(see Section~2.1 of \cite{BS2}).

We first review the notion of partially hyperbolic 
filtrating Markov partitions of saddle type which was introduced in \cite{BS2}. For simplicity, we use the phrase 
``filtrating Markov partitions'' 
in the sense of partially hyperbolic filtrating Markov partitions.
For more information about the 
definition and its basic properties, 
see Section~1.2 and Section~2 of \cite{BS2}.

Throughout this article, $M$ denotes a closed 
(compact and boundaryless) smooth manifold 
of dimension $3$.
A compact subset $C$ of $M$ is said to 
be a \emph{rectangle} if it is $C^1$-diffeomorphic 
to a cylinder $\mathbb{D}^2 \times [0, 1]$. 
$\partial_l C$ denotes the subset of $C$ corresponding
to $\mathbb{D}^2 \times \{ 0, 1\}$, 
called \emph{lid boundary} and
$\partial_s C$ to 
$(\partial \mathbb{D}^2) \times [0, 1]$,
called \emph{side boundary}.

Given a cone field on a rectangle, we 
have the notion of vertical cone field.
A cone field $\mathcal{C}$ on a rectangle $C$ is \emph{vertical} 
if there is a $C^1$-diffeomorphism $\psi$ which sends 
$C$ to the standard cylinder $\mathbb{D}^2 \times [0, 1]$
and $d\phi (\mathcal{C})$ contains $\partial/\partial z$
and is transverse to the plane
$ \langle \partial/\partial x, \partial/\partial y\rangle$, where 
$(x, y, z)$ denote the local coordinate functions of $\phi$.
A cone field is said to be \emph{unstable} if it is strictly invariant (i.e., the image of the closure of the cone field is 
contained in the interior of it on each fiber of the
projective bundle) 
and if every vector in it is uniformly expanding
(see Section~2.2 of \cite{BS2} for details).

Now we are ready to state the definition of filtrating Markov
partitions.
\begin{defi}\label{d.fmp}
Let $f:M \to M$ be a $C^1$-diffeomorphism and
let $\mathbf{R} \subset M$ be a compact set. 
We say that it is a (partially hyperbolic) 
filtrating Markov partition if the following holds:
\begin{itemize}
\item It is a filtrating set: $\mathbf{R} = A \cap R$ for 
an attracting set $A$ and a repelling set $R$.
\item It is a disjoint union of finitely many rectangles:
$\mathbf{R} = \cup C_i$.
\item For each $C_i$, its side boundary is 
contained in $\partial A$ and 
its lid boundary is 
contained in $\partial R$.
\item $\mathbf{R}$ has a vertical,  strictly invariant 
unstable cone field $\cC$.
\item For every $(i, j)$, we have that $f(C_i) \cap C_j$
consists of finitely many vertical rectangles in $C_j$.
Roughly speaking, a sub rectangle $C' \subset C$ is called vertical if it properly crosses $C$ for the precise definition, 
see Section~2.3 of \cite{BS2}.  
\end{itemize}
\end{defi}
A filtrating Markov partition is a capsuled set of information about 
a filtrating set which behaves Markovian way 
keeping the shape of rectangles differential-topologically. 
Recall that having a filtrating Markov partition is a $C^1$-robust 
property 
and for a 
Markov partition $\mathbf{R}$ we have the notion 
of refinements: $f^{-m}(\mathbf{R}) \cap f^{n}(\mathbf{R}) $ turns 
to be a filtrating Markov partition (see Corollary~2.14 of \cite{BS2}). 
We call it the $(m, n)$-refinement of 
$\mathbf{R}$ and denote it by $\mathbf{R}_{(m, n)}$. 
See also Section~\ref{s.elem} of this paper.
We also use the notation $\mathbf{R}_{(m, n ;f)}$ when
we want to indicate the map used to take the refinement.

In \cite{BS1}, we defined the notion of 
an \emph{$\varepsilon$-flexible periodic 
point}. It is a periodic point for which we can find 
a convenient $\varepsilon$-small deformation. 
For the precise definition, 
see Section~\ref{ss:flex} of this paper.
In the same article, we showed that the existence of 
$\varepsilon$-flexible points is abundant 
among chain recurrence classes satisfying certain conditions.
Let us recall the result. In the following, 
$C(p)$ denotes the chain recurrence class of 
a hyperbolic periodic point $p$ of $f$.
\begin{defi}
Let $f$ be a $C^1$-diffeomorphism of a three dimensional 
manifold and $p$ be its hyperbolic periodic point of
stable-index two. Consider the following conditions 
for a chain recurrence class $C(p)$:
\begin{itemize}
\item There is a filtrating Markov partition containing $p$ 
having a large stable manifold. 
We say that a hyperbolic periodic point 
of stable index two in a filtrating Markov 
partition has a large stable manifold 
if $W^s(p)$ cuts the cylinder which $p$ belongs to, see Definition~2.16 of \cite{BS2}.
\item There is a hyperbolic periodic point $p_1$ homoclinically 
related to $p$ such that $p_1$ has a stable non-real eigenvalue.
\item It has a robust heterodimensional cycle  
(see Proposition~5.1 of \cite{BS1} for the definition of 
robust heterodimensional cycles).
\end{itemize}
In this paper, we say that $C(p)$ satisfies condition $(\ell)$ if 
it satisfies all the conditions above.
\end{defi}

In \cite{BS2}, we showed that diffeomorphism 
satisfying the condition ($\ell$) is a wild diffeomorphism
(see Corollary~1.2 of \cite{BS2}).
We proved it by showing that $C^1$-generically 
there is accumulation of isolated saddles nearby. 
The aim of this paper is to show that it has 
stronger pathological behavior. To state it, 
we prepare a definition.

\begin{defi}[\cite{Bo}]
A property $\mathcal{Q}$ about chain recurrence classes
containing a hyperbolic periodic point is called 
\emph{$C^r$-viral} if for every 
$C^r$-diffeomorphism $f$ and 
every hyperbolic periodic point $p$ of $f$ whose chain recurrence class
$C(p; f)$ satisfies property $\mathcal{Q}$ the following holds:
\begin{itemize}
\item There is a $C^r$-neighborhood $\mathcal{U}$ of $f$ such that 
$C(p ; g)$
also satisfies property $\mathcal{Q}$ for every $g \in \mathcal{U}$,
where $C(p; g)$ denotes the chain recurrence class of the 
continuation of $p$ for $g$. In other words, $\cQ$ is 
$C^r$-robust.
\item For every $C^r$-neighborhood $\mathcal{V}$ of $f$ and every 
neighborhood $V$ of $C(p; f)$, there exists $g \in \mathcal{V}$ such that the following holds:
	\begin{itemize}
	\item $g$ has a hyperbolic periodic point $p'$ with 
	$C(p';g) \subset V$, and
	\item $C(p'; g)$ satisfies property $\mathcal{Q}$ and 
	$C(p'; g) \neq C(p; g)$.
	\end{itemize}
\end{itemize}
\end{defi}
Now let us give a main result of this paper.
\begin{theo}
\label{theo:viral}
Having a partially hyperbolic filtrating Markov partition 
containing a chain recurrence class satisfying condition 
($\ell$) is a $C^1$-viral property. 
\end{theo}

By carefully investigating the proof of result, 
we can obtain the following.

\begin{theo}\label{t.aperi}
Let $f \in \mathrm{Diff}^1(M)$ having a
filtrating Markov partition 
containing a chain recurrence class $C(p)$
satisfying property $(\ell)$.
Consider an open neighborhood 
$\mathcal{O} \subset \mathrm{Diff}^1(M)$ of $f$ 
where we can define the continuation of $C(p)$
keeping the property $(\ell)$.
Then, $C^1$-generic diffeomorphisms in $\mathcal{O}$
has an aperiodic class.
\end{theo}


\subsection{Main result with precise information}
\label{ss.results}
Theorem~\ref{theo:viral} implies the creation of new chain 
recurrence classes for filtrating Markov partitions 
with condition ($\ell$) up to a $C^1$-small perturbation. 
While this result is easy to understand, our construction 
indeed gives more information about the structure of the 
chain recurrence classes ejected. In this subsection, 
we formulate it. 

We begin with a general definition. 
Let $f$ be a $C^1$-diffeomorphism of $M$.
For a point $q$, by $\cO(q)$ we denote the 
orbit of $q$,
namely, $\cO(q) = \{f^i(q)\}_{i \in \mathbb{Z}}$.
By a \emph{circuit of points} we mean a collection of
finitely many 
hyperbolic periodic orbits $\{\mathcal{O}(q_i)\}$ and 
finitely many 
transverse homo/heteroclinic orbits $\{ \mathcal{O}(Q_i)\}$ 
connecting among them. Consider a directed graph whose 
vertices are periodic orbits $\{\mathcal{O}(q_i)\}$
and edges are the collection of 
homo/heteroclinic orbits 
connecting the vertices. 
If we consider a $C^1$ diffeomorphism $g$
sufficiently $C^1$-close to $f$, 
then we can consider 
the continuation of $S$, which we denote by $S_g$.

We say that a circuit of points is \emph{transitive} if 
its corresponding directed graph is transitive 
(every two vertices can be connected by a sequence of edges).
In this article, we only consider circuits which are transitive. 
Thus, throughout this paper by a circuit we mean a transitive one.
We define a similar notion for filtrating sets.
Given a filtrating set $R$, suppose that it has finitely many 
connected components. This is always the case for filtrating Markov partitions. 
We say that $R$ is \emph{$\mathrm{c}$-transitive} if 
for every pair of connected components $R_1$, $R_2$ of 
$R$, there is a sequence of components $(S_i)_{i=1,\ldots, k}$
such that $S_1=R_1$, $S_k=R_2$ and 
$f(S_i) \cap S_{i+1} \neq \emptyset$ holds for every 
$i=1,\ldots, k-1$.

Let $f$, $g$ be $C^1$-diffeomorphisms of $M$.
Let $\Lambda_f$, $\Lambda_g$ 
be $f$, $g$-invariant subset of $M$ respectively
and $\delta >0$. 
We say that $\Lambda_f$ and $\Lambda_g$ 
are \emph{$\delta$-similar}
if there is a homeomorphism 
$h:\Lambda_f \to \Lambda_g$ which is
$C^0$-$\delta$-close to the identity
(that is, for every $x \in \Lambda_f$ we have
$d(x, h(x)) < \delta$ holds, where $d$ is a distance 
function) such that $h$ is a 
conjugacy between $f$ and $g$, that is, 
$ h \circ f = g \circ h$ holds on $\Lambda_f$.

Suppose that
we have a hyperbolic set $\Lambda \subset C(p)$. 
We say that  $\Lambda$ is 
\emph{$\varepsilon$-coarsely expulsible from $C(p)$} 
if the following holds: 
For any $\delta >0$ there is a $C^1$-diffeomorphism 
$g= g_{\delta}$ which 
is $C^1$-$\varepsilon$-close to $f$ such that the following holds: 
\begin{itemize}
\item There is a $g$-invariant set $\Lambda_g$ such that 
$\Lambda_f$ and $\Lambda_g$ are $\delta$-similar.
\item There is a filtrating set $R'$ containing $\Lambda_g$ 
such that $R'$ does not contain $p$. 
\end{itemize}

We say that $\Lambda$ is 
\emph{$\varepsilon$-expulsible from $C(p)$}
if $R'$ can be chosen arbitrarily close to $\Lambda$, 
that is, for any neighborhood $U'$ of $\Lambda$
we can choose $g$ in such a way that $R'$ is contained 
in $U'$. 

We introduce one more definition 
about filtrating Markov partitions. 
It is called the \emph{robustness}.
Roughly speaking, a filtrating Markov partition is 
said to be $\alpha$-robust if 
persists for every $C^1$-diffeomorphism $g$ which 
is $\alpha$-close to $f$ and 
coincides with $f$ outside the filtrating Markov partition. 
We give the precise definition in the 
next section, see Definition~\ref{d.robu}.

Then, our refined statement is the following:
\begin{theo}\label{theo:expu}
Let $f$ be a $C^1$-diffeomorphism of a closed 
three manifold having a 
chain recurrence class $C(p)$ contained in a 
filtrating Markov partition $\mathbf{R}$ which is 
$\alpha$-robust. Let $S \subset \mathbf{R}$ be 
a circuit of points which does not contain 
$\mathcal{O}(p)$. Assume that every periodic 
orbit of $S$ is $\varepsilon$-flexible (where $\varepsilon$
satisfies $2\varepsilon <\alpha$) and has a large stable manifold.
Then, $S$ is $2\varepsilon$-expulsible with a filtrating set 
$\mathbf{R}'$ which is also a filtrating Markov partition. 
Furthermore, $\mathbf{R}'$ can be chosen in such a 
way that it is $\mathrm{c}$-transitive, 
$(\alpha -2 \varepsilon)$-robust and every periodic orbit
of $S$ is $\varepsilon$-flexible with a large stable manifold 
in $\mathbf{R}'$.
\end{theo}

We give one more statement of the expulsion result containing
more information about the new filtrating Markov partition 
and the process of the perturbation. 
In the following, given a filtrating Markov partition, 
we often consider not the whole set of rectangles but 
a sub family of rectangles. We call them a \emph{sub Markov partition of} $\mathbf{R}$.
Note that a sub Markov partition may fail to be a filtrating set. 
We are mainly interested in the sub Markov partition 
$\mathbf{R}(S)$, where $\mathbf{R}(S)$ denotes 
the set of rectangles
having non-empty intersection with $S$. 

In section~\ref{s.elem}, we define the notion of 
\emph{affine Markov partitions}. It roughly means that 
the dynamics restricted to there is given by affine maps 
and the shape of the cylinders respects the affine structures.
It enables us to investigate the bifurcation of the dynamical 
systems there in terms of two dimensional dynamics. 
For the precise definition, see Definition~\ref{d.affi}.

In the following, by the \emph{support of a diffeomorphism $g$ with respect to $f$}, denoted by $\mathrm{supp}(g, f)$,
 we denote the closure of the set 
 $\{x \in M \mid f(x) \neq g(x)\}.$ 
The next result is the first step of the proof.

\begin{theo}\label{t.affine}
Let $f$ be a $C^1$-diffeomorphism 
having a filtrating Markov partition $\mathbf{R}$ containing
a circuit of points $S$ such that every periodic orbit of $S$ has 
a large stable manifold. Then for any neighborhood 
$W$ of $S$ and any $C^1$-neighborhood $\cU$ of $f$ 
there is a diffeomorphism $f_1 \in \cU$ 
such that the following holds:
\begin{itemize}
\item The support $\mathrm{supp}(f_1, f)$ is contained in
 $W$.
\item For $f_1$, all the orbit of $S_{f_1}$ has 
the same orbit with the same derivatives along the 
orbits as $S$.
\item 
For every 
sufficiently large $m$ and $n$
(to be precise, there are $m_0, n_0$ such that if 
$m\geq m_0$ and $n \geq n_0$ then)
 the rectangles 
$\mathbf{R}_{(m, n; f_1)}(S)$ is an affine Markov partition.
\end{itemize}
\end{theo}
Thus, roughly speaking, up to an arbitrarily small perturbation 
which preserves the local dynamics along the 
periodic orbits we may assume that 
for a sufficiently fine refinement we have affine property.

The following is one of the 
main steps of the proof of Theorem~\ref{theo:expu}.
For a filtrating Markov partition, 
$\mathbf{R}$,
we say that it is \emph{generating} if 
for any two rectangles $C_1, C_2$ of $\mathbf{R}$,
$f(C_1) \cap C_2$ has at most one connected component.
See Section~\ref{ss.gene} for more information.

\begin{theo}\label{t.circ}
Let $f \in \mathrm{Diff}^1(M)$ having 
a circuit of points $S$ in a generating filtrating Markov partition
$\mathbf{R}$
such that every periodic orbit is $\varepsilon$-flexible 
and has a large stable manifold in $\mathbf{R}$.
Assume that $\mathbf{R}$ is
$\alpha$-robust for $\alpha > 2 \varepsilon$
and $\mathbf{R}(S)$ is an affine Markov partition. 
Then 
for every sufficiently large $n$, there is a diffeomorphism 
$f_n$ which is $2\varepsilon$-close to $f$ 
and whose support $\mathrm{supp}(f_n, f)$ is contained in 
the interior of $\mathbf{R}(S)$ such that the 
following holds:
\begin{itemize}
\item $f_n$ has a transitive filtrating Markov partition
$\mathbf{R}_n'$ 
containing a circuit of points $S_n$ which 
is similar to $S$ and satisfying the following:
\begin{itemize} 
\item We have 
$\mathbf{R}(S) = \mathbf{R}(S_n)$ and we can require 
that the points of $S$ and $S_n$ which are conjugated 
under the conjugacy belong to the same rectangle.
\item All the periodic orbits of $S_n$ have large stable manifolds
in $\mathbf{R}_n'$ and they are all $\varepsilon$-flexible.
\item The periodic orbits of $S_n$ have the same orbits as $S$.
\end{itemize}
For $\mathbf{R}_n'$, we have
\begin{itemize}
\item Each rectangle of $\mathbf{R}_n'$ is 
a vertical sub rectangle of some 
rectangle of $\mathbf{R}_{(0, n; f_n)}(S_{n})$.
Especially, $\mathbf{R}_n'$ is contained in
$\mathbf{R}_{(0, n; f_n)}(S_{n})$. 
\item Each rectangle of $\mathbf{R}_{(0, n; f_n)}(S_{n})$
contains one and only one rectangle of $\mathbf{R}_n'$.
\item The cone field of $\mathbf{R}_n'$ is the restriction
of the one of $\mathbf{R}_{(0,n;f_n)}$.
In particular, $\mathbf{R}_n'$ is $(\alpha -2 \varepsilon)$-robust.
\end{itemize}
\end{itemize}
\end{theo}

Note that, in general $f_n$ is so far from 
$f$ that we may fail to have a continuation of $S$.
This theorem claims the non-trivial existence 
of the continuation of $S$.

This result, together with the abundance result of the 
flexible points with large stable manifolds implies 
Theorem~\ref{theo:viral}. We will discuss the derivation
of Theorem~\ref{theo:viral} in Section~\ref{s.3dim}.

%

Let us briefly see the idea of the proof.
The proof is divided into two steps.
The first part is that we describe such perturbation results
in the context of \emph{Markov iterated function systems} (referred as Markov IFSs), 
which are the abstraction of the information of 
affine Markov partitions. 
Theorem~\ref{t.main} (see Section~\ref{s.ifs})
is the technical core of this paper. It states that the ejection 
described above is possible in the level of IFSs: 
Given a circuit consisting of $\varepsilon$-flexible periodic 
orbits with large stable manifolds related
by heteroclinic/homoclinic orbits, one can eject 
a hyperbolic set containing this circuit away from 
a given class, by an $\varepsilon$-perturbation.
Then, we will transfer 
the result for iterated function systems to 
filtrating Markov partitions 
by giving a perturbation technique 
which reduces to the original problem to the study of 
Markov IFSs. We will discuss more about the proof of 
Theorem~\ref{t.main} later (see Section~\ref{ss.back})

\subsection{Organization of this paper}
Finally, let us explain the structure of this paper. 
In Section~\ref{s.elem}, we introduce several notions 
related to rectangles of filtrating Markov partitions 
containing a circuit of points. We discuss the effect of 
taking refinements for such rectangles. We also give a 
of linearization result 
(Theorem~\ref{t.affine}) for the dynamics around them. This 
enables us to reduce the proof of 
Theorem~\ref{theo:viral} and Theorem~\ref{theo:expu}
into the problem of two dimensional dynamics.
In Section~\ref{s.ifs}, we introduce the definition of 
Markov IFSs and discuss their elementary properties
such as their periodic orbits, refinements, 
attracting/repelling regions. Based on these preparations
we give the statement of Theorem~\ref{t.main}, 
the main perturbation result stated in terms of Markov IFSs. 
In Section~\ref{s.3dim}, 
after preparing several preliminary perturbation techniques
which are essentially given in past papers \cite{BS1, BS2},
we prove Theorem~\ref{theo:expu} and \ref{t.circ}
admitting Theorem~\ref{t.main}. 
We also see how we derive 
Theorem~\ref{theo:viral} and \ref{t.aperi}
from Theorem~\ref{theo:expu} and \ref{t.circ}. 
The rest of the paper is 
dedicated to the proof of Theorem~\ref{t.main}.
We first prove Theorem~\ref{t.IFS}, which is a simplified 
version of Theorem~\ref{t.main}.
In Section~\ref{s.game}, 
we introduce several notions such as 
retarded family, wells and obstructions. 
They are extractions of some important 
information of the Markov IFSs 
for the proof of Theorem~\ref{t.IFS}.
In Section~\ref{s.solution},
we complete the proof of Theorem~\ref{t.IFS}. 
Finally, in Section~\ref{s.expla} we explain how to 
deduce Theorem~\ref{t.main} from the proof of 
Theorem~\ref{t.IFS}.
\bigskip

{\bf Acknowledgments.} This work is supported by the 
 JSPS KAKENHI Grant Numbers 18K03357.
KS is grateful for the hospitality of Institut de Math\'ematiques de Bourgogne of Universit\'e de Bourgogne during his visit.


\section{Local linearization of Markov partitions}\label{s.elem}
In this section, we prove several elementary results which reduce 
the investigation of Markov partitions into a simpler one
up to small perturbations. 

\subsection{Basic notions and refinements}
For a filtrating Markov partition $\mathbf{R}= \cup C_i$ and 
a point $x \in \mathbf{R}$, we denote the (unique)
rectangle containing $x$ by $C_x$.

Consider a periodic point contained in $\mathbf{R}$.
Note that the assumption that $\mathbf{R}$ is 
a filtrating set implies that the orbit is contained in 
the interior of $\mathbf{R}$. 
In general, a periodic 
orbit may contain two points which belongs to the same 
rectangle. For us it would be convenient if each point 
belongs to different rectangles. 
To formulate this we prepare a definition.
In the following, for a hyperbolic 
periodic point $q$ in $\mathbf{R}$, 
by $W_{\mathrm{loc}}^u(q)$
(resp. $W_{\mathrm{loc}}^s(q)$) we denote the connected
component of $W^u (q) \cap C_q$
(resp. $W^s(q) \cap C_q$)
containing $q$.
\begin{defi}
Let $\mathbf{R} = \cup C_i $ be a filtrating Markov partition. 
Suppose that $\mathbf{R}$ contains a periodic point $q$
of period $\pi$. 
We say that $\{C_{f^i(q)}\}_{i=0,\ldots,\pi -1}$ is 
a cycle of periodic rectangles for $q$ if $\{C_{f^i(q)}\}$ ($i=0,\ldots,\pi -1$)
are mutually different.
\end{defi}
In the following, we use the alphabet $K$ to notify
that it is a rectangle in a cycle of periodic rectangles.
Let us define a similar notion for homo/heteroclinic points.
\begin{defi}
\label{d.deparr}
Let $\mathbf{R}$ be a filtrating Markov partition satisfying 
the following:
\begin{itemize}  
\item $\mathbf{R}$ contains a periodic point $q_{j}$ 
and there is a cycle of periodic rectangles
$\{K_{f^i(q_j)}\}_{i=0,\ldots,\pi_{j}-1}$, where 
$\pi_j$ is the period of $q_j$, for $j=1, 2$.
\item There is a point 
$Q \in W_{\mathrm{loc}}^u(f^d(q_1)) \cap W^s(q_2)$, 
where $d$ is some integer,
such that $\{f^i(Q)\}_{i=1,\ldots T-1}$ is 
disjoint from any rectangles $\{K_{f^i(q_j)}\}$ and 
$f^T(Q) \in W_{\mathrm{loc}}^s({f^{a}(q_2)}) \cap W^u(q_1)$ holds for some integer $a$
and $T>0$. 
\end{itemize}
We say that the family of rectangles 
$\{C_{f^i(Q)}\}_{i=1,\ldots, T-1}$ is 
\emph{a path of transition rectangles} 
if they are all distinct.
The integer $T$ is called its \emph{transition time}.

We call $K_{f^d(q_1)}$ the \emph{departure rectangle} and 
$K_{f^a(q_2)}$ the \emph{arrival rectangle}.
Note that we do not exclude the case where $q_1$ and $q_2$
has the same orbit.
\end{defi}

In the following, we consider the set of periodic rectangles
connected by paths of transition rectangles. Let us formulate 
it. 

\begin{defi}\label{d.cirrec}
Let $S$ be a circuit of points consisting of a set of periodic 
points $\{q_j\}$ and homo/heteroclinic points $\{Q_l\}$
contained in a filtrating Markov partition $\mathbf{R}$. 
We say that the sub Markov partition $\mathbf{R}(S)$ 
(see Section~\ref{ss.results} for the definition)
is a \emph{circuit} of rectangles for $S$ 
if the set of 
rectangles $\{K_{f^i(q_j)}\}_{0\leq i \leq \pi_j-1}$
are mutually disjoint cycles of periodic rectangles and 
the rectangles
$\{L_{f^i(Q_l)}\}_{1 \leq i \leq T_l -1 }$
are paths of transition rectangles.
\end{defi}

\begin{rema}\label{r.circ-str}
When we consider a circuit of rectangles, it may be that 
two paths have common rectangles.
We only require that for each path the rectangles are 
distinct, 
and we do not require such 
conditions among two different paths. 
\end{rema}

Recall that for a filtrating Markov partition $\mathbf{R}$,
we can define a new one by taking \emph{refinements}.
Let $\mathbf{R}= A \cap R$ where
$A$ is an attracting set and $R$ a repelling set. 
Then, the set
$\cap_{k=-m}^n f^k(\mathbf{R})$ turns to be 
a filtrating 
Markov partition with an attracting set $f^n(A)$ and 
a repelling set $f^{-m}(R)$
(see Corollary~2.14 of \cite{BS2}). In this paper, 
we call it \emph{an $(m, n)$-refinement} of $\mathbf{R}$
and denote it by $\mathbf{R}_{(m, n)}$ or $\mathbf{R}'$
when we do not need to indicate $(m, n)$. 
We write $\mathbf{R}_{(m, n; f)}$ when we want to 
indicate with which map we took the refinement.

If $\mathbf{R}$ has a cycle of periodic rectangles or 
a path of transition rectangles, 
then one can naturally associate
new ones for $\mathbf{R}'$. 
Suppose we have a cycle of periodic rectangles 
$\{K_{f^i(q)}\}_{i=0,\ldots,\pi-1}$ in a filtrating Markov 
partition $\mathbf{R}$.
Then take its $(m, n)$-refinement. 
In the refinement, there are rectangles containing 
$f^i(q)$ ($i=0,\ldots,\pi -1$) and one can check that
they form a cycle of periodic rectangles for $q$ as well. 
We call it the \emph{corresponding cycle of 
periodic rectangles} in the refinement.
We denote the corresponding rectangles by $\{K'_{f^i(q)}\}$.
Similarly, consider a filtrating Markov partition 
with cycles of periodic rectangles
$\{ K^j_{f^i(q_j)}\}_{j=0,\ldots,\pi_j-1}$
($j=1,2$) and a path of transition
rectangles
$\{L_{f^i(Q)}\}_{i=1,\ldots, T-1}$ with 
respect to a homo/heteroclinic point $Q$
having the departure rectangle $K^1_{f^d(q_1)}$
and the arrival rectangle $K^2_{f^a(q_2)}$.
For the $(m, n)$-refinement, by considering homo/heteroclinic
points $f^{-m}(Q)$ we have transition rectangles with
transition time $T+m+n$, with 
the departure rectangle $(K^1)^{\prime}_{f^{d-m}(q_1)}$
and the arrival rectangle $(K^2)^{\prime}_{f^{a+n}(q_2)}$.
We call the such rectangles 
\emph{corresponding transition rectangles}.
Note that similar construction holds for a circuit of rectangles.

\subsection{Choosing rectangles}
In this subsection, we prove that for a circuit of points 
in a filtrating Markov partition,
by taking some refinements we can obtain a circuit 
of rectangles, if every periodic point has 
a large stable manifold.

\begin{prop} \label{prop:rec_choix}
Let $\mathbf{R}= \cup C_i$ be a 
filtrating Markov partition and $S$ be 
a circuit of points 
such that every periodic point has 
a large stable manifold. Then for every sufficiently large 
$m$ and $n$, $\mathbf{R}_{(m, n)}(S)$ 
is a circuit of rectangles.
\end{prop}

\begin{proof}
We first prove that in every sufficiently fine refinement
the periodic points in the periodic orbit 
are in distinct rectangles. 

Notice that taking backward refinements makes 
the height of the rectangles uniformly (exponentially)
small. Also, the largeness of the stable manifold of periodic points
ensures that taking forward refinement makes the 
width of the rectangles which the periodic points 
belong to uniformly small. Thus, by taking 
sufficiently large (both in forward and backward) refinements
we can assume that the rectangles which 
contain a point of an orbit of a periodic point 
are uniformly small. In particular,
none of them can coincide. 

Now let us see how to construct the paths of transition rectangles.
First, by the argument of the first step we can assume that 
for each rectangle of a cycle of periodic rectangles, they 
contain homo/heteroclinic points only in the local 
stable/unstable manifolds by taking sufficiently fine refinements.
Then, notice that while taking refinements increases 
the number of rectangles, each 
newly created rectangles are contained 
in the initial periodic rectangles and 
they contain at most one homo/heteroclinic 
point of a given homo/heteroclinic orbit.
Thus we only need to prove that by taking refinements we can
separate the points of homo/heteroclinic points outside
the periodic rectangles. 
Since the diameter of periodic rectangles tends 
to zero by taking refinements and the 
transition rectangles in the refinements are 
images of periodic rectangles, we can 
assume that the diameters of 
rectangles containing the point of  
homo/heteroclinic points outside the periodic rectangles 
also goes to zero as we take finer 
refinements. Thus they are separated into different rectangles 
for every sufficiently fine refinements.

This completes the proof.
\end{proof}

\begin{rema}\label{r.tiny}
This proof shows that if every periodic orbit of 
$S$ has a large stable manifold, then given a 
neighborhood $W$ of $S$, for every sufficiently large
$(m, n)$ we have $\mathbf{R}_{(m, n)}(S) \subset W$.
\end{rema}

\subsection{Refinements and itinerary}
\label{ss.gene}
Let $S$ be a circuit of points 
contained in a filtrating Markov partition
$\mathbf{R}$ such that $\mathbf{R}(S)$ is a circuit 
of rectangles. As in Remark~\ref{r.circ-str} it may be that 
two paths of transition rectangles of $\mathbf{R}(S)$ 
share some rectangles. In this subsection, we show that
two different homo/heteroclinic points cannot have totally the 
same itinerary under a mild condition.

First, let us recall the definition of the generating property
which we defined in Introduction.

\begin{defi}
Let $\mathbf{R}$ be a filtrating Markov partition. 
We say that 
$\mathbf{R}$ is \emph{generating} if 
for any two rectangles $C_1, C_2$ of $\mathbf{R}$,
$f(C_1) \cap C_2$ has at most one connected component.
\end{defi}

\begin{lemm}\label{l.gene}
For a filtrating Markov partition, its refinement is generating.
\end{lemm}

\begin{proof}
We prove it for $(0, 1)$-refinement.
The general case is similar.
Let $C'_1$, $C'_2$ be rectangles in 
$\mathbf{R}\cap f(\mathbf{R})$ and assume 
$f(C'_1) \cap C'_2 \neq \emptyset$. We will show
$f^{-1}(f(C'_1) \cap C'_2) =  C'_1 \cap f^{-1}(C'_2)$
is connected. Let us take a rectangle $C_1$ 
of $\mathbf{R}$ which contains $C'_1$. Then, 
$C'_1$ is a vertical subrectangle of $C_1$ and 
$f^{-1}(C'_2)$ is a horizontal subrectangle of $C_1$
(see Section~2.3 of \cite{BS2} for the detail).
Thus they intersects and there is unique connected 
component. 
\end{proof}

\begin{lemm}\label{l.shadow}
Let $S$ be a circuit of points contained in
a generating filtrating Markov partition $\mathbf{R}$. 
Let $x, y$ be its homo/heteroclinic orbits from 
$q_1$ to $q_2$, where $q_1, q_2$ are periodic points in $S$.  
If the transition time of $x$ and $y$ are the same 
and $f^i(x)$ and $f^i(y)$ belong to the same 
rectangle for every $0 \leq i \leq T$ where $T$
is the common transition time, then $x = y$. 
In other words, two different homo/heteroclinic
orbit must have different itineraries.
\end{lemm}

\begin{proof}
Let $\sigma$ be the local unstable manifold 
of $f^{d}(q_1)$. It 
contains $x$ and $y$. Let $C_i$ be the rectangle 
which contains $f^i(x)$ and $f^i(y)$ for $i=0,\ldots,T$.  
By using the invariance of the cone field and the 
generating property of $\bfR$, we see that the 
part of the $\sigma$ whose image under $f$ is in $C_1$
is a connected curve. Inductively, the part
of $\sigma$ whose image under $f^i$ is in $C_i$ 
for every $0\leq i \leq k$ is a connected curve for 
every $0\leq k \leq T$. Let us denote the curve by $\sigma_k$, and consider $\sigma_T$. By definition, 
we know $x, y \in \sigma_T$ and $f^T(\sigma_T)$ has 
unique intersection with the local stable manifold of 
$f^a(q_2)$ due to the invariance 
of the cone field. Thus we know 
$f^T(x) = f^T(y)$ and consequently $x=y$.
\end{proof}

\subsection{Linearization of periodic rectangles}
In this and next subsections, we discuss perturbation 
techniques which transform the dynamics near 
a circuit of rectangles into a simpler form. 

Let us prepare some definitions.
We say that a compact set in $\mathbb{R}^3$
is a \emph{product rectangle} if it contains the origin in the 
interior and has the form
$D \times I$
where $D \subset \mathbb{R}^2$ is a compact set 
$C^1$-diffeomorphic to the round disc $\mathbb{D}^2$ and
$I$ is a closed interval.
\begin{defi}\label{d.lincyc}
Let $\mathbf{R}=\cup C_i$ be a filtrating Markov partition of 
a diffeomorphism $f$
and $\{K_{f^i(q)}\}$ be a cycle of periodic rectangles of 
a periodic point $q$. We say that the cycle
$\{K_{f^i(q)}\}_{i=0,\ldots,\pi-1}$ 
is \emph{linearized} if the following holds:
\begin{itemize}
\item For each $0\leq i\leq \pi -1$, there exists a coordinate
neighborhood $(U_i,\phi_i)$ containing $K_{f^i(q)}$
such that $\phi_i(K_{f^i(q)})$ is a product rectangle
$D_i \times I_i \subset \mathbb{R}^2 \times \mathbb{R}$.
We set $(U_\pi,\phi_\pi) :=(U_0,\phi_0)$.
\item For each $0\leq i\leq \pi -1$, 
let $J_{f^i(q)}$ be the 
connected component of 
$K_{f^i(q)} \cap f^{-1}(K_{f^{i+1}(q)})$
containing $f^i(q)$.
Then the map $\phi_{i+1} \circ f \circ \phi^{-1}_{i}$
restricted to $\phi_i(J_{f^i(q)})$ 
is an affine map preserving the product structure 
$\mathbb{R}^2 \times \mathbb{R}$ in such 
a way that $\mathbb{R}^2$, $\mathbb{R}$
corresponds to $E^{cs}$, $E^u$ directions, respectively.
\end{itemize} 
\end{defi}
\begin{rema}
For a linearized cycle, if we take a refinement then 
the corresponding cycle is also linearized. 
\end{rema}
The following result says that up to an arbitrarily $C^1$-small
perturbation and a refinement one can make a cycle
of periodic rectangles with a large stable manifold 
being linearized.
\begin{prop}\label{prop:pelipeli}
Let $\mathbf{R}=\cup C_i$ be a filtrating 
Markov partition of a diffeomorphism $f$
and $\{K_{f^i(q)}\}$ be a cycle of periodic rectangles of 
a periodic point $q$. Suppose that $q$ has a large 
stable manifold. Then, for every $C^1$-neighborhood 
$\cU$ of $f$ and
every neighborhood $W$ of $\mathcal{O}(q)$
there exists a diffeomorphism 
$g \in \cU$ such that the following holds:
\begin{itemize}
\item $q$ is a periodic point of $g$ with
the same orbit as for $f$, and
the derivatives of $f$ and $g$ along $\cO(q)$ 
are the same.
\item The support $\mathrm{supp}(g, f)$ 
is contained in $W$. 
\item For $g$, $\mathbf{R}$ is a filtrating Markov partition 
such that for every sufficiently large $(m, n)$,
the corresponding cycle of periodic rectangles
in $\mathbf{R}_{(m,n;g)}$ 
is linearized.
\end{itemize}
\end{prop}

For the proof, we need the Franks' lemma, which enables 
us to linearize the dynamics locally. 
See for instance \cite{BD-sta} for more information.
\begin{lemm}[Local linearization by Franks' lemma]\label{llin-map}
Let $f \in \mathrm{Diff}^1(M)$, $\dim M = m$, $x \in M$
and $\phi : U \to \mathbb{R}^m$, $\psi : V \to \mathbb{R}^m$
be two coordinate neighborhood of $x$, $f(x)$ respectively
such that $\phi(x)$, $\psi(f(x))$ are the origin of  $\mathbb{R}^m$.
Then for any $\varepsilon >0$ and any neighborhood $U'$ of $x$, 
there exist a neighborhood $\tilde{U}$ of $x$ contained in $U'$ 
and $\tilde{f} \in \mathrm{Diff}^1(M)$ such that 
$\tilde{f}$ is $\varepsilon$-$C^1$-close to $f$,
$\tilde{f}$ coincides with $f$ on $M \setminus U'$ and
the map $\psi \circ \tilde{f}\circ \phi^{-1}$ coincides with 
a linear map given by $d(\psi \circ f \circ\phi^{-1})$ on $\tilde{U}$.
\end{lemm}

\begin{proof}[Proof of Proposition~\ref{prop:pelipeli}]
The proof is similar to \cite[Proposition~3.6]{BD-sta}.

First, we apply the Franks' lemma along the orbit of $q$.
More precisely, we take a diffeomorphism $f_1$
that satisfies the following:
\begin{itemize}
\item For each $j \in \mathbb{Z}$, 
$f^j_1(q) = f^j(q)$. 
Especially, $q$ is a 
periodic point for $f_1$ with the same orbit. 
\item $\mathrm{supp}(f_1, f) \subset W$.
\item For each $j$, there exists a coordinate
neighborhood $V_j$ of 
$f_1^j(q) = f^j(q)$ such that the dynamics of $f_1$ 
on $V_j$ is given by the linear map $df(f^j(q))$.
\end{itemize} 
We fix such $f_1 \in \cU$.
This is the only part we perform the perturbation 
along $\cO(q)$. 
Note that this  
does not change the derivatives along $\mathrm{orb}(q)$.

By choosing $f_1$ sufficiently close to $f$, 
we may assume that $\mathbf{R}$ is still a filtrating 
Markov partition and
$q$ has a large stable manifold for $f_1$ as well.
Now let us take the refinements. As we take finer 
refinements, due to the fact that $q$ has a 
large stable manifolds, the rectangles which $f_1^i(q)$
belongs to shrinks to $\{f_1^i(q)\}$.
Thus we may assume that each corresponding periodic 
rectangle is contained in the linearized coordinates.
Notice that in this coordinates 
$f_1^i(q)$ is mapped to the origin.
We assume that the $xy$-plane coincides with $E^{cs}$ 
direction and the $z$-axis coincides with the $E^u$ direction.

We show that by taking the refinement and slightly 
perturbing $f_1$,
we have that the periodic rectangles are product rectangles
in the linearized coordinates.
First, let us see how to make 
lid boundary flat.
For each $i$, the intersection between 
$W_{\mathrm{loc}}^u(f_1^i(q))$ and 
$\partial_l(K_{f^i(q)})$ consists of two points.
We denote them by $y_{i, +}$ and $y_{i, -}$.
Now we perform a $C^1$-perturbation whose support 
is contained in a small neighborhood of
$f_1^{-1}(y_{i, \pm})$  
so that for $f_2$ (the perturbed diffeomorphism),
the image of the lid boundary near $y_{i \pm}$ 
is parallel to the $xy$-plane.
Note that the size of the perturbation can be chosen 
arbitrarily small by taking refinement (due to the partial 
hyperbolicity near the periodic point) and 
by taking sufficiently fine refinement in advance
we can guarantee that the 
support of the perturbation is contained in $W$.

Since the periodic 
rectangles converges to $W_{\mathrm{loc}}^u(f_2^i(q))$ 
by taking the forward refinement 
and the property that $R'$ is flat near $y_{i, \pm}$ are not 
affected by taking the forward refinement (note that
forward refinement does not change the 
repelling set but just replace the attracting set), 
we have that 
the periodic rectangles for $\{f_2^i(q)\}$ has flat 
lid boundary 
by taking sufficiently fine refinement.

Next, let us see
how to make the side boundary flat.
The argument is essentially the same.
Let us consider the 
rectangle $K_{f_2^i(q)}$ in the linearized coordinates. 
Then in the linearized coordinates 
$W_{\mathrm{loc}}^s(f_2^i(q))$ is 
a flat plane which coincides with 
the $xy$-plane locally.
For $K_{f_2^i(q)}$, we fix a $C^1$-circle 
$B_i := \partial_s K_{f_2^i(q)}\cap W_{\mathrm{loc}}^s(f_2^i(q))$. 
Now near $f^{-1}_{2}(B_i)$ 
we perform a $C^1$-small perturbation 
so that for $f_3$ (the perturbed map) the 
image of the side boundary near $B_i$ are flat. 
Now we take backward refinement: Notice that taking 
backward refinement does not destroy the flatness of the 
lid boundary. 
On the other hand, by the uniform contraction 
property we know that as the number of the refinement tends
to infinity we have that the rectangle containing 
$K_{f_3^i(q)}$ converges 
to the local stable manifold of $f_3^i(q)$. Thus, at some moment 
all the side boundaries of $\{K_{f_3^i(q)} \}$ 
turn to be flat. 
In particular, the cylinder containing $f_3^i(q)$ is now flat. 

Thus, letting $g = f_3$, we obtain the conclusion.
\end{proof}

\begin{rema}
In the above construction, we take refinements and 
adding perturbations several times. One may wonder 
if the refinement of $\bfR$ with respect to $g$ 
consists of the rectangles we constructed. 
This is true as long as the support of the perturbation 
is contained in the interior of the rectangles. 
More precisely, let us 
consider the refinement $\bfR_{(m,n;f)}$ 
and a perturbation $g$ of $f$. If 
$\supp (g, f) \subset \bfR_{(m, n; f)}$ 
then for $m', n' \geq 0$ we have
\[
[\bfR_{(m, n; f)}]_{(m', n';g)} = 
[\bfR_{(m, n; g)}]_{(m', n';g)} =
\bfR_{(m+m', n+n'; g)}
\]
Thus assuming this holds at each step we can conclude 
the coincidence of two refinements. To be precise, to obtain 
the coincidence we need to confirm this property but 
since it is easy and appearing many times, for simplicity 
we omit this kind of argument. 
\end{rema}
%

\subsection{Linearization of transition rectangles}
Let us discuss the linearization for paths of 
transition rectangles.

\begin{defi}\label{d.tracyc}
Let $\mathbf{R}= \cup C_i$ be a filtrating Markov partition 
of a diffeomorphism $f$
and $\{K_{f^i(q_j)}\}$ ($j=1, 2, i=0,\ldots,\pi_j-1$,
where $\pi_j$ is the period of $q_j$) 
be a cycle of periodic rectangles of 
a periodic point $q_j$. Let $Q$ be a homo/heteroclinic point 
from $q_1$ to $q_2$
and $\{L_{f^k(Q)}\}_{k=1,\ldots, T-1}$ be a path of 
transition rectangles of $Q$, where $T$ is the transition 
time of $Q$.
We assume that $\{K_{f^i(q_j)}\}$ are linearized with 
the local coordinates $\{(U_{j,i}, \phi_{j,i})\}$.
We say that the path $\{L_{f^k(Q)}\}$ 
is \emph{linearized} if the following holds:
\begin{itemize}
\item For each $k=1,\ldots, T-1$, there exists a coordinate
neighborhood $(V_k,\psi_k)$ containing $L_{f^k(Q)}$ 
such that 
$\psi_k(L_{f^k(Q)})$ is a product rectangle.
In the following, we set 
 $(V_0,\psi_0) = (U_{d,1}, \phi_{d,1})$ and
$(V_T,\psi_T) = (U_{a,2}, \phi_{a,2})$, 
where $d$, $a$ are the integers for the departure 
and the arrival rectangles of $Q$ respectively 
(see Definition~\ref{d.deparr}). 
\item For each $k=0,\ldots,T-1$, 
let $J_{f^k(Q)}$ be the 
connected component of 
$L_{f^k(Q)} \cap f^{-1}(L_{f^{k+1}(Q)})$
containing $f^k(Q)$, where we set $L_0$ to be the 
departure rectangle and $L_T$ the arrival rectangle of $Q$
respectively.
Then the map $\psi_{k+1} \circ f \circ \psi^{-1}_{k}$ is 
an affine map preserving the product structure 
$\mathbb{R}^2\times \mathbb{R}$ on 
$\psi_{k}(J_{f^k(Q)})$ in such 
a way that $\mathbb{R}^2$, $\mathbb{R}$
corresponds to $E^{cs}$, $E^u$ directions respectively.
\end{itemize} 
\end{defi}

\begin{rema}
For a path of transition rectangles which is linearized, 
if we take a refinement then the corresponding path
is also linearized.
\end{rema}

\begin{prop}\label{prop:trali2}
Let $\mathbf{R}= \cup C_i$ be a filtrating Markov partition of a diffeomorphism $f$
and $S$ be a circuit of points in $\bfR$. 
Suppose that 
every periodic point of $S$ has a large stable manifold, 
$\bfR(S)$ is a circuit of rectangles for $S$ and 
every cycle of periodic rectangles is linearized.
Then, for any neighborhood $W$ of the orbit of 
homo/heteroclinic points and any $C^1$-neighborhood 
$\cU$ of $f$, there exists $g \in \cU$ such that the 
following holds:
\begin{itemize}
\item $\supp(g,f) \subset W$ and it is disjoint from 
the orbits of the periodic points in $S$.
Especially, the derivatives $Df, Dg$ are the same along 
every periodic orbit of $S$.
\item $f=g$ on $S$. Especially, $S$ is a circuit 
of points for $g$ as well.
\item For every sufficiently large $m$ and $n$,
$\bfR_{(m,n;g)}(S)$ is a circuit of rectangles for $S$ and 
every cycle of periodic rectangles is linearized.
\item Every path of transition rectangles of 
$\bfR_{(m,n;g)}(S)$ is linearized. 
\end{itemize}
\end{prop}
To prove this, we first prove the following: 
\begin{prop}\label{prop:trali}
Let $\mathbf{R}= \cup C_i$ be a 
filtrating Markov partition of a diffeomorphism $f$
and $\{K_{f^i(q_j)}\}$ be a linearized 
cycle of periodic rectangles of 
the periodic point $q_j$ for $j=1, 2$. 
Let $Q$ be a homo/heteroclinic point from $q_1$ to $q_2$
with a path of 
transition rectangles $\{L_{f^k(Q)}\}_{k=1,\ldots, T-1}$.

Suppose that $q_1$ and $q_2$ have large stable manifolds.
Then, for every $C^1$-neighborhood $\cU$ of $f$ 
and every neighborhood
$W$ of $\mathcal{O}(Q)$,
there exists a diffeomorphism $g \in \cU$ 
satisfying the following:
\begin{itemize}
\item The support $\mathrm{supp}(g, f)$ 
is disjoint from $\mathcal{O}(q_i)$ and contained in $W$.
\item $f=g$ on $\cO(q_1) \cup \cO(q_2) \cup \cO(Q) $ 
\item For every sufficiently large $(m, n)$,
the refinement
$\mathbf{R}_{(m,n;g)}$ 
satisfies the following:
\begin{itemize}
\item The corresponding cycle of 
periodic rectangles for $q_j$ in $\mathbf{R}_{(m,n; g)}$ 
is linearized for $j=1, 2$, and 
\item the corresponding 
path of transition rectangles for $Q$ in
$\mathbf{R}_{(m,n; g)}$
is also linearized.
\end{itemize}
\end{itemize}
\end{prop}

\begin{proof}
The proof is very similar to the proof of 
Proposition~\ref{prop:pelipeli}, but it requires 
some extra care.

{\bf Step 1. A perturbation along the transition map.}
First, we give an auxiliary perturbation. 
Let us consider the $(m, n)$-refinement of $\mathbf{R}$. 
Then $Q$ is substituted by $f^{-m}(Q)$ and 
its transition time is $T+m+n$. 
Then consider a perturbation whose support is contained 
in a small neighborhood of $f^{-m}(Q)$ and
$f^{M+n-1}(Q)$ such that $D(f^{m+M+n})(Q)$ preserves
the center stable direction and the strong unstable direction 
of the linearized coordinates. 
Note that, by letting $m$ and $n$ large, 
the $C^1$-size of the perturbation tends to zero, 
thanks to the hyperbolicity near the periodic 
orbits. 
Thus by giving this perturbation 
we may assume that the transition map $Df^T(Q)$ 
preserves the center stable direction and the 
strong unstable direction 
of the linearized coordinates
from the very beginning. Note that this property is 
preserved by taking refinements.

\medskip

{\bf Step 2. Linearizing local dynamics.}
Now let us consider the orbit $\{f^k(Q)\}_{k=0,\ldots,T-1}$.
By applying Franks' lemma along $\{f^k(Q)\}$,
we obtain a diffeomorphism $f_1$ close to $f$
for which we have
linearized coordinates around $\{f_1^k(Q)\}$ for every $k$.
By taking the support sufficiently small, we may assume that
the support is contained in $W$ and  
the perturbation does not disturb the linearization 
property of the periodic rectangles.
We also assume that in the linearized 
coordinates $f_1^k(Q)$ is mapped to the origin, 
the $xy$-plane coincides with the center-stable 
direction and the $z$-axis coincides with the unstable direction.
Note that since $Df^T(Q)$ preserves the strong unstable 
and the center-stable direction in the linearized 
coordinates (see Step 1), these new coordinates are 
compatible with the ones in the periodic rectangles.

{\bf Step 3. Obtaining the product rectangles.}
Then we take refinements.
Due to the largeness of the stable manifolds of $q_j$,
the diameter of transition rectangles goes to zero
as we take refinements. Thus by taking sufficiently fine
refinements the corresponding transition rectangles outside 
the periodic rectangles (the ones before the refinement)
are contained in the domain of linearized coordinates.
In particular, we may assume that they are in $W$.

One thing which is different from the proof of 
Proposition~\ref{prop:pelipeli} is that 
taking refinements increases the number of 
transition rectangles. Notice that these new rectangles are 
contained in the periodic rectangles
(the ones before the refinement).
Thus for these newly produced rectangles we can 
furnish linearized coordinates just by restricting the 
linearized coordinates for the cycle of periodic rectangles and
thus the increase of rectangles does not cause any problem 
for constructing linearized coordinates.

Now let us see how to make the lid and the 
side boundaries of the rectangles flat. 
The argument is almost the 
same as the proof of 
Proposition~\ref{prop:pelipeli}
so we discuss only for the lid boundary.
First, for each $L_{f^k_1(Q)}$, we take the intersection 
of $W_{\mathrm{loc}}^u(f_1^k(Q))$ and the lid
boundary. There are 
two such points. Then we slightly perturb $f_1$ 
so that the lid boundary 
near the intersection points is flat plane. 

Then we take a forward refinement. Notice that taking 
$(0, n)$-refinement increases the number of transition rectangles by $n$. 
On the other hand, the rectangles which are newly created 
are included in the linearized region of the periodic rectangles. 
Thus we know that their lid boundaries are flat. 
For the rest of the rectangles, since taking the forward refinement
shrinks the cylinders in the center stable direction,
we know that up to some sufficiently large 
forward refinement 
we have that all of the transition rectangles 
have flat lid boundaries.

By the same argument, we can make the side boundary of the transition 
rectangles flat as well.
\end{proof}

Now let us discuss the proof of Proposition~\ref{prop:trali2}.
\begin{proof}[Proof of Proposition~\ref{prop:trali2}]
For a circuit of rectangles $\bfR(S)$,
by applying Proposition~\ref{prop:trali}
to each path one by one, 
we linearize all the homo/heteroclinic orbits. 
For that, we need to confirm that 
we can 
apply Proposition~\ref{prop:trali} without destroying the 
linearized coordinates which are already obtained. 
Let us explain how to avoid the interference.

First, in the very beginning we apply the perturbation on 
along each homo/heteroclinic orbit so that the derivative 
of each transition map preserves the center-stable and the 
unstable direction of the linearized coordinates for periodic 
rectangles. 

Then, let us consider a path of transition 
rectangles $\{L_{f^j(Q)}\}_{1\leq j \leq T-1}$
for a homo/heteroclinic point $Q$.
If $\{L_{f^j(Q)}\}$ does not contain any transition rectangles 
which are already linearized, 
then by applying Proposition~\ref{prop:trali} we can 
linearize the homo/heteroclinic orbit $\{ f^j(Q) \}$. 

If not, we consider the refinements. 
Since we assume that every periodic point of the circuit 
has a large 
stable manifold, as we take refinements every rectangle
shrinks to a point. Thus we may assume that the (finitely many) 
rectangles $L'_{Q}, \cdots, L'_{f^T(Q)}$
(these $Q$ and $T$ are the same as the initial ones, we do not take 
the corresponding points and transition times)
are distinct and none of them contain any other 
homo/heteroclinic points of the circuit.

Then, we follow the procedure of Proposition~\ref{prop:trali}. 
By applying Franks' lemma along 
$Q, \cdots, f^T(Q)$, we linearize the local 
dynamics along $Q$. Note that taking refinements increases 
the number of transition rectangles, but by assumption we know 
that for these newly created rectangles we can endow 
linearized coordinates just by taking restrictions. Thus, at this moment, 
the dynamics along the orbit of $Q$ is linearized.
Then we need to make the boundaries of the transition rectangles flat, but 
this can be done in the same way as in the proof of 
Proposition~\ref{prop:trali}.

Repeating this argument, we can obtain the desired coordinates
for every rectangle of the circuit.
\end{proof}

\subsection{On the shape of other rectangles}
In this subsection, we discuss perturbation techniques
which make the shape of rectangles easier to handle.
We begin with a definition.

\begin{defi}\label{d.adap}
Let $\mathbf{R}=\cup C_i$ be 
a filtrating Markov partition of a diffeomorphism $f$.
Let $N$ be a rectangle of $\mathbf{R}$ which is linearized 
(that is, it is either a rectangle in a cycle of periodic 
rectangles 
or a path of transition rectangles which is linearized). 
We say that it is \emph{adapted} 
if every connected component of $f(\mathbf{R}) \cap N$
is a product rectangle in the linearized coordinates.
\end{defi}

\begin{rema}
Let $f \in \mathrm{Diff}^1(M)$, $\mathbf{R}$ be a 
filtrating Markov partition and $\{K_{f^i(q)}\}$
be a cycle of periodic rectangles. If every rectangle in 
$\{K_{f^i(q)}\}$ is adapted then the same holds for 
the corresponding rectangles in the refinements: 
It is obvious for backward refinements. For forward 
refinements, it follows since the rectangles and image
rectangles in the refinements are images of
rectangles and image rectangles, respectively.
 
Similarly, suppose there is a linearized path of transition rectangles
between two cycles of periodic linearized rectangles such that  
all the rectangles involved are adapted. 
Then, the same is true for the corresponding path of 
transition rectangles in the refinements.
\end{rema}

\begin{prop}\label{prop:ado}
Let $\mathbf{R}=\cup C_i$ be a 
filtrating Markov partition of a diffeomorphism $f$
and $S$ be a circuit of points consisting of 
periodic points $\{q_j\}$ and
homo/heteroclinic points 
$\{Q_k\}$ in $\bfR$. Suppose that $\bfR (S)$ is 
a circuit of rectangles for $S$, every $q_j$ has
a large stable manifold, and every cycle of periodic 
rectangles and every path of transition rectangles are linearized.

Then, given a neighborhood $W$ of $S$ 
and a $C^1$-neighborhood $\cU$ of $f$,
there is $g \in \cU$ such that the following holds:
\begin{itemize}
\item $\mathrm{supp}(g, f) \subset W$ and it is 
disjoint from $\{\cO(q_j)\}$ and $\{\cO(Q_k)\}$.
\item For every sufficiently large $(m, n)$, every 
rectangle in $\mathbf{R}_{(m,n)}(S)$ is  
adapted (and linearized).
\end{itemize}
\end{prop}

\begin{proof}
Let us see the perturbations which make rectangles 
in a cycle of periodic rectangles and a path of
transition rectangles adapted. Then applying these  
perturbations one by one, we obtain the conclusion.

First, note that due to the 
largeness of the stable manifolds of periodic 
points in $S$,
every sufficiently fine refinements $\mathbf{R}'$ 
satisfies  
$\mathbf{R}'(S) \subset W$. 
Thus we may assume that this holds from the very 
beginning.

Let us see how to make the cycles of periodic rectangles 
adapted.
Let $\{K_{f^i(q)}\}$ be a cycle of periodic rectangles 
for $q$.
We consider $f(\mathbf{R}) \cap K_{f^i(q)}$. 
Recall that as we take backward refinements, 
the height of the rectangle which 
$f^i(q)$ belongs to decreases and tends to zero,
while $W^{s}_{\mathrm{loc}}(q)$ does not 
change. 
By assumption, the image rectangles 
$f(\bfR) \cap K_{f^i(q)}$ containing a point of $S$ 
are product rectangles. 
For the other image rectangles, they automatically have flat 
lid boundaries but possibly with non-flat side 
boundary.
Note that, if we take enough forward refinement, 
each rectangle in $f(\mathbf{R}') \cap K'_{f^i(q)}$
are almost vertical thanks to the partial hyperbolicity 
near $\cO(q)$. 

Then by giving a $C^1$-small perturbation,
we can construct a diffeomorphism $f_1$ close to $f$
such that each connected component of 
$f_1(\mathbf{R}') \cap K'_{f_1^i(q)}$ has 
flat boundaries
near $W_{\mathrm{loc}}^s(f_1^i(q))$. 
Since these perturbations are only for the 
boundary of rectangles 
which are not periodic rectangles and transition rectangles, 
we may assume that the support is contained in $W$ and 
disjoint from $S$.
Also, we may assume that 
these perturbations do not disturb the 
linearization property on the periodic rectangles and 
transition rectangles. 

Now, by taking backward refinement,
we have the product property 
for $f_1(\mathbf{R}^{\prime\prime}) \cap K^{\prime\prime}_{f_1^i(q)}$.
Thus we can obtain the adaptedness for a cycle 
of periodic rectangles. 
By repeating this perturbation, we may
assume that every cycle of periodic rectangles is adapted. 

The proof for the transition rectangles can be done 
similarly. 
The only thing we need to pay extra attention 
is that by taking refinements
the number of transition rectangles increases. 
Meanwhile, once we have the adaptedness for periodic
rectangles, this does not bring any problem 
for the following reasons:
\begin{itemize}
\item Taking the $(0, n)$-refinement increases the number 
of transition rectangles by $n$. It adds $n$ rectangles
in the cycle containing the arrival rectangle. 
However, the first $n$ transition
 rectangles are images of periodic rectangles which are adapted. 
 Thus the number of rectangles which are not adapted 
 are the same.
 \item Taking $(m, 0)$-refinement 
 increases the number 
of transition rectangles by $m$. It adds $m$ 
rectangles in the head. But the first $m$ transition
 rectangles are contained in periodic rectangles which are adapted. 
 Thus the number of rectangles which are not adapted 
 are the same.
\end{itemize}
In short, these newly added rectangles are automatically 
adapted. Hence, in order to obtain the adaptedness 
for transitions we only need to repeat the argument 
for finitely many rectangles. 
%
%
%
%
%
%
\end{proof}
%

Now we are ready to state the definition of affine Markov 
partitions.
\begin{defi}
\label{d.affi}
Let $f \in \mathrm{Diff}^1(M)$, $\mathbf{R}$
be a filtrating Markov partition and $S$ be a circuit 
of points contained in $\mathbf{R}$. We say that 
$\mathbf{R}(S)$ is an affine Markov partitions if we 
have the following:
\begin{itemize}
\item $\mathbf{R}(S)$ is a circuit of rectangles for 
$S$ (see Definition~\ref{d.cirrec}).
\item Every cycle of periodic rectangles in $\mathbf{R}(S)$
is linearized and adapted (see Definition~\ref{d.lincyc} and
\ref{d.adap}).
\item Every path of transition rectangles in $\mathbf{R}(S)$
is linearized and adapted (see Definition~\ref{d.tracyc} and
\ref{d.adap}).
\end{itemize}
\end{defi}
Note that the arguments in this section,
more precisely, 
Proposition~\ref{prop:rec_choix}, \ref{prop:pelipeli}, 
\ref{prop:trali2} and \ref{prop:ado}
conclude Theorem~\ref{t.affine}.
Indeed, 
\begin{itemize}
\item Proposition~\ref{prop:rec_choix}
guarantees that if we take sufficiently fine refinements 
then $\bfR(S)$ is a circuit of rectangles.
\item By applying Proposition~\ref{prop:pelipeli}
to each cycle of periodic rectangles, by an arbitrarily small 
perturbation we can linearize the cycle, up to some refinements.
\item By applying Proposition~\ref{prop:trali2}
by an arbitrarily small 
perturbation we can linearize all the paths, up to some refinements.
\item Finally, by Proposition~\ref{prop:ado}, by 
an arbitrarily small perturbation we obtain the adaptedness 
for every rectangles, up to some refinements.
\end{itemize}
Thus we obtain the conclusion.

\subsection{Robustness of filtrating Markov partitions}
In this subsection, we clarify the definition of the  
robustness of a filtrating Markov partitions which 
we proposed in Introduction. Also, we discuss the relation 
between the perturbation technique we discussed 
in this section and the robustness. 

Let us begin with the definition.
\begin{defi}\label{d.robu}
Let $\mathbf{R}$ be a filtrating Markov partition of 
$f \in \mathrm{Diff}^1(M)$. We say that $\mathbf{R}$
is \emph{$\alpha$-robust} if there is a cone field
$\mathcal{C}$ over $\mathbf{R}$ which satisfies 
the definition of filtrating Markov partition (see Definition~\ref{d.fmp}) and for any $C^1$-diffeomorphism
$g$ which is $\alpha$-close to $f$ the cone field 
$\mathcal{C}$ is strictly-invariant and unstable.
\end{defi}

Note that, in the definition of a filtrating Markov partition
(see Definition\ref{d.fmp}), 
except the last condition every condition refers some properties
about the behavior of $f$ on $\mathbf{R}$. Thus 
if $\mathbf{R}$ is $\alpha$-robust, then for every 
$g$ which is $\alpha$-close to $f$ and whose support 
is contained in the interior of $\mathbf{R}$, we know 
that $\mathbf{R}$ is a filtrating Markov partition for 
$g$ with the same coordinates and the same cone field.

The $\alpha$-robustness gives a sufficient condition for 
the persistence of a filtrating Markov partition, but a priori
it may be that $\mathbf{R}$ persists under a perturbation 
whose size is larger than $\alpha$.

\begin{rema}\label{r.concon}
Let us discuss the robustness of refinements of 
filtrating Markov partitions. 
\begin{itemize}
\item Suppose that $\mathbf{R}$ is $\alpha$-robust.
Then recall that for the refinement $\mathbf{R}_{(m, n)}$
the cone field for $\mathbf{R}$ also satisfies the 
assumption of Definition~\ref{d.fmp}. Thus we know that 
$\mathbf{R}_{(m, n)}$ is also $\alpha$-robust
with the same cone field (see Proposition~2.10 in \cite{BS2}).
\item Now consider a filtrating Markov partition
$\mathbf{R}$ which is $\alpha$-robust and suppose that 
for a refinement $\mathbf{R}_{(m, n)}$, the rectangles
$\mathbf{R}_{(m, n)}(S)$ are affine for some circuit $S$.
Then, each rectangle in $\mathbf{R}_{(m, n)}(S)$
has linearized coordinates. 
In general, we do not know if
the restriction of the cone field for 
$\mathbf{R}_{(m, n)}(S)$ satisfies
the condition of Definition~\ref{d.fmp} with respect 
to the linearized coordinates. However, by the 
partial hyperbolicity on the circuit, we know that the 
cone field must contain $z$-direction and does not contain
$x, y$-direction over the points of $S$. Then, by the 
continuity of the cone field, we know that if the linearized 
coordinates are defined in a sufficiently small neighborhood 
of $S$, then we have the compatibility between the 
linearized coordinates and the cone field.
Thus, by taking sufficiently fine refinement
we know that the restriction of the cone field of 
$\mathbf{R}$ to $\mathbf{R}_{(m, n)}(S)$ gives 
a vertical, strictly invariant unstable cone field
with respect to the linearizing coordinates.  
\end{itemize}
\end{rema}

\subsection{Realizing 2 dimensional perturbation in dimension 3}
In this subsection, we consider the following perturbation result.
\begin{prop}\label{prop:3dim}
Given a filtrating Markov partition $\mathbf{R} = \cup C_i$, 
suppose that 
there are rectangles $N_i$ ($i=1, 2$) which are linearized. 
Let $\phi_i$ be the linearization coordinates,
$\phi_i(N_i) = D_i \times I_i$ and assume that 
$N_1 \cap f^{-1}(N_2) \neq \emptyset$. Let
$F:D_1 \to D_{2}$ be the corresponding 
two dimensional maps on $D_1 \times J$, 
a connected component of 
$\phi_1^{-1}\left( N_1 \cap f^{-1}(N_2) \right)$.

Suppose that we have a $C^1$-diffeomorphism
$G:D_1\to D_2$ such that
\begin{itemize}
\item $G$ coincides with $F$ near the boundary of $D_1$,
\item the $C^0$-distance between $F$ and $G$ is less than $\varepsilon_0$,
\item the $C^1$-distance between $F$ and $G$ is less than $\varepsilon_1$.
\end{itemize}
Then, there exists a $C^1$-diffeomorphism $g$ which 
is ($\varepsilon_1 + K \varepsilon_0$)-$C^1$-close to $f$ 
(where $K$ is some constant which depends only on the 
choice of rectangles) 
such that the following holds:
\begin{itemize}
\item The support $\supp (g, f)$ is contained in the interior of
$ N_1$.
\item $g$ keeps the product structure for $\phi_i$ and 
the two dimensional map over $D_1 \times J$ is given by $G$.
\end{itemize}
\end{prop}
\begin{proof}
The construction of $g$ can be done by standard arguments 
involving the
partition of unity and the closeness of $G$ to $F$ in 
the $C^1$-distance. We just give a sketch of the proof. 

Assume that on $D_1 \times J$ the map $\phi_2 \circ f \circ \phi_1 ^{-1}$ is 
given by
\[
(x, y, z) \mapsto (F(x,y), \lambda z)
\]
in some neighborhood of $D_1 \times J$. 
Then, we choose a $C^1$-function $\rho$ defined in the 
interval $J'$ which contains $J$ in the interior and 
satisfying the following: 
\begin{itemize}
\item $0 \leq \rho (z) \leq 1$.
\item $\rho(z) \equiv 1$ on $J$.
\item $\rho(z)\equiv 0$ near the endpoints of $J'$.
\end{itemize}
Then, given $G$ consider the following map:
\[
(x, y, z) \mapsto ((1-\rho(z))F(x,y) + \rho (z)G(x,y),  \lambda z).
\]
Considering the fact that $F \equiv G$ near the 
boundary of $D_1$,
This map is equal to $\phi_i \circ F \circ \phi_i ^{-1}$
near the boundary of $D_1 \times J'$, thus extends to outside 
$D_1 \times J'$ so that it coincides with the unperturbed map. If $G$ is sufficiently $C^1$-close to $F$, then 
one can check that this defines a diffeomorphism on each
slice by $xy$-plane.
The surjectivity of the map is the consequence of standard 
algebraic topological argument.
The injectivity follows if we 
choose $F$ sufficiently close to $G$. 

Now we measure the $C^1$-size of this perturbation. 
By a direct calculation, the difference of the derivatives of
$f$ and $g$ in the local coordinates is given by
\[
\begin{pmatrix}
\rho (z)D_{x, y}(G(x, y) -F(x,y)) & 0 \\
\rho' (z)(G(x, y) -F(x,y)) & 0
\end{pmatrix},
\]
where $D_{x, y}$ denotes the Jacobi matrix with respect to 
$x$ and $y$.
This calculation shows that the $C^1$-distance is given 
by two terms $|\rho (z)D_{x,y}(G(x, y) -F(x,y))|$ and 
 $|\rho' (z)(G(x, y) -F(x,y))|$. The supremum norm of 
the first one is proportional 
to the $C^1$-distance between $F$ and $G$, and the 
second one is to $|\rho' (z)|$ times the $C^0$-distance
between $F$ and $G$. 
Thus, the $C^1$-distance of the perturbation itself is 
given by the form $\varepsilon + K \delta$, where $K$ 
is determined by $\rho'$, 
which depends only on the shape of 
$J'$.
\end{proof}


\section{Markov IFS}
\label{s.ifs}

To prove Theorem~\ref{theo:expu},
we investigate two dimensional iterated function systems (IFSs) 
where the iteration is chosen in a Markovian way. 
In this section, we give the precise definition of it 
and discuss their elementary properties. 
In Section~\ref{ss.corres} we discuss the relation between 
Markov IFSs and affine Markov partitions. This enables us to
use Markov IFSs for the investigation of the bifurcation 
of filtrating Markov partitions.

\subsection{Definition}

By a \emph{disc} we mean
a subset of $\mathbb{R}^2$ which 
is $C^1$-diffeomorphic to 
a two dimensional round disc. 
Let $\cD=\coprod_{i=1}^m D_i$,
the disjoint union of a finite set of discs $D_i$.
We denote the boundary of $D_i$ by $\partial D_i$,
put $\mathrm{Int}(D_i) = D_i \setminus\partial D_i$ and call it \emph{(geometric) interior}.

\begin{rema}
We introduce the topology induced from $\mathbb{R}^2$ for each $D_i$. 
Accordingly, each $D_i$ itself is an open set
and the topological boundary of $D_i$ is empty. 
This seemingly strange topology will be convenient for instance when we define 
the notion of relatively repelling regions, see Section~\ref{ss.rrr}. 
\end{rema}

\begin{figure}[h]
\begin{center}
\includegraphics[width=8cm]{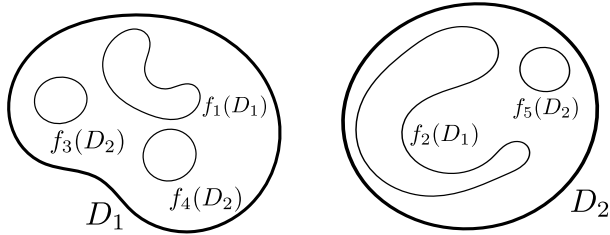}
\caption{An example of Markov IFS. It consists of two discs 
$\cD = D_1 \coprod D_2$ and five diffeomorphisms $F = \{f_1,\ldots, f_5\}$
on their images.
The non-trivial restriction on Markov IFS is that the images of discs have 
empty overlaps.}
\label{f.IFS}
\end{center}
\end{figure}

The following is the formal definition of a Markov IFS, 
see also Figure~\ref{f.IFS}.

\begin{defi}
A \emph{Markov IFS} on $\cD=\coprod_{i=1 }^m D_i$ is a family of 
finitely many local diffeomorphisms (where 
a local diffeomorphism means a diffeomorphism on its image)
$F = \{f_j\}_{1 \leq j \leq k}$ such that the followings hold:
\begin{itemize}
 \item For every integer $j \in [1, k]$ there are
 integers $\mathrm{dom}(j) \in [1,m]$ 
 and $\mathrm{im}(j) \in [1, m]$
 such that $D_{\mathrm{dom}(j)}$ 
 is the domain the definition of $f_j$ (called the domain disc of $f_j$) and 
 $\mathrm{im}(f_j) :=f_j(D_{\mathrm{dom}(j)})$  
 (called the image of $f_j$) is contained in 
 $\mathrm{Int}(D_{\mathrm{im}(j)})$ ($D_{\mathrm{im}(j)}$ is called the target disc of $f_j$). 
 \item
 The images $\{\mathrm{im}(f_i)\}_{i \in [1,k]}$ are pairwise disjoint.
\end{itemize}
\end{defi}
\begin{rema}\label{r.inverse} 
We put $F(\mathcal{D}) := \coprod_{1\leq j\leq k} \mathrm{im}(f_j)$.
The collection $\{f_j^{-1}\}_{1\leq j\leq k}$ defines a uniquely 
defined inverse map from $F(\mathcal{D})$ to $\mathcal{D}$. 
We denote it by $F^{-1}$.
\end{rema}

%

\subsection{Periodic points and homo/heteroclinic orbits}
\label{ss.peri}
In this subsection, we introduce several basic definitions related 
to Markov IFSs.

\subsubsection{Periodic points}\label{sss.per}
Let $(\cD=\coprod_{i=1 }^m D_i, F = \{f_j\}_{1 \leq j \leq k})$
be a Markov IFS. 
We consider words whose letters are in $\cI = [1, k]$.
We say that a non-empty 
word $ \omega=j_1\cdots j_n$
is \emph{admissible} if $f_{j_m}(D_{\mathrm{dom}(j_m)}) \subset
D_{\mathrm{dom}(j_{m+1})}$
holds for every $m = 1, \ldots, n-1$.
For an admissible word $\omega$, 
we put $F_\omega:= f_{j_n}\circ\cdots\circ f_{j_1}$.
We say that a point $p \in \cD$ 
is \emph{periodic} if $F_\omega(p) = p$ holds for some 
admissible $\omega$. 
The period of $p$ is the least length of non-empty
word $\omega$ for which 
$f_{\omega}(p) =p$ holds. 
As a straightforward consequence of Remark~\ref{r.inverse}, we have the following: 
\begin{rema}If $p$ is a periodic point of period $n$, 
there is a unique word $\omega(p)$ of length $n$ 
such that $p$ is a fixed point of $F_{\omega(p)}$.  
We call $\omega(p)$ the \emph{itinerary} of $p$.
If $p$ is a fixed point of $F_{\omega'}$ where
$\omega'$ is another word, then $\omega'$ is a concatenation 
of several copies of $\omega(p)$.  
\end{rema}

The \emph{periodic orbit} of the periodic point $p$,
denoted by $\mathrm{orb}(p)$, 
is the set of points $p_i=f_{j_i}\circ\cdots\circ f_{j_1}(p)$, 
$i\in\{1,\dots, \pi(p)\}$, 
where $\pi(p)$ is the period of $p$ and we put 
$\omega(p) = j_1\cdots j_{\pi(p)}$. We set $p_0=p$ and $f_{j_{\pi +1}} = f_{j_1}$.  
In the following, by abuse of notation we write
$F^i(p)$ in the sense of $p_i$.
\begin{rema} Since $F^{-1}$ is a 
well-defined map, we have 
$p_i=F^{-\pi(p)+i}(p)$ for $1 \leq i \leq \pi(p)$.
\end{rema}

For $x \in \cD$, we denote
the disc of $\cD$ which contains $x$ by $D_x$. 
A periodic point $p\in \cD$ is 
called a \emph{hyperbolic periodic point} 
if it is a hyperbolic fixed
point of $F_{\omega(p)}: D_p \to D_p$ 
and its \emph{$s$-index} is the dimension of its stable 
manifold.
Suppose that $DF_{\omega(p)}|{_{T_pD_p}}$
has two eigenvalues $0<\lambda_1<1, \lambda_1 < \lambda_2$
(we allow the case $\lambda_2 \leq 1$).
Then the \emph{local strong stable manifold} of $p$, 
denoted by $W^{ss}_{\mathrm{loc}}(p)$, 
is the strong stable manifold of $p$ of $F_{\omega(p)}$ in $D_p$
tangent to the eigenspace of $\lambda_1$ at $p$. 
A periodic orbit $\mathrm{orb}(p)$
is called \emph{separated} if $\{D_{p_i}\}_{i=0,\ldots,\pi(p)-1}$ are all distinct. 
We say that $p$ is separated if $\mathrm{orb}(p)$ is.
For a separated periodic point $p$, by $F_p$ we denote the map 
$\coprod D_{p_i} \to \coprod D_{p_i}$ defined by
$F_{p}|_{D_{p_i}} = f_{j_{i+1}}$
(recall that we put $f_{j_{\pi(p)+1}} = f_{j_1}$). 
Let $p'$ be another periodic point whose orbit is not 
equal to that of $p$. We say that $p$ and $p'$ are 
\emph{mutually separated} if there is no disc in $\cD$
which contains points of $\mathrm{orb}(p)$ and 
$\mathrm{orb}(p')$.

We say that a periodic point $p\in D_p$ 
has a \emph{large stable manifold} 
if the whole disc $D_p$ is contained in the local
stable set $W^{s}_{\mathrm{loc}}(p)$
for $F_{\omega(p)}$, where we put 
$W^{s}_{\mathrm{loc}}(p) := \{ y \in D_p \mid 
(F_{\omega(p)})^{n}(y) \to p \quad (n \to \infty)\}$.

\subsubsection{Homo/heteroclinic points}
\label{sss.hh}
Let $p$ be a separated periodic point.
A point $P \in W^{s}_{\mathrm{loc}}(p_i) \setminus 
f_{j_{i-1}}(D_{p_{i-1}})$ 
is called a 
\emph{$u$-homoclinic point of $p$} if 
there is an integer $k \geq 0$ such that 
$F^{-k}(P)=p_l$ holds for some $l$. 
If $p$ has the strong stable manifold $W^{ss}_{\mathrm{loc}}(p)$, 
then $P$ is a \emph{$u$-homoclinic point of $p$} if 
$P \in W^{ss}_{\mathrm{loc}}(p_i)$ and it is a $u$-homoclinic point of $p$, too. 
For a $u$-homoclinic point $P$, there exists a word $\omega$ such that 
$(F_{\omega})^{-1}(P) = p_l$ holds. One can check that
there is the unique shortest word among 
such words. We denote it by $\omega(P)$ and call it the \emph{itinerary} of $P$. 
For $u$-(strong) homoclinic points, the backward orbit
$F^{-i}(P)$ makes sense for $i \geq 0$. 
Also, for $i\geq 0$ we define $F^i(P) := (F_p)^i(P)$.
We put $\mathrm{orb}(P) := \{F^{i}(P) \}$ and 
call it the \emph{homoclinic orbit of $P$}.
Given two periodic points $p_1$ and $p_2$ having different orbits, 
we also define the notion of \emph{$u$-heteroclinic points} 
in a similar way.
The notion of the itinerary and the heteroclinic orbit are defined similarly.

Let $P\in W^s_{\mathrm{loc}}(p)$ be a $u$-homo/heteroclinic 
point of a periodic point $p$.
We say that 
$P$ is \emph{$p$-free} if the following holds: Let $k$ be the 
least positive integer such that $F^{-k}(P) \in \mathrm{orb}(p)$. Then 
for every $i=1,\ldots,k-1$, we have 
$F^{-i}(P) \not \in \cup_{0\leq i \leq \pi(p)-1} D_{p_i}$

\subsubsection{Perturbations of IFSs}
Given a Markov IFS $(\cD, F = \{f_i\}_{i=1,\ldots,k})$, 
consider 
a family of $C^1$-diffeomorphisms 
$\tilde{F} = \{\tilde{f}_i\}_{i=1,\ldots,k}$
such that the following holds:
\begin{itemize} 
\item Near the boundary of the domain disc, $\tilde{f}_i \equiv f_i$ for every $i$.
\item Each $\tilde{f}_i$ is a $C^1$-diffeomorphism with the same image as $f_i$
for $1\leq i\leq k$. 
\end{itemize}
Then, $(\cD, \tilde{F})$ is also a Markov IFS with the 
same set of admissible words. We call it a
\emph{perturbation} of $(\cD, F)$.
For each $i$, 
we define $\mathrm{supp}(\tilde{f}_{i})$
to be the closure of the set 
$\{x \in D_{\mathrm{dom}(i)} \mid \tilde{f}_{i}(x)\neq f_i(x)\}$ and call it 
the \emph{support} of $\tilde{f}_{i}$. We put 
$\mathrm{supp}(\tilde{F}) := \cup_{i=1}^k\mathrm{supp}(\tilde{f}_{i})$. 

Let $G:=\{g_i\}$ be a perturbation of $F$ and $p$ be a 
periodic point of $F$. $G$ is called a \emph{perturbation along the 
orbit of $p$} if and only if
they only differs for $\{f_{j_i}\}_{i=1,\ldots,\pi(p)}$, 
where $j_i$ is a letter appearing $\omega(p)$.

\subsection{Refinement of a Markov IFS}
Let $(\cD=\coprod_i D_i, F=\{f_j\})$ be a Markov IFS. For $n>0$, 
consider the disjoint union of the 
images of the discs 
$$F^n(\cD):=\coprod_{\substack{|\omega|=n,\\ \omega: \mbox{\tiny{admissible}}}} F_\omega(D_{\omega}),$$
where $|\omega|$ denotes the length of 
the word $\omega$ and $D_{\omega}$ denotes the domain of $F_{\omega}$.  
Now, consider the collection of local 
diffeomorphisms 
\[\wedge_n F=\{f_i|_{F_\omega(D_{\omega})} \mid |\omega|=n,\,\, \omega: \mbox{admissible},\,\,
F_\omega(D_{\omega}) \subset D_{\mathrm{dom}(i)}\}.\] 
Then the pair 
$(F^n(\cD), \wedge_n F)$ defines a Markov IFS.
We call it an 
\emph{$n$-refinement of $(\cD,F)$}. 

\begin{rema}  
\begin{itemize}
\item $(\wedge_n F)^{-1}$ is the restriction of $F^{-1}$ to $F^{n}(\cD)$. 
\item A point $x \in F^n(\cD)$ is a periodic point of $\wedge_n F$ if and only if it is periodic for $F$. 
In such a case, the periods of $x$ for $F$ and 
$\wedge_n F$ are same.
\item  A periodic point $x \in \mathcal{D}$ 
has a large stable manifold for $\wedge_n F$ 
if and only if it has large stable manifold for $F$. 
The equivalence follows by noticing that 
the disc in $F^n{\cD}$ which contains $x$ is the image of the disc in $\mathcal{D}$ 
which contains $F^{-n}(x)$.
\item Suppose we have a $u$-homoclinic point $P$ of 
a periodic point $p$ with a large stable manifold. 
Then for the refinement $(\wedge_1 F)$
we take the point $f_{j_P}(P)$ (where $j_P$ is the 
letter of $\omega(p)$ such that $x\in \mathrm{dom}(f_{j_P})$ holds) 
and call it the
\emph{$u$-homoclinic point corresponding to $x$} 
for $(\wedge_1 F)$. Inductively we define the 
corresponding homoclinic point for $(\wedge_n F)$.
We define the same notion for heteroclinic points.
\end{itemize}
\end{rema}

\subsection{Flexible periodic points}
\label{ss:flex}

In \cite{BS1}, we defined the notion of an
\emph{$\varepsilon$-flexible periodic point}.
It is a periodic point of a diffeomorphism whose
invariant manifold is so flexible that 
its configuration in an prescribed
fundamental domain can be deformed into 
arbitrarily chosen shape
by an $\varepsilon$-small perturbation. 
Let us recall precise definition of it
(see \cite{BS1} for further information).

\begin{defi}\label{d.fleco}
Let $(A_i)_{i = 0, \ldots ,n-1}$ be a
two dimensional linear cocycle,
that is, let $A_i \in \mathrm{GL}(2, \mathbb{R})$ for every $i$. We say that $(A_i)$ is an
\emph{$\varepsilon$-flexible cocycle} if 
there exists a continuous path of 
linear cocycles 
$\mathcal{A}_t = (A_{i, t})_{t \in [-1, 1]}$
such that the 
following holds:
\begin{itemize}
\item $\mathrm{diam}(\mathcal{A}_t) < \varepsilon$, that is, 
for every $i$ we have 
$\max_{-1\leq s<t\leq 1}\|A_{i,s} - A_{i,t}\| < \varepsilon$.
\item $A_{i, 0} = A_i$ for every $i$.
\item For every $t \in (-1, 1)$, 
the product $A_{t} := A_{n-1, t} \cdots A_{0, t}$
has two distinct positive contracting eigenvalues.
\item $A_{-1}$ is a contracting homothety.
\item Let $\lambda_t$ be the smallest eigenvalue
of $A_t$. Then 
$\max_{-1\leq t \leq 1} \lambda_t <1$.
\item $A_{1}$ has an eigenvalue equal to $1$.
\end{itemize}
\end{defi}

In regard to Markov IFSs, 
a periodic point $p$ of
a Markov IFS $(\mathcal{D}, F=\{f_i\})$ 
is called $\varepsilon$-flexible if the 
linear cocycle of linear maps 
$(Df_{j_i})_{i=1,\ldots, \pi(p)}$ 
between tangent spaces over the orbit of $p$ 
is $\varepsilon$-flexible, where we put $\omega(p) = j_1\cdots j_{\pi(p)}$. 

\subsection{Affine circuits and Markov IFSs}
\label{ss.corres}
Let us clarify the relation between Markov IFSs and 
affine Markov partitions. 
\begin{defi}
Suppose that we have a filtrating Markov partition 
$\mathbf{R}$ with a circuit of points $S$ such 
that $\mathbf{R}(S)$ is an affine Markov partition. 
Assume that:
\begin{itemize}
\item Each (linearized) rectangle 
has the form $D_i \times I_i$ in the linearized coordinates.
\item For each connected component of $C_i \cap f^{-1}(C_j)$ containing a point of $S$, 
where $C_i$, $C_j$ are rectangles of $\mathbf{R}(S)$,
the dynamics in two dimensional direction 
is given by a map $F_{(i, j), k}:D_i \to D_j$.
\item For each $C_i$, the image rectangles $C_i \cap f(\mathbf{R})$ 
which does not contain any point of $S$
has the form $E_{i, l} \times I_i$ in the linearized coordinates.
\end{itemize}
Then, we can define a Markov IFS as follows:
\begin{itemize}
\item The set of discs are $\{D_i\} \cup \{E_{i,l}\}$.
\item The set of maps are $\{F_{(i, j), k}\} \cup \{\mathrm{id}_{i,l}\}$,
where $\mathrm{id}_{i,l}:E_{i,l} \to D_i$ is the restriction of the identity map.
\end{itemize}
This Markov IFS is called the \emph{corresponding Markov
IFS} of $\mathbf{R}(S)$ and we denote it
by $\cM(\mathbf{R}(S))$.
\end{defi}
In the proof of Theorem~\ref{theo:expu}, the information 
of $f$ on the adapted rectangles are not important. The only 
information we need is the shape of the image of the rectangles. 
Thus for the maps on $\{E_{i,l}\}$ we put the identity maps.

Also, we can define the periodic 
points $\{q_i\}$ and homo/heteroclinic points $\{Q_j\}$ for $\cM(\mathbf{R}(S))$.
For $\{q_i\}$ we just take the projections of them in
the above Markov IFS. For homo/heteroclinic
points, instead of dealing with $\{Q_j\}$
we consider the projections of $\{f^{T_j}(Q_j)\}$,
where $T_j$ is the transition time of $Q_j$. Then 
they gives $u$-homo/heteroclinic points in the Markov IFS.

Note that 
\begin{itemize}
\item If $\{q_i\}$ are $\varepsilon$-flexible 
or have large stable manifolds, then the same
property holds 
for the projected periodic points.
\item If $\mathbf{R}(S)$ is a circuit of rectangles, 
then for the associated Markov IFS, every periodic 
orbit is separated, every pair of periodic orbits 
are mutually separated and every homo/heteroclinic orbits 
are free from every periodic orbit (see Section~\ref{sss.per}
for the definitions).
Furthermore, if $\mathbf{R}$ is generating, then 
each pair of homo/heteroclinic orbits has different 
itineraries (see Lemma~\ref{l.shadow}).

\item The operation of taking forward refinements is functorial 
with respect to the operation $\cM$. 
More precisely, if $\mathbf{R}(S)$
is an affine circuit, then we have
\[
\cM(\mathbf{R}_{(0, n)}(S)) = \wedge_n\cM(\mathbf{R}(S)).
\] 
\end{itemize}

\subsection{Relatively repelling regions}
\label{ss.rrr}
We are interested in constructing a repelling region 
by giving small perturbation to a Markov IFS. 
Since what we deal with is not a single diffeomorphism 
but an IFS, the formulation of the notion of repelling/attracting sets requires extra care.
In the following subsections we will discuss their definitions.

\begin{figure}[h]
\begin{center}
\includegraphics[width=8cm]{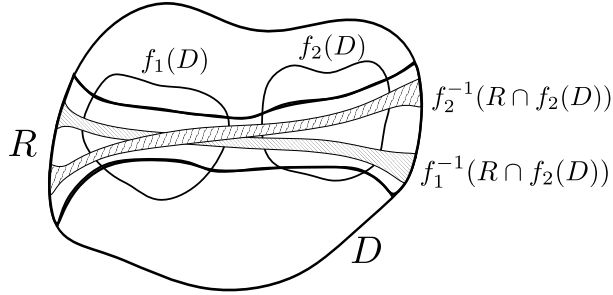}
\caption{An example of a relatively repelling region. The Markov IFS 
consists of one disc $D$ and two diffeomorphisms $F = \{f_1,f_2\}$
on their images. The relative repelling region $R$ 
has attracting property for $f_1^{-1}$ restricted to 
$f_1(D) \cap R$ and for $f_1^{-1}$ to $f_2(D) \cap R$. }
\label{f.relarep}
\end{center}
\end{figure}

\begin{defi}\label{d.relare}
Let $(\cD, F)$ be a Markov IFS.
We say that a compact set 
$R \subset \mathcal{D}$ is a \emph{relatively repelling region} if 
we have
$$F^{-1}(R \cap F(\cD))\subset \mathrm{Int}_{t}(R),$$
where $\mathrm{Int}_{t}(R)$ denotes the 
topological interior of $R$ 
with respect to the topology of $\cD$, see Figure~\ref{f.relarep}. 
\end{defi}

The definition may look strange, but
since $F^{-1}$ is a well-defined map, it is natural to define the repelling property 
as an attracting property for $F^{-1}$.  

While we do not use the following lemmas in this paper, to have 
better understanding of the notion of relatively repelling regions, 
let us prove the followings.

\begin{lemm}\label{l.repellingiterate} 
If $R \subset \cD$ is a relatively repelling region, then 
$R_n := F^{-n}(R \cap F^n(\cD)) \subset \cD $ is a
relatively repelling region for $(\cD, F)$.   
\end{lemm}

\begin{proof}
We prove the case $n =1$, that is, 
$F^{-1}(R_1 \cap F(\cD) ) \subset \mathrm{Int}_t( R_1)$. 
The general case follows by induction: 
Notice that 
\begin{align*}
R_{n+1} &= F^{-1}(F^{-n}(R \cap F^{n+1}(\cD)  ))\\
&= F^{-1}(F^{-n}(R\cap F^{n}(\cD) \cap (F^{n+1}(\cD)))=F^{-1}(R_n\cap F(\cD)).
\end{align*}
 
Consider a point $x \in F^{-1}(R_1 \cap F(\cD) )$.
By definition, there is a point $y \in \cD$ such that $x =F^{-2}(y)$.
Then consider the point $F^{-1}(y)$, which belongs to 
$R_1 \cap F(\cD)$.  
We know $F^{-1}(y) \in F(\cD)$. 
By the fact that $R$ is a relatively
repelling region and $F^{-1}(y) \in R_1$, 
we know that $F^{-1}(y)  \in \mathrm{Int}_t(R)$.
Then, by taking the inverse image, we have $x =F^{-2}(y) \in \mathrm{Int}_t(F^{-1}(R \cap F(\cD))) = \mathrm{Int}_t(R_1)$.
\end{proof}

\begin{lemm}
\label{l.gouge}
 If $R \subset F^n(\cD)$ is a relatively repelling region for the $n$-refinement 
 $(F^n(\cD),\wedge_n F )$ of $(\cD,F)$ then  
 $R\cup \overline{(\cD\setminus F^n(\cD))}$ is a relatively repelling region for 
 $(\cD, F)$, where $\overline{X}$ denotes the closure of $X$. 
\end{lemm}

\begin{proof}
We prove the case $n=1$. Then the general case follows by induction. 
Consider the set $\tilde{R}:=R\cup \overline{(\cD\setminus F(\cD))}$ and take a point 
$x \in \tilde{R}$. If $x \not\in F(\cD)$, then there is nothing 
we need to prove. Suppose $x \in F(\cD)$. Then, we have 
either $x \in F^2(\cD)$ or not. For the first case, by the definition of 
relatively repelling region, we know that $F^{-1}(x) \in 
\mathrm{Int}_t(R) \subset \mathrm{Int}_t(\tilde{R})$. 
For the latter case $F^{-1}(x) \in 
\cD \setminus F(\cD) \subset \mathrm{Int}_t(\tilde{R})$. 
This shows relatively repelling property of $\tilde{R}$. 
\end{proof}

%
%
%

\subsection{Contracting invariant curves}
In this subsection, we formulate the notion of \emph{(normally 
repelling/attracting) contracting invariant curves} and 
discuss some related notions. 

In this article, by a \emph{curve}
we mean an image of a $C^1$-embedding of an interval to $\cD$ satisfying the following two 
conditions:
\begin{itemize}
\item The image intersects with the boundary of the disc transversely at the endpoints.
\item Except the endpoints, the image has no intersection with the boundary of the disc. 
\end{itemize}
By a family of curves we mean a union of finitely many $C^1$-curves
(some of them may have non-empty intersections to the others).

\begin{defi}\label{d.nhiv}
Let $(\mathcal{D} = \coprod D_j, F =\{f_i\})$ be a Markov IFS and
$\Gamma = \cup \gamma_i \subset \mathcal{D}$ 
be a family of curves. 
We say that $\Gamma$ is a family 
of \emph{invariant curves} 
if the following holds:
\begin{center}
 $\Gamma$ is invariant under $F^{-1}$:
 For every $x \in \Gamma \cap F(\cD)$ we have 
 $F^{-1}(x) \in \Gamma$.
\end{center}
Furthermore, we say that $\Gamma$ is \emph{contracting}
if the following holds:
 There is $k_0 >0$ such that for every $x \in \Gamma \cap F^{k_0}(\cD)$ we have 
 \[
 \|DF^{-k_0}|_{T\gamma_j(x)}\| >1,
 \] 
 where 
 $DF^{-k_0}|_{T\gamma_j(x)}$ denotes the differential map
 restricted to $T\gamma_j$ and $\gamma_j$ is any curve
 of $\Gamma$ containing $x$.

We say that $\Gamma$ is \emph{univalent}
if for every $D_i$ the curve $D_i \cap \Gamma$ is 
empty or consists of single regular (unbranched) curve.
\end{defi}

We prepare a definition.

\begin{defi}\label{d.ntan}
A family of invariant curves $\Gamma$
in a Markov IFS $(\cD, F)$
is \emph{normally expanding} (resp. repelling)
if there is $k_1>0$ such that 
$TF^{-k_1}|_{\mathcal{N}\gamma_i}$ can be chosen 
uniformly greater (resp. smaller) than one
at every point where $F^{-k_1}$ is defined. 
In this definition, the normal derivative 
is the linear map induced on 
the quotient bundle 
$\mathcal{N}\gamma_i := T\mathcal{D} / T\gamma_i $.
\end{defi}

Since a family of normally repelling invariant curves are expanding in the 
normal direction, one can see the following:

\begin{rema}\label{r.relati}
If $\Gamma$ is a family of 
univalent normally repelling contracting curves, 
then there exists 
a neighborhood $R \subset \cD$ of $\Gamma$ which is a relative repelling region. 
\end{rema}

\begin{defi}\label{d.nnor}
Let $0< \eta <1$.
A family of invariant curves $\Gamma$ 
for $(\cD, F)$ is $\eta$-weak if there exists $k_1>0$
such that 
the normal derivative 
$TF^{-k_1}|_{\mathcal{N}\gamma_i}$  belongs to the open interval
$( (1-\eta)^{k_1}, (1-\eta)^{-k_1})$ at every point where $F^{-k_1}$ is defined.
We refer the number $\eta$ as the \emph{normal strength} of $\Gamma$.
\end{defi}

In the following, we want to construct a family of invariant curves
with arbitrarily weak normal strength (close to 0).
Roughly speaking, the importance of  
weak normal strength is that, by adding $\eta$-$C^1$-perturbation, 
we can produce attracting/repelling behavior. 

Let us formulate the notion of 
an attracting region.

\begin{defi}\label{d.filtrating}
Let $(\cD, F)$ be a Markov IFS and $R \subset \cD$ a relatively repelling region.
We say that a compact set $A \subset \cD$ is an 
\emph{attracting region with respect to $R$} if for any $i$ and for any
connected component $A_j$ of $A$ contained in $D_{\mathrm{dom}(i)}$ either 
$f_i(A_j)$ is contained in the interior of $A$ or $f_i(A_j)\cap R=\emptyset$ holds.
\end{defi}

This definition may seem tricky. The importance of it can be seen 
when we discuss three dimensional systems, see Section~\ref{ss.baibai}.

\begin{rema}\label{r.attrel}
If $A$ is an attracting region with respect to $R$
contained in the interior of $R$ and its connected 
components are all discs, then the family of restrictions 
$f_i|_{A_j}$ satisfying $f_i(A_j)\subset R$ defines a Markov IFS.
\end{rema}

\subsection{Constructions for weak invariant curves}
In the following, we discuss the construction of 
a relatively repelling region and an attracting region with
respect to it near a family of univalent,
contracting invariant curves with small normal strength. 

We prepare one fundamental perturbation result.

\begin{prop}\label{prop:cont}
Suppose that $(\mathcal{D}, F)$ is a Markov IFS
having a family of univalent 
invariant curves $\Gamma$ and 
let $\kappa$ be some real number.
Then there exists a set of diffeomorphism
 $\{\xi_i :D_i \to D_i\}$ which is $C^1$-$|\kappa|$-close to 
 the identity map such that the following holds:
\begin{itemize} 
\item The support of $\xi_i$ is contained in 
the arbitrarily small neighborhood of $\gamma_i$ in $D_i$.
In particular, if $\Gamma \cap D_i = \emptyset$, then 
$\xi_i$ is the identity map on $D_i$.
\item The support of $\xi_i$ is contained in the interior of $D_i$.
\item $\xi_i|_{\gamma_i} = \mathrm{id}|_{\gamma_i}$, 
where $\gamma_i$
is the connected component of $\Gamma$ in $D_i$.  
\item On $\gamma_i$, $D\xi_i|_{\mathcal{N} \gamma_i} 
= 1+ \kappa$
except some small neighborhood of the endpoints
of $\gamma_i$. This neighborhood 
can be chosen arbitrarily small.
\end{itemize}
Furthermore, we can choose the $C^0$-distance 
of $\xi_i$ and the identity map arbitrarily small.
\end{prop}

\begin{proof}
For each connected component $\gamma_i$ of $\Gamma$,
we take a smooth vector field satisfying the following:
\begin{itemize}
\item It is perpendicular to $\gamma_i$ and have the form
$\displaystyle \frac{dx}{dt} = [\log (1+ \kappa)]x$ except near the endpoints.
\item The support of the vector field is in a 
small neighborhood of $\gamma_i$.
\item The vector field is null near the endpoints.
\end{itemize}
Choosing the vector field adequately, 
the time-$1$ map 
of this vector field satisfies the conclusion.
\end{proof}

\begin{rema}
In the proof of Proposition~\ref{prop:cont},
the size of the perturbation needed to make the curves 
expanding depends only on the normal strength
and it is 
independent of the geometry of $\Gamma$, 
since we are working on the $C^1$-topology. 
\end{rema}

Now, let us state a perturbation result in the 
form of a proposition.

\begin{prop}\label{prop:exp}
Let $(\mathcal{D}, F=\{f_j\})$ be a Markov IFS having a
family of univalent invariant contracting curves $\Gamma$
with normal strength $\eta>0$, where $\eta$ is 
sufficiently close to zero. Then, there exists a set 
of diffeomorphisms $\tau_i: D_i \to D_i$ 
such that 
each $\tau_i$ is $6\eta$-close to the identity map in the 
$C^1$-topology and arbitrarily $C^0$-close to the identity 
such that the new Markov IFS 
with the discs $\mathcal{D}$ and 
the maps 
$\tilde{F}=\{\tau_{\mathrm{im}(f_j)} \circ f_j\}$
satisfies the following:
\begin{itemize}
\item $\Gamma$ is a family of univalent, contracting
invariant curves for $\tilde{F}$, too.
\item There exists a relatively repelling region 
$R \subset \mathcal{D}$ containing $\Gamma$.
$R$ can be chosen in such a way that it is contained in
an arbitrarily small neighborhood of $\Gamma$.
\item There is an attracting region $A = \cup A_i$
with respect to $R$ such that each $A_i$ 
is a $C^1$-disc which contains 
exactly one component of $\Gamma$ and $A$ contains 
$\Gamma$. Furthermore, the Markov IFS  
$(\cA = \coprod A_i, \hat{F} = \tilde{F}|_{\cA})$ 
(see Remark~\ref{r.attrel})
is uniformly contracting 
in the sense that there exists $k_1 \geq 1$ such that 
for every admissible words $\omega$ with $|\omega|=k_1$,
the inequality $\|D\hat{F}_{\omega}\| <1$ holds.
\end{itemize}
\end{prop}

\begin{rema}
The uniform contraction property of $\hat{F}$ implies 
that every periodic orbit in $(\cA, \hat{F})$ has a large stable manifold.
\end{rema}

\begin{proof}
The proof consists of two steps.

\textbf{Step 1: Construction of a relatively repelling region.}
We apply Proposition~\ref{prop:cont} to $\Gamma$
letting $\kappa = 2\eta$. We obtain a family of 
diffeomorphisms $\{\xi_i\}$ satisfying the conclusion.
Then, by direct calculation one can check that for
$F_1:= \{\xi_{\mathrm{im}(f_j)} \circ f_j \}$,
$\Gamma$ 
is a family of normally repelling curves 
with normal strength at most $3\eta$,
if $\eta$ is sufficiently close to $0$. 
Notice that the condition that $\xi_i$ is not necessarily
expanding near the endpoints does not affect the conclusion, 
since the points near the endpoints goes out from the 
discs by the backward iteration. 
By Remark~\ref{r.relati}, 
this gives us a repelling region $R$ near $\Gamma$ and
it can be chosen arbitrarily close to 
$\Gamma$. 

\textbf{Step 2: Construction of an attracting region in $R$.}
Then, for this $F_1$ we apply 
Proposition~\ref{prop:cont} letting $\kappa = -3\eta$. 
It gives another set of diffeomorphisms $\{\theta_i\}$ such that 
for the Markov IFS $F_2 := 
\{ \theta_{\mathrm{im}(f_j)}\circ \xi_{\mathrm{im}(f_j)} \circ f_j  \mid f_j \in F \}$, 
$\Gamma$ is a family of univalent,
normally attracting, contracting invariant curves. 
We choose $\{\theta_i\}$ 
in such a way that its support is contained in 
the relatively repelling region $R$ which we constructed 
in Step~1.
Since $\Gamma$ is contracting in the tangential direction, 
for each $\gamma_i$ we can find a 
disc $A_i$ in $R$ 
containing $\gamma_i$ such that $\cup A_i$
is a relative attracting region with respect to $R$ and
$F_2$ is uniformly contracting. 
Thus by letting $\tau_i := \theta_i \circ \xi_i$ we complete
the proof.
\end{proof}

\begin{rema}\label{rema:ide}
In the application of the Proposition~\ref{prop:exp},
we add one more perturbation to $\{\tau_i\}$
(see Section~\ref{ss.baibai}). 
Suppose $\Gamma$ contains several periodic orbits
$\{q_i\}$ which has large stable manifolds for $F$.
Then, we may choose $\tilde{\tau}_i$ which is 
$6\eta$-close to $\tau_i$ such that the following holds:
\begin{itemize}
\item The support of $\tilde{\tau}_i$ is arbitrarily close to 
$\gamma_i$ and the 
$C^0$-distance between $\tilde{\tau}_i$
and $\tau_i$ is arbitrarily small. 
\item For $\tilde{F}$, $q_i$ are still periodic points and 
they have large stable manifold. 
\item Near $q_i$, $\tilde{\tau}_i$ is the identity map.
\end{itemize}
Such $\tilde{\tau}_i$ can be obtained just by composing the 
inverse of $\tau_i$ near $q_i$. Since $q_i$ had large 
stable manifold for $F$ and $\tau_i$ is just a contraction 
perpendicular to $\Gamma$, 
One can make such a local deformation
keeping the largeness of the stable manifold.
\end{rema}

\subsection{Coordinate change of Markov IFS}
In the following sections, for simplifying the presentation 
we perform change of the coordinates of Markov IFSs. 
Let us briefly discuss what it exactly means. 
\begin{defi}\label{d.coord}
Let $(\mathcal{D} = \coprod D_i, F =\{f_j\})$ be a Markov IFS. 
Suppose that we have a family of diffeomorphisms $\phi_k : D_k \to D'_k$ for each $k$. 
Then, we can check that the followings also define an Markov IFS:
\begin{itemize}
\item $\mathcal{D}' = \coprod D'_i$,
\item $\{ f'_j = \phi_{\mathrm{im}(j)}  \circ f_j \circ (\phi_{\mathrm{dom}(j)})^{-1} \}$
\end{itemize} 
We refer the map  $(\coprod \phi_k): \mathcal{D} \to \mathcal{D}'$ as a \emph{coordinate change}
between Markov IFSs $(\mathcal{D}, F)$ and $(\mathcal{D}', F')$.
\end{defi}

Notice that coordinate change preserves information of dynamical systems. 
For example:
\begin{itemize}
\item If $R'$ is a relatively repelling region for $(\mathcal{D}', \{ f'_j\})$, 
then $ (\coprod \phi_k)^{-1}(R)$ is a relative repelling region for $(\mathcal{D}, \{ f_j\})$,
where $(\coprod \phi_k): \mathcal{D} \to \mathcal{D}'$ is the map 
which is defined by $(\phi_k)$ in a natural way.
\item If $\Gamma'$ is a normally repelling curves for $(\mathcal{D}', \{ f'_j\})$, 
then $ (\coprod \phi_k)^{-1}(\Gamma')$ is a relative repelling curves for $(\mathcal{D}, \{ f_j\})$.
\item If $\{g'_{k, m}\}$ is a sequence of perturbations converging to $\{ f'_j\}$, then
$\{ (\phi_{\mathrm{dom}(j)})^{-1} \circ g'_{j, m} \circ \phi_{\mathrm{dom}(j)}\}$
is a sequence of perturbations converging to $\{f_j\}$.
\end{itemize}

\subsection{Statement of the main  perturbations result}

Now we are ready to state our main result.
\begin{theo}\label{t.main}
Let $(\cD, F)$ be a Markov IFS 
 and  $\varepsilon>0$, $\eta>0$ be given. 
 Assume that $(\cD, F)$ has the following objects:
\begin{itemize}
 \item $\varepsilon$-flexible points $\{q_i\}_{i\in[1,k]}$ 
 with large stable manifolds. Each $q_i$ is separated and
 for any pair of $q_i$ and $q_j$ ($i\neq j$), they have different 
 orbits and are mutually separated (see Section~\ref{sss.per}).
 \item  $\{Q_{\ell}\}_{\ell \in[1,m]}$, a finite set of 
 $u$-homo/heteroclinic points between $q_{\ell(0)}$ and $q_{\ell(1)}$ 
 for some $\ell(0), \ell(1) \in [1, k]$ ($\ell(0)$ and $\ell(1)$ may coincide). Each $Q_l$ is $q_i$-free for all $q_i$ 
 (see Section~\ref{sss.hh}) and any pair of homo/heteroclinic 
 orbits has different itineraries.
\end{itemize}
Then, for any $\varepsilon_0>0$ and every sufficiently 
large integer $m$,
there is a $C^1$-$\varepsilon$-perturbation $G=\{g_{i}\}$ of $F$ 
which is $\varepsilon_0$-$C^0$-close to $F$ 
such that:
\begin{itemize}
 \item $\{q_i\}$ are $s$-index $1$ periodic points for $G$ whose 
 orbits coincide with that of $F$. For each $q_i$,
 the derivatives $(DG_{q_i})$ (see Section~\ref{sss.per}
 for the definition of $G_{q_i}$) along $q_i$
 coincides with
 $(B_{i; j, 1})$ where $(B_{i; j, t})$
 is some $\varepsilon$-flexible cocycle 
 (see Section~\ref{ss:flex}). Thus it has one 
 contracting eigenvalue and the other is equal to one. 
 \item Each $q_i$ has a large stable manifold 
 (it is not uniformly contracting but just topologically 
 attracting).
 \item $\{Q_{\ell}\}$ are still $u$-homo/heteroclinic points between the 
 same periodic points with the same itineraries as for $F$.
 \item There is a family of $C^1$-curves $\Gamma$ 
 containing $\{q_i\}$ and $\{Q_{\ell}\}$ such that
 \begin{itemize}
 \item $\Gamma$ is of normal strength $\eta$. 
 \item $\Gamma$ is contracting and contains $W^s_{\mathrm{loc}}(q_i)$ for every $q_i$.
 \end{itemize}
 \end{itemize}
 Furthermore, we have the following:
\begin{itemize}
\item $\Gamma \cap (G)^m(\cD)$ is univalent. 
\item The discs in $\Gamma \cap (G)^m(\cD)$
which has non-empty intersection with 
$\Gamma$ contains a point of orbits 
of $\{q_i\}$ or 
the homo/heteroclinic points $\{Q_{\ell}\}$
(see Section~\ref{sss.hh} for the definition of the 
orbit of a homo/heteroclinic point).
\end{itemize}
\end{theo}

\subsection{A simplified result}
The proof of Theorem~\ref{t.main} is 
one of the main topic of this paper.
It involves several flexible points 
and homo/heteroclinic 
points. Because of the plurality of the objects, 
a direct proof of Theorem~\ref{t.main} will be 
complicated. Thus, for the sake of simplicity 
we only give the proof of the case where
only one flexible point and 
only one homoclinic point are involved. 
Below we give it in the form of a theorem.
In Section~\ref{s.expla} we explain how we deduce
Theorem~\ref{t.main} by the proof of Theorem~\ref{t.IFS}.

\begin{theo}\label{t.IFS}
Let $(\cD, F)$ be a Markov IFS and $\varepsilon >0$, $\eta>0$
be given. 
Assume that it has the following objects: 
\begin{itemize}
\item $q$, a separated, $\varepsilon$-flexible point with 
a large stable manifold. 
\item $Q$, a $u$-homoclinic point of $q$ such that it is $q$-free.  
\end{itemize}
Then given $\varepsilon_0>0$ and every sufficiently large 
integer $m$,
there is a $C^1$-$\varepsilon$-perturbation $G$ of $F$ 
which is $C^0$-$\varepsilon_0$-close to $F$ such that:
\begin{itemize}
\item $q$ is an $s$-index 1 periodic point for $G$
with the same orbit as $F$. The derivative cocycle 
$(DG_{q})$ along the orbit of $q$ 
is $(B_{i, 1})$ of some $\varepsilon$-flexible 
cocycle $(B_{i, t})$. Furthermore, $q$ has a large stable 
manifold for $G$.
\item The point $Q$ is a $u$-heteroclinic point of $q$ with the same itinerary for $G$. 
\item There is a family of $\eta$-weak, contracting 
invariant curves $\Gamma = \cup \gamma_i$
which contains $q$, $W_{\mathrm{loc}}^s(q)$ and $Q$. 
\end{itemize}

Furthermore, $\Gamma \cap (G)^m(\cD)$ is 
univalent and each disc in $(G)^m(\cD)$
having non-empty intersection with $\Gamma$ contains 
a point of orbit of $q$ or that of $Q$. 
\end{theo}


\section{Expulsion in dimension three}
\label{s.3dim}
In this section, we complete the proof of theorems 
which we presented in Introduction assuming the 
main technical result Theorem~\ref{t.main}. 

In Section~\ref{ss.prepara} we give some auxiliary 
perturbation results, which are essentially proved in 
papers such as \cite{BCDG, BS1, BS2}.
In Section~\ref{ss.baibai} we prove Theorem~\ref{t.circ}.
The proof goes as follows: First we reduce the 
problem into the one about Markov IFSs. We apply 
Theorem~\ref{t.main} to the reduced problem 
and, by applying Proposition~\ref{prop:3dim}, 
we realize the perturbation 
in Theorem~\ref{t.main} in dimension three,
which concludes Theorem~\ref{t.circ}.
In Section~\ref{ss:viral} we prove 
Theorem~\ref{theo:expu} and \ref{theo:viral} using 
Theorem~\ref{t.affine}, \ref{t.circ} and 
the results in Section~\ref{ss.prepara}.
Finally, in Section~\ref{ss:appear} we complete 
the proof of Theorem~\ref{t.aperi} by improving the 
proof of Theorem~\ref{theo:viral}. 

\subsection{Auxiliary results}
\label{ss.prepara}
In this subsection, we prepare some preparatory results.

\subsubsection{Abundance of flexible points}
In \cite{BS2} we proved a result about the 
existence of flexible points with large stable manifolds.
We cite it with a small modification.
In order to state it, we prepare some definitions.
Let $U$ be a subset of closed three dimensional 
manifold $M$ and $p, q$ be two hyperbolic 
periodic points of the same $s$-index whose orbits 
are contained in $U$.
We say that $p, q$ are 
\emph{homoclinically related in $U$} if
there are 
heteroclinic orbits contained in $U$ from $p$ to $q$ 
and vice versa. 
The \emph{relative homoclinic class} $H(p, U)$
is the closure of the periodic points in $U$ which are homoclinically related to $p$ in $U$. 

For a hyperbolic periodic point $p$ of $s$-index $2$, we say 
that it is $\varepsilon$-flexible if the linear 
cocycle obtained by restricting the derivative cocycle 
to the stable direction along $\cO(p)$ is $\varepsilon$-flexible, 
see Section~\ref{ss:flex} for the definition of 
$\varepsilon$-flexible cocycles.

Let $p$ be an $s$-index two hyperbolic periodic point of  
$f \in \mathrm{Diff}^1(M)$.
We say that 
$p$ has \emph{a robust heterodimensional cycle in $U$} 
if the following holds (see also Proposition~5.1 of \cite{BS1}): There are hyperbolic basic sets 
$\Lambda$ and $\Sigma$ in $U$ such that: 
\begin{itemize}
\item $\Lambda$ is $s$-index two and $\Sigma$ is 
$s$-index one.  
\item There is a $C^1$-neighborhood $\cU$ of $f$ 
such that for every $g \in \cU$ the continuations 
$p_g$, $\Lambda_g$ and $\Sigma_g$ are defined and  
contained in $U$. Furthermore, $p_g \in \Lambda_g$ holds.
\item For every $g \in \cU$, there are 
heteroclinic points in $W^s(\Lambda) \cap W^u(\Sigma)$ 
and $W^s(\Sigma) \cap W^u(\Lambda)$ whose orbits 
are contained in $U$.
\end{itemize}
Now we are ready to state the result.

\begin{defi}
Let $f$ be a $C^1$-diffeomorphism of a three dimensional 
manifold having a filtrating Markov partition $\mathbf{R}$. 
Let $\mathbf{W}$ be a sub Markov partition of $\mathbf{R}$
(i.e., a collection of rectangles of $\bfR$)
such that there is a hyperbolic periodic point $p$
whose orbit is contained in $\mathbf{W}$.
We say that the relative homoclinic class $H(p, \mathbf{W})$
satisfies property $(\ell_{\mathbf{W}})$ if
the following holds:
\begin{itemize}
\item $p$ has a large stable manifold in $\bfR$.
\item There is a hyperbolic periodic point $p_1$ whose orbit 
is contained in $\mathbf{W}$ such that $p$ and $p_1$ are homoclinically 
related in $\mathbf{W}$
and $p_1$ has a stable non-real eigenvalue.
\item $p$ has a robust heterodimensional cycle in $\mathbf{W}$.
\end{itemize}
\end{defi}

Now let us give the result.
\begin{prop}[See Proposition~5.1 in \cite{BS1} and 
Lemma~3.8 in \cite{BS2}]
\label{l.flexiblelarge} 
Let ${\bf R} = \bigcup C_i$ be a filtrating Markov 
partition of a diffeomorphism $f$ and
$\cU$ a $C^1$-neighborhood of $f$ in which one can find a 
continuation of ${\bf R}$.
Assume that there are hyperbolic periodic point $p\in {\bf R}$ 
and a sub Markov partition $\mathbf{W}$ of ${\bf R}$
such that 
for every $\tilde{f}\in \cU$ the 
relative homoclinic class 
$H(p_{\tilde{f}}, \mathbf{W};\tilde{f})$ 
satisfies condition $(l_{\mathbf{W}})$.

Then, for any $\varepsilon>0$ there is a $C^1$-open and 
dense subset $\cD$ of $\cU$ such that every  
diffeomorphism $g\in \cD$ has a periodic point 
$x\in H(p_g, \mathbf{W}; g)$ of $s$-index two, 
with a large stable manifold, $\varepsilon$-flexible, 
homoclinically related with $p_g$ in $\mathbf{W}$, 
and whose orbit is $\varepsilon$-dense 
in $H(p_g, \mathbf{W}; g)$.
\end{prop}

\begin{rema}
In Proposition~\ref{l.flexiblelarge}, we may assume that 
for any $N>0$ the flexible point $x$ has period larger than $N$
(choose $x$ letting $\varepsilon$ sufficiently small).
\end{rema}

The difference from the original statement is that 
the assumption of the condition $(l)$ is stated for  
relative homoclinic classes and the conclusion holds
for the relative homoclinic classes.  Since the 
argument used in the proof is local ones, one can obtain 
this result just by following the proof line-by-line, including 
the result Proposition~5.2 of \cite{BS1}.

By adding a small modification, we obtain the following:

\begin{coro}\label{c.flexlarge}
Let $f$, $\mathbf{R}$ and $H(p,\mathbf{W};f)$ be as in the
assumption of Proposition~\ref{l.flexiblelarge}. 
Then, there is a $C^1$-diffeomorphism $h$ which 
is $C^1$-arbitrarily close to $f$ such that 
$h$ satisfies the conclusion of the 
Proposition~\ref{l.flexiblelarge} and 
$h \equiv f$ holds outside $\mathbf{W}$.
\end{coro}

\proof{
First, by applying Proposition~\ref{l.flexiblelarge},
take a sequence of diffeomorphisms $(f_n)$ which 
converges to $f$ such that each $f_n$
satisfies the conclusion. Then, $(f_n)^{-1} \circ f$ is a $C^1$-diffeomorphism 
which converges to the identity map in the $C^1$-topology. 
For $\mathbf{R}$, we take a $(1, 1)$-refinement 
$\mathbf{R}' := f(\mathbf{R}) \cap f^{-1}(\mathbf{R})$
and set $\mathbf{W}' = \mathbf{R}' \cap \mathbf{W}$. 
Note that, due to the filtrating property of 
$\mathbf{R}$, $\mathbf{W}'$ contains 
$H(p_{f_n}, \mathbf{W};f_n)$ for sufficiently large $n$. 

Then, we consider the diffeomorphism $(f_n)^{-1} \circ f$
for sufficiently large $n$. 
In the 
following, we will show that for $n$ sufficiently large, we can find 
a diffeomorphism $g_n$ satisfying the following:
\begin{itemize}
\item On $\mathbf{W}'$, $g_n$ is the identity map,
\item $g_n = (f_n)^{-1} \circ f $ outside $\mathbf{W}$,
\item $(g_n)$ converges to the identity map
in the $C^1$-topology as $n\to \infty$. 
\end{itemize}
For the time being, assuming the existence 
of such $g_n$ let us conclude 
the proof. Consider the diffeomorphism $h_n :=  f_n\circ g_n$.
By definition, one can see that 
$h_n = f_n \circ (f_n)^{-1} \circ f  = f$ outside 
$\mathbf{W}$. Furthermore, 
on $\mathbf{W}'$ we have $h_n = f_n$. Since the relative 
homoclinic class $H(p, \mathbf{W}; h_n)$ is 
determined by the behavior of the dynamics on $\mathbf{W}'$, 
we know that $H(p, \mathbf{W}; h_n) = H(p, \mathbf{W}; f_n)$. Note that 
$\varepsilon$-flexibility is a local property and the largeness of the 
stable manifold can be determined by the behavior on $\mathbf{W}'$.
Thus for $h_n$ we still have the periodic point
which is $\varepsilon$-flexible, 
having a large stable manifold and
$\varepsilon$-dense in $H(p, \mathbf{W}; h_n)$.

Now let us construct $(g_n)$. For that we 
first fix a smooth bump function 
$\kappa :M \to [0, 1]$ which takes value $1$ outside $\mathbf{W}$ and $0$
on $\mathbf{W}'$.
Now, we follow the classical construction of representing a 
diffeomorphism close to the identity map by a vector field, 
for instance see \cite{Les}.  
We fix a 
smooth Riemannian metric. Then we consider a map 
$TM \to M \times M$ which sends $(x, v) \in M \times T_xM$
to $(x, \mathrm{exp}_x(v))$. This is a diffeomorphism
in the neighborhood of the image of the zero section in $TM$. 
Now given a 
$C^1$-diffeomorphism $f$ which is sufficiently $C^1$-close
to the identity, we can associate a $C^1$-vector field $F$ 
such that for every $x \in M$ we have
$\mathrm{exp}_x(F(x)) = f(x)$ holds. 

By applying this construction
to $(g_n)$, 
we take a sequence of $C^1$ vector fields 
$(G_n)$ whose image under the 
exponential map is $(g_n)$. Now consider the vector field 
$ \kappa  G_n$ and take its image under above correspondence. 
Note that applying $\mathrm{exp}(\,\cdot\,)$ is continuous 
with respect to the $C^1$-topology.
Thus this defines the desired sequence 
of diffeomorphisms $h_n$.
}

\begin{rema}
In general, a relative homoclinic class $H(p, U)$ may behave badly 
under perturbation. 
In this article, we deal with the case where $U$
is a sub Markov partition $\mathbf{W}$ 
of a filtrating Markov partition $\mathbf{R}$. 
In such a case, due the filtrating property of $\mathbf{R}$,
we know that every point of $H(p, \mathbf{W})$ has uniform distance 
from the boundary of $U$ and it enables us to deal 
$H(p, \mathbf{W})$
as if it is a homoclinic class in an ambient manifold.
\end{rema}

\subsubsection{Flexibility implies 
condition ($\boldsymbol{\ell}$)}
We present a result which recover the property $(\ell)$ 
for an invariant set containing a flexible periodic point.
We prove it under a local setting (that is, for property 
$(\ell_{\mathbf{W}})$).
For the purpose of this paper, 
the version for property $(\ell)$ is enough, 
but for the future use we provide the proof under 
more general settings.

The following result can be proved 
by arguments based on 
the results in \cite{BS1} and \cite{BCDG}.  
\begin{prop}\label{prop:recep}
Let $f \in \mathrm{Diff}^1(M)$ having
an $\varepsilon$-flexible periodic point $p$ in a filtrating 
Markov partition $\mathbf{R}$ with 
a large stable manifold.
Let $\mathbf{W}$ be a sub Markov partition 
of $\mathbf{R}$.
Assume that $\mathbf{R}$ is $4\varepsilon$-robust and
$H(p, \mathbf{W})$ is non-trivial.
Then, given $\delta >0$
there is $g = g_\delta \in \mathrm{Diff}^1(M)$ 
which is $C^1$-$4\varepsilon$-close to $f$ 
such that the following holds:
\begin{itemize}
\item $p$ is still an $\varepsilon$-flexible point 
with the same orbit for $g$.
\item $H(p_g, \mathbf{W};g)$ 
satisfies condition $(\ell_{\mathbf{W}})$. 
\item The support of $g$ is contained in $\mathbf{W}$.
Especially, $\mathbf{R}$ is still a filtrating Markov 
partition and $\mathbf{W}$ is its sub Markov partition.
\item Suppose that $p$ is contained in a circuit $K$. 
Then, for any $\delta>0$ by choosing 
appropriate $g = g_\delta$ there is a circuit $K_g$ which is 
$\delta$-similar to $K$. 
\end{itemize}
\end{prop}
\begin{proof}
For the proof we need to construct several objects
(a periodic point with non-real eigenvalue 
and a robust heterodimensional cycle) 
by small perturbation keeping the largeness of the 
stable manifold and 
the similarity of the circuit small. Such a process is already 
well described for instance in the proof of 
Proposition~5.2 and Corollary~5.4 of \cite{BCDG}. 
Thus we only give 
the sketch of the proof.
First we explain how to construct these objects. In the 
end of the proof we will see how to guarantee
the largeness of the stable manifold and 
the smallness of the similarity of the circuit.

First, let us see how to construct a periodic point 
with complex eigenvalues.
Recall that the flexibility of $p$
guarantees the existence of a path of linear hyperbolic 
cocycles of which connects 
$\left( Df(f^i(p)) \right)$ and 
a cocycle whose product has non-real eigenvalues.
Thus, if we deform $f$ along $p$ which gradually 
realizes the path, we can change $p$ so that it 
has non-real eigenvalues into $E^{cs}$-direction.
It guarantees the existence of $f_1$ which is an 
$\varepsilon$-perturbation along $\mathcal{O}(p)$
such that $p$ has the 
same orbit and the eigenvalues to the 
$E^{cs}$-direction are non-real, keeping the non-triviality 
of the homoclinic class.
Note that due to the existence of such $p$, we know 
that $H(p, \mathbf{W}; f_1)$ does not admit any dominated
splitting in $E^{cs}$-direction.

Now we add a perturbation whose $C^1$-size 
can be chosen arbitrarily small such that 
in the relative homoclinic class
$H(p, \mathbf{W})$ there is a hyperbolic periodic point 
which is not (the continuation of) $p$ and
has non-real eigenvalues in $E^{cs}$-direction. 
This process is explained in the proof of 
Proposition~5.2 of \cite{BCDG}. 
The only difference is again that we work on relative 
homoclinic classes, but the local feature of the argument 
enables us to prove this. We denote the perturbed 
diffeomorphism by $f_2$.

Now, we give another perturbation around $p$
to obtain a diffeomorphism $f_3$ which returns $p$ 
into an $\varepsilon$-flexible point, keeping the existence
of periodic points with non-real eigenvalues. 
Such a perturbation can be done
by following the path of cocycles used in the previous 
step in the opposite direction. Up to now, the amount 
of the size of the perturbation is $2\varepsilon$.

\medskip

In the following, we give another sequence of
perturbations to obtain 
a robust heterodimensional cycle.
First, using the $\varepsilon$-flexibility of the periodic point $p$, we obtain a
diffeomorphism $f_4$ such that $p$ is almost an $s$-index 1
hyperbolic periodic point.
This can be done by following the path of cocycles to the 
direction of $t=1$.
 Then, by Proposition~5.2 
of \cite{BS1}, we know that, up to an arbitrarily small 
perturbation we can find an 
$\varepsilon'$-flexible point homoclinically 
related to $p$, say $r$, where $\varepsilon'$ 
can be chosen arbitrarily close to zero.

Then, as we obtained $f_3$ from $f_2$,
we give an $\varepsilon$-perturbation around $p$ 
to obtain a diffeomorphism $f_5$ such
that $p$ is again $\varepsilon$-flexible, without disturbing 
the existence of periodic orbit having 
complex eigenvalues and 
the $\varepsilon'$-flexible periodic 
point $r$. The size of the perturbation from $f_3$
to $f_5$ is also bounded by $2\varepsilon$.

Now, we construct a robust heterodimensional cycle: 
Using the $\varepsilon'$-flexibility of $r$, 
we give an $\varepsilon'$-perturbation 
around $r$ such that $r$ is a stable index $1$
periodic point whose strong stable 
manifold has a non-empty intersection with the 
unstable manifold of some hyperbolic 
periodic point of $s$-index 2, say $p'$, homoclinically related 
to $p$. Such a perturbation is possible due to the flexibility, 
see Theorem~1.1 of \cite{BS1}. 

Then, by \cite{BD2}, we may assume that 
the heterodimensional cycle turns to be robust up to an 
arbitrarily small perturbation (since the homoclinic 
class of $p'$ is non-trivial, we can apply 
Theorem~5.3 of \cite{BD2}). 
Thus the relative
homoclinic class of $p$ now $C^1$-robustly satisfies 
the condition $(\ell_{\mathbf{W}})$.
The size of the last perturbation is $\varepsilon'$ and 
it can be chosen arbitrarily close to zero. As a result, the 
total amount of the size of the perturbation is less than
$2\varepsilon + 2\varepsilon = 4\varepsilon$.

\medskip

Finally, let us see how to obtain the similarity of the circuit
and the largeness of the periodic points after perturbation.
For the $\delta$-similarity of the circuit, notice that 
the perturbation we performed is either arbitrarily small
or a perturbation around the periodic point using the 
flexibility of the point. For the first 
one by decreasing the size of the perturbation we can guarantee 
the $\delta$-similarity by continuity. For the second one 
we use the 
``adapted perturbation'' (see Section~3 and the proof 
of the Corollary~5.4 of \cite{BCDG}), which 
preserves bounded part of the invariant manifolds. 

A circuit consists of periodic 
orbits and 
hetero/homoclinic orbits.
The periodic orbits are preserved under perturbation 
based on the flexibility property. For the other perturbation 
the size can be chosen arbitrarily small, thus the change 
of the orbits can be made arbitrarily small, too.
For the hetero/homoclinic 
points, note that they converge to periodic orbits
(both forwardly and backwardly). Thus the effect 
of the perturbation on the homo/heteroclinic orbits 
near the periodic orbits are always bounded. 
As a result, we see that 
the invariance of bounded part of the invariant
manifold is enough to guarantee the smallness of the variation 
of the position of the orbits. Thus, by requiring the invariance
of the invariant manifold in the fixed region we can
guarantee the $\delta$-similarity of the circuit. 

Also, notice that the largeness of the 
stable manifold is determined by the information of the invariant 
manifold contained in a bounded part. Thus, 
by using the adapted perturbation we 
can keep the largeness of the manifold. 

This concludes the proof. 
\end{proof}

By using the same argument in the proof of 
Proposition~\ref{c.flexlarge}, we have the following:
\begin{coro}
In the same hypothesis, we can choose a $C^1$-diffeomorphism
$h$ which is $4\varepsilon$-close to $f$ and coincides 
with $f$ outside $\mathbf{W}$ such that 
$H(p_h, \mathbf{W};h)$ satisfies the conclusion.
\end{coro}

\subsubsection{Recovering the flexibility}

The following result is used to recover the 
flexibility of the periodic points. 
\begin{prop}\label{p.reco}
Let $p$ be a periodic point of a Markov IFS
$(\mathcal{D}, F=\{f_i\})$ having a 
large stable manifold. Suppose that 
the cocycle of the differentials $(DF_{p_i})$ coincides with
$(B_{i, 1})$ where $(B_{i, t})$ is some $\varepsilon$-flexible
cocycle. Then, there exists an $\varepsilon$-perturbation 
$G$ of $f$ along the orbit of $p$ such that the following holds:
\begin{itemize}
\item The support of the perturbation is contained in 
an arbitrarily small neighborhood of $\{p_i\}$.
\item $p$ is a periodic point for $G$ with the same itinerary.
\item $p$ has a large stable manifold for $G$, too.
\item $DG_{p_i} = B_{i, 0}$ for every $i$.
\end{itemize}
\end{prop}
The proof can be done by repeating the argument 
in the proof of Proposition~4.1 of \cite{BS1}, 
which also appeared in 
the proof of Proposition~\ref{l.flexiblelarge}. 
Thus we just give a short account of the proof. 
For the detail, see Section~4 of \cite{BS1}.

By the definition of the $\varepsilon$-flexible cocycle, 
there is a path of cocycles $(B_{i,t})$ which connects 
$(B_{i, 0})$ and $(B_{i, 1})$ such that it is uniformly hyperbolic 
for $0 \leq t \leq 1$ and has the size smaller than $\varepsilon$.
Then, by realizing this path slowly, we can deform $(B_{i, 1})$
into $(B_{i,0})$, keeping the largeness of the stable manifold.
manifold. It gives the perturbation $G$ we desired.

\subsection{Expulsion of the circuit}
\label{ss.baibai}
Using the results in Section~\ref{ss.prepara}, together with the linearization result, 
we prove Theorem~\ref{t.circ}. 
\proof[Proof of Theorem~\ref{t.circ}]{
\textbf{Step 1: reduction to Markov IFSs.} 
Let $f \in \mathrm{Diff}^1(M)$ having 
a filtrating Markov partition $\bfR$ of 
robustness $\alpha$, and
suppose that we have a circuit $S$ in $\bfR$ 
consisting of $\varepsilon$-flexible 
periodic orbits $\{\mathcal{O}(q_i)\}$ with large 
stable manifolds
and 
homo/heteroclinic orbits $\{\mathcal{O}(Q_j)\}$ such that
$\mathbf{R}(S)$ is an affine Markov partition.

Then as is explained in Section~\ref{ss.corres}, 
we construct a Markov IFS $\cM(\bfR(S))$.
We denote it by $(\mathcal{D}, F)$.
We have corresponding periodic points 
and homo/heteroclinic points for it. 
Recall that these periodic points are separated, 
any pair of them are mutually separated a
and every homo/heteroclinic orbit is free from 
the periodic orbits (see Section~\ref{ss.corres}).
Also, the generating property of the Markov partition 
ensures that every pair of 
homo/heteroclinic points has different itineraries
(see Lemma~\ref{l.shadow}).

\bigskip

\textbf{Step 2: Solving the two dimensional problem.} 
Then we apply Theorem~\ref{t.main}. We set 
$\varepsilon$ to be the size of the flexibility of $\{q_i\}$
and $\eta$ to be sufficiently small number (which 
we will fix later). Then 
for every sufficiently large integer $n$ we can find 
an $\varepsilon$-perturbation $G_{0, n} = G_0$ of 
the two dimensional maps $F$ such that
\begin{itemize}
\item $\{q_i\}$ are still periodic points with the same orbits.
\item There is a family of contracting 
invariant curves $\Gamma$ 
containing $q_i$ and $Q_j$.
\item $\Gamma$ is univalent 
in the $n$-refinement $G_0^n(\mathcal{D})$. 
\item The normal strength of $\Gamma$ is smaller than 
$\eta$.
\item The $C^0$-size of the perturbation is less than $\varepsilon_0$, which can be chosen arbitrarily small.
\end{itemize}

Now we apply Proposition~\ref{prop:exp} to this 
family in the $(0,n)$-refinement: Then in 
$(G_0^n(\cD), \wedge_n G_0)$ we can find a family of diffeomorphisms
$\{\tau_i\}$ which is $6\eta$-close to the identity map such that 
$G_1:= \{\tau_{\mathrm{im}(g_j)} \circ f_j\}$ has
a relatively repelling region $R$ and an attracting region 
$\cA=\cup A_i$ with respect to $R$
such that $(\cA, G_1)$ defines a new Markov IFS
containing $\{q_i\}$ and $\{Q_j\}$.
Since $G_1$ is contracting on $\cA$, 
we have that each $q_i$ has a large stable manifold 
in $\cA$.

Remark~\ref{rema:ide} tells us that
by performing another $6\eta$-$C^1$-small perturbation,
we may assume that $\tau_{i}$ is the identity map 
near the orbit of $q_i$, keeping the largeness of the stable 
manifold. We denote the perturbed IFS by $G_2$.
Finally, using Proposition~\ref{p.reco} we perform another 
$\varepsilon$-small perturbation
which makes $q_i$ to be $\varepsilon$-flexible, keeping 
the largeness of the stable manifold. We denote the 
obtained IFS by $H_n$. Note that the support 
of $H_n$ is contained in $G_0^n(\cD)$. Thus one 
can consider $H_n$ as a perturbation of $F$ as well.  
Then the amount of 
the size of the perturbation between $F$ and $H_n$ is less than
\[
\varepsilon + 6\eta + 6\eta + \varepsilon 
= 2 \varepsilon +12\eta,
\] 
and $\eta$ can be arbitrarily small. Thus, 
letting $\eta$ sufficiently small, the size 
of the perturbation is less than $2 \varepsilon$.
Note that the $C^0$-distance between $F$ and $H_n$ can be chosen arbitrarily small. We denote it by $\delta$.

\bigskip

\textbf{Step 3: Expulsion in dimension 3.} 
Let $f$ denote the three dimensional map
in Step 1. 
Now we perform a perturbation to the three 
dimensional diffeomorphism $f$ by 
Proposition~\ref{prop:3dim}. 
We apply Proposition~\ref{prop:3dim} to the
Markov IFS $(\mathcal{D}, H_n)$ which we obtained in 
Step 2:
We can find $h_n$ which 
is $(2\varepsilon + K\delta)$-close to 
$f$ such that 
$h_n$ still keeps the product structure on the 
rectangles and the corresponding Markov IFS
is given by $H_n$. Since the size of the perturbation 
in the $C^0$-distance can be made arbitrarily small and 
$K$ is already fixed,
we may assume that $h_n$ is indeed 
$2\varepsilon$-$C^1$-close to $f$. 
Throughout this proof, this is the last part where we perform
the perturbation. 

Note that by the property of $H_n$ we have:
\begin{itemize}
\item The support $\supp (h_n, f)$ 
is contained in $\mathbf{R}(S)$,
\item $h_n$ has the continuation of the circuit $S_n$, 
whose periodic orbits has the same orbits as in $S$,
\item the periodic orbits of $S_n$ are all $\varepsilon$-flexible.
\item $S_n$ is similar to $S$; 
for the homo/heteroclinic points of $S_n$, we only need 
to collect the points which corresponds to $f^{T_j}(Q_j)$.
Since the points $\{f^{T_j-k}(Q_j)\}_{k>0}$ does 
not get any change under the perturbation and 
$\{(h_n)^{T_j+k}(Q_j)\}_{k\geq0}$ belong to the 
local stable manifold of some periodic points, we can 
choose the conjugacy in such a way that the corresponding 
points in $S$ and $S_n$ belong to the same rectangle.
\end{itemize}
For this $h_n$ we can construct a new filtrating set. 
First, recall that the filtrating Markov partition
$\mathbf{R} = \cup C_i$
has an attracting set $A$ and a repelling 
set $R$ such that $\cup C_i = A \cap R$.
Then we consider its $(0, n)$-refinement
$\mathbf{R}_{(0,n)}$ with respect to $h_n$,
which has the form $(h_n)^n(A) \cap R$.
We denote the rectangles of it by $\{D_{j}\}$.

The two dimensional dynamics of $h_n$ on 
these rectangles is given by the iterated function system $H_n$.
Recall that $H_n$ has a relative repelling region for the 
$(0, n)$-refinement. 
Thus, for each $D_j$ there is a corresponding three dimensional
set which projects to the relative repelling region. 
We denote it by $\widehat{D_{j}}$. 

Now we define a repelling set as follows:
\[
R' = \big( R \setminus (\cup D_i) \big) 
\cup (\cup \widehat{D_j}).
\] 
By construction one can check that $R'$ is a repelling set for 
$h_n$. Thus $(h_n)^n(A) \cap R'$ is a 
filtrating set containing $S_n$.

Now we choose the attracting set.
Recall that for $H_n$
we have an attracting set with respect to the repelling 
set corresponding to $\{\widehat{D_j}\}$ contained in 
$\{D_j\}$. We denote them by $\{A_j\}$, take the corresponding 
three dimensional sets and denote them by 
$\{\widehat{A_j}\}$. 

Then, one can see that 
\[
A' = \big( \overline{(h_n)^n(A) \setminus (\cup \widehat{D_i})} \big) \cup (\cup \widehat{A_j})
\] 
is an attracting set for $h_n$. Now, after smoothing the
corners of $R'$ and $A'$ appropriately, one can see that 
$\mathbf{R}'' = R' \cap A' = \cup \widehat{A_j}$ 
satisfies the condition of filtrating Markov partitions
except the existence of the cone field. 
To confirm the existence of the cone field, 
we need to consider 
the robustness of $\mathbf{R}$. 
Note that for $\mathbf{R}_{(0,n)}$ and $f$,
there is an invariant cone field inheriting from $\mathbf{R}$
whose robustness is $\alpha$ which is greater than $2\varepsilon$.
Then, since 
each rectangle in $\{\widehat{A_j}\}$ is product 
rectangle in the linearized coordinate of $\mathbf{R}_{(0,n)}$,
together with Remark~\ref{r.concon} we see that the restriction 
of the cone field of $\mathbf{R}$ 
to $\mathbf{R}''$ gives the vertical 
cone field of robustness $\alpha -2\varepsilon$. 

Now we conclude that $\mathbf{R}''$ is 
a filtrating Markov partition.
Recall that the attracting region in the Markov IFS is
obtained as the neighborhood of the 
family of normally contracting 
invariant curves $\Gamma$
in the Markov IFS and it is univalent in the $n$-refinement,
that is, it has one and only one connected component.
This shows that
each $\mathbf{R}_{(0,n)}(S_n)$ contains 
one and only one rectangle in $\{\widehat{A_j}\}$.
Also, note that by the property of $H_n$ all the 
periodic orbits in $S_n$ have 
large stable manifolds in $R' \cap A'$.

Let us confirm the c-transitivity of the rectangles. 
For that we only need to confirm for each pair of
rectangles $C_1$, $C_2$ containing a periodic point,
we can find a path of connected components connecting 
them
(for other rectangles, it must contain homo/hetero clinic point
but they can be connected to rectangles having periodic 
rectangles).
Let $q_1$, $q_2$ be the periodic points which $C_1$, $C_2$ contain 
respectively. If they have the same orbit the conclusion is
straightforward. If not, 
by the definition of the transitivity of the circuit of points
 (recall that 
we assume that every circuit of points is transitive), 
there are sequence of homo/heteroclinic
points $\{x_{j,l}\}$ and periodic points $\{q_{j,l}\}$ 
contained in the circuit connecting 
$C_1$ and $C_2$. Then, by following these connection we can 
find desired chain of connected components.  
It concludes the proof.}

%

\subsection{Proof of the viralness}
\label{ss:viral}
In this subsection let us see how
to obtain Theorem~\ref{theo:expu}
using Theorem~\ref{t.affine} and Theorem~\ref{t.circ}.

\begin{proof}[Proof of Theorem~\ref{theo:expu}]
Given a chain recurrence class $C(p)$, 
suppose that we have a circuit 
$S$ contained in the filtrating Markov partition $\mathbf{R}$
satisfying the assumption. Then, first we 
apply Theorem~\ref{t.affine} to $f$ such that 
up to an arbitrarily small perturbation we may assume
in a sufficiently fine refinement $\mathbf{R}'$
we have that $\mathbf{R}'$ is generating and
$\mathbf{R}'(S)$ is an affine Markov partition. 
Note that, by using the largeness of the stable manifolds
for every periodic orbits in $S$, we may assume that
$\mathbf{R}'(S)$ do not contain the orbit of (continuation)
of $p$ and the diameter of each rectangle of 
$\mathbf{R}'(S)$ is less than $\delta$,
see Remark~\ref{r.tiny}.
Then we apply Theorem~\ref{t.circ}:
It gives, up to an $2\varepsilon$-perturbation
us a new filtrating Markov partition $\mathbf{R}''$
containing a circuit $S'$ which is $\delta$-similar to $S$
where $\delta$ can be chosen arbitrarily small.
By construction, we know that this new $\mathbf{R}''$
satisfies all the conclusions claimed, that is, c-transitivity,  
the largeness of the stable manifold and the 
$\varepsilon$-flexibility of periodic points.
\end{proof}

Now the proof of Theorem~\ref{theo:viral} is immediate.
\begin{proof}[Proof of Theorem~\ref{theo:viral}]

Suppose we have a $C^1$-diffeomorphism $f$ having a chain 
recurrence class $C(p)$ contained in a filtrating Markov partition
$\mathbf{R}$
satisfying property $(\ell)$. 
Note that this implies $H(p, \bfR)$ satisfies property $(\ell_{\bfR})$.
By choosing small $\varepsilon>0$
we may assume that $\mathbf{R}$ is $6\varepsilon$-robust.
By Proposition~\ref{l.flexiblelarge}, up to an
arbitrarily small perturbation we may assume that 
$H(p; \mathbf{R})$ contains an $\varepsilon$-flexible point
with a large stable manifold which is not equal 
to $p$, say $q$. Since both $p$ and 
$q$ have large stable manifolds, we know that
$H(q; \mathbf{R})$ is not trivial. In particular, we can 
find a circuit $S$ which consists of the orbit of $q$ and a 
homoclinic orbit of $q$. 
Then we apply Theorem~\ref{theo:expu} to $S$:
We can find a $2\varepsilon$-perturbation of $f$ such that 
there is a filtrating Markov partition $\mathbf{R}'$
containing the continuation $S'$ of $S$ and disjoint from 
the continuation of $p$ such that the continuation of $q$
has a large stable manifold and is $\varepsilon$-flexible.
We have that $\mathbf{R}'$ 
is $6\varepsilon-2\varepsilon = 4\varepsilon$
robust.

Now we apply Proposition~\ref{prop:recep}: 
Up to $4\varepsilon$-perturbation we know that 
the relative homoclinic class $H(q, \mathbf{R}')$,
satisfies the property $(\ell_{\mathbf{R}'})$.
Hence the chain recurrence class $C(q)$ satisfies
the property $(\ell)$ 
the filtrating Markov partition $\mathbf{R}'$.

In short, up to $6\varepsilon$-perturbation, we have 
constructed a new chain recurrence class $C(q)$ satisfying 
the property $(\ell)$ and $\varepsilon >0$ can be chosen
arbitrarily small.
It concludes the proof of the viralness.
\end{proof}

\begin{rema}\label{r.largeperi}
In the above proof, we may assume that the 
period of $q$ is larger than any prescribed integer. 
This is possible by letting $\varepsilon$ close to zero,
for $\varepsilon$-flexible points with small $\varepsilon$
must have a large period.
\end{rema}

\subsection{Construction of aperiodic classes}
\label{ss:appear}
Finally, let us give the proof of Theorem~\ref{t.aperi}.
For the construction of the aperiodic classes, we prove the 
following:
\begin{prop}\label{prop:baiin}
Let $f$ be a diffeomorphism having a filtrating Markov partition
containing a 
chain recurrence class $C(p)$ satisfying the property $(\ell)$.
Let $\mathcal{W}$ be 
a $C^1$-neighborhood of $f$ in which we have 
the continuation of $C(p)$ keeping the property $(\ell)$. Then, there exists a 
nested sequence of $C^1$-open sets $\{\mathcal{O}_n\}$
satisfying the following:
\begin{itemize}
\item Each $\mathcal{O}_n$ is dense in $\mathcal{W}$.
\item For each $n \geq 1$, there exists a function 
$\mathcal{U}_{n} :\mathcal{O}_n \to \mathcal{K}(M)$
where $\mathcal{K}(M)$ is the set of
all compact subsets of $M$
such that the following holds:
\begin{itemize}
\item $\mathcal{U}_{n}$ is locally constant in 
$\mathcal{O}_n$.
\item $\mathcal{U}_{k+1}(f) \subset \mathcal{U}_{k}(f)$
for $k\geq 1$, $k \leq n$ and $f \in \mathcal{O}_n$. 
\item $\mathcal{U}_{n}(f)$ is a  
$\mathrm{c}$-transitive
filtrating set of 
$f$. 
Every connected component of $\mathcal{U}_{n}(f)$
has the diameter less than $1/n$.
\item There is a periodic point $q_n$ in 
$\mathcal{U}_{n}(f)$ such that $C(q_n)$ satisfies the 
property $(l)$.
\item For every $r \in \mathcal{U}_{n}(f)$, 
if it is a periodic point of $r$, then the period is larger than $n$.
We call this property \emph{$n$-aperiodic}.
\end{itemize}
\end{itemize}
\end{prop}
First, let us complete the proof of Theorem~\ref{t.aperi}
assuming Proposition~\ref{prop:baiin}.
\begin{proof}
Suppose we have such $\{\mathcal{O}_n\}$. 
Then $\mathcal{R} := \cap \mathcal{O}_n$ is 
residual in $\mathcal{W}$. Take $f \in \mathcal{R}$.
Then there is an infinite nested sequence of filtrating regions 
$\mathcal{U}_1(f) \supset \mathcal{U}_2(f) \supset \cdots$.

Now consider $C' :=\cap \mathcal{U}_i(f)$.
The c-transitivity for each $\mathcal{U}_k(f)$ 
and the smallness of 
the connected component for each $\mathcal{U}_k(f)$
for $k$ large implies 
that $C'$ is chain-transitive. Furthermore, since 
each $\mathcal{U}_i(f)$ is a filtrating set, we know that 
$C'$ is a chain recurrence class. Finally, the condition 
about the period implies the aperiodicity of $C'$. 
\end{proof}

Now let us prove Proposition~\ref{prop:baiin}.
A key step is a version of Theorem~\ref{theo:expu} 
which delivers stronger conditions. 

\begin{prop}\label{prop:virain}
In Theorem~\ref{theo:viral}, 
given $k > 0$ and $\delta>0$
we may assume that the 
expelled chain recurrence class satisfies the following:  
The filtrating Markov partition 
for the new chain recurrence class is $k$-aperiodic and 
each connected component has diameter less than $\delta$.
\end{prop}

\begin{proof}
In this proof, we use the notations from the proof of 
Theorem~\ref{theo:viral}.
The last condition can be obtained asking the size of the 
expelled Markov partition very small, which is contained in the 
definition of the $\varepsilon$-expulsibility. 
For the first condition, 
we choose $q$ in such a way that its period is greater than $k$,
see Remark~\ref{r.largeperi}. 
Then we expel the circuit $S$, letting the new filtrating partition
very close to $S'$. 
If it is close enough, then we may assume that 
there are more than $k$ consecutive homoclinic points outside
the neighborhood of the periodic orbit. Thus given a point 
in the filtrating Markov partition, 
\begin{itemize}
\item if it is near the periodic orbit, its period must be larger
than $k$, because to come back to the initial point following the periodic points, we need more than $k$ times
iteration and following the homoclinic orbit does not give 
much shortcut.
\item If it is near the homoclinic point, it needs at least 
$k$ iteration to come back, because the length of the 
homoclinic orbit is larger than $k$.
\end{itemize}
Thus the proof is done.
\end{proof}

Using Proposition~\ref{prop:virain} 
let us conclude the proof of Proposition~\ref{prop:baiin}.
\begin{proof}[Proof of Proposition~\ref{prop:baiin}]
We explain how we obtain $\mathcal{O}_{n+1}$ from 
$\mathcal{O}_n$. We also need to confirm the existence of 
$\mathcal{O}_1$ but it can be done by following the 
induction step.

Given $f \in \mathcal{O}_n$, we have a locally constant 
nested sequence of filtrating regions $\{\mathcal{U}_k(f)\}_{k=1,\ldots,n}$
such that each $\mathcal{U}_k(f)$ satisfies the conclusion. 
Then, we apply Proposition~\ref{prop:virain} to $(f, \mathcal{U}_n(f))$,
which produces a new smaller filtrating region $\mathcal{U}_{n+1}$
which robustly satisfies the conditions in Proposition~\ref{prop:baiin}.
This gives the candidate for $\mathcal{U}_{n+1}(f)$. We need to extend
$\mathcal{U}_{n+1}$ to some open and dense set of $\mathcal{O}_n$.

We proceed as follows: First we choose a countable dense set 
$(g_m) \subset \mathcal{O}_n$. Then, for each $g_m$ we apply
Proposition~\ref{prop:virain}. It gives an open neighborhood 
$\mathcal{W}_m$ of $g_m$ in 
$\mathcal{O}_n$ where we have a filtrating set 
$\mathcal{V}_m$ satisfying the conclusion. 
We define two sequence of open sets $(\mathcal{A}_m)$ 
and $(\mathcal{B}_m)$ as follows:
\begin{eqnarray*}
\mathcal{A}_1 = \mathcal{W}_1, \quad 
\mathcal{B}_1 = \mathcal{W}_1, \quad
\mathcal{A}_{n+1} = \mathcal{W}_{n+1} \setminus 
\overline{\mathcal{B}_n}, \quad
\mathcal{B}_{n+1} = \mathcal{B}_n \cup \mathcal{A}_{n+1}. 
\end{eqnarray*}
Then one can check that 
$\mathcal{O}_{n+1} :=  \cup \mathcal{B}_m$ is an open and dense 
set of $\mathcal{O}_n$. On $\mathcal{A}_m$ we define 
$\mathcal{U}_{n+1}$ to be $\mathcal{V}_m$. This defines a map
on $\mathcal{B}_m$ and consequently on $\mathcal{O}_{n+1}$.

This finishes the construction of $\mathcal{U}_{n+1}$.
\end{proof}


\section{Reduction to local problem}
\label{s.game}

In this section we begin the proof of Theorem~\ref{t.IFS}. 
We provide several arguments which enables us to 
reduce the problem into a local perturbation problem.

\subsection{Overview of the proof of Theorem~\ref{t.IFS}}
\label{ss.back}
Before the proof of Theorem~\ref{t.IFS}, we give 
some overview of the proof of Theorem~\ref{t.IFS}.

In the assumption of Theorem~\ref{t.IFS}, we have an
$\varepsilon$-flexible point $q$ having a large stable manifold 
and a homoclinic point $Q$. 
Then our goal is to construct  
a family of invariant curves  
containing $q$ and $Q$ by a small perturbation
such that in the normal direction 
the dynamics exhibits almost neutral behavior,
in other words, small normal strength. 

In \cite{BS2}, we established a perturbation technique
with which we can expel $q$ using its flexibility.
The main idea (here translated in the IFS language) of the proof is in the same direction: By using the flexibilty of $q$ we can perform a perturbation so that $q$ becomes neutral in one 
diretion, and then we again use the flexibility of $q$ for controlling its local strong stable manifold. In particular, we 
modify the strong stable manifold to be far from the 
components of the image of the IFS in a fundamental 
domain of the center stable manifold. 
These components are consider as
obstructions to avoid. Then turning 
the neutral direction of $q$ to be 
repelling expels the point $q$ from the class. 

In our present problem we have extra difficulty. 
When we expel $q$ we need to keep 
$Q$ as a homoclinic point of $q$. 
Let us explain the strategy for that. 
Once again, we deform $q$ to be neutral by performing a 
perturbation and we want to control its strong stable manifold. 
However, now we furthermore want that this 
curve coincides with its return 
corresponding to the homoclinic point $Q$ (in backward iteration), 
in a neighborhood of $Q$. 
This will be the notion of \emph{pre-solution} (see 
Section~\ref{ss.well-preso}). More precisely, we want that the 
strong stable manifold of $q$ will be disjoint from all the 
components of the images of the IFS (called \emph{the 
obstructions}) but the one 
which corresponds to the homoclinic point $Q$. 
This component will be called the \emph{well}. 
In a well we also define two 
kinds of nested sequence of discs. 
The \emph{transition well} is the 
successive nested images of the IFS along the orbit of the
transition. 
The \emph{periodic well} is the 
nested images along the periodic orbit of $q$. 

For a pre-solution, we require that the (backward) return of 
$W^{ss}_{\loc}(q)$ following the transition itinerary of $Q$ coincides with $W^{ss}_{\loc} (q)$ exactly on one 
of this nested images, avoiding all the obstacles and image of the obstructions out of this well. 
The nested sequence of the wells allow to define the \emph{depth} of the pre-solution.   
It is an integer which indicates in how small region the 
coincidence property holds.
The construction of the invariant curve with this coincidence 
property is one of the big issues 
of this theorem and of this paper. We will do it in 
Section~\ref{s.solution}.

For the time being assume that we can construct such an 
pre-solution and  consider the second point, that is, 
obtaining the almost neutrality of the normal behavior 
of the invariant curve. Recall that we have the neutral behavior of the periodic point $q$, as it was assumed to be flexible. 

We want to spread the neutrality of $q$ 
on the whole invariant curve, but that is not always possible,
as we have no a priori knowledge about the behavior of the 
intermediate dynamics.  
This will be possible only for pre-solution of very large depth, 
that is, for which  the return of $W^{ss}_{\loc}(q)$ coincides with $W^{ss}_{\loc}(q)$ only on 
a very small neighborhood of the orbit of $q$ and $Q$.  
For such pre-solutions the orbits in the maximal invariant set 
in that curve spend most of the time 
near the orbit of $q$, 
and therefore inherit of its neutral behavior. 
For pre-solutions with enough profound depth, 
we will see that the neutralization of the normal dynamics 
indeed happens, thus constructing them will conclude the 
theorem, see Section~\ref{ss.weakinv} and Section~\ref{ss.finfin}. 




The construction of the pre-solution with arbitrarily large depth will be the aim of Section~\ref{s.solution}. An important issue in the construction consists of proving that the \emph{cost} of a pre-solution does not depend on its depth. Let us explain this: Small perturbations in a neighborhood of a flexible point $q$ allow us to get an arbitrary number of successive fundamental domains around $q$ where the diffeomorphism is an homothety. This is the notion of \emph{retarded families} 
that we introduce in Section~\ref{ss.reta}.    
The more one has homothetic fundamental domains,
the more one has freedom for performing perturbations to
modify the the strong stable manifold of $q$. 
The \emph{cost} of a pre-solution is the number of these fundamental domains one needs for getting the chosen invariant curve, and an important point is that this cost remains bounded when the required depth tend to infinity.  
We will discuss more details in Section~\ref{s.solution}.



\subsection{Retarded families}
\label{ss.reta}
Let us start the proof. 

In \cite{BS1} the authors defined the 
notion of an $\varepsilon$-flexible periodic point and 
proved one of its eminent properties: By an 
$\varepsilon$-$C^1$-small perturbation 
in an arbitrarily small 
neighborhood of its orbit, one 
can choose the position of one fundamental 
domain of its stable manifold as one wishes.  
Unfortunately, this result is not enough 
for our purpose and we need to come back
to some essential step of its proof. 

The proof of flexibility property of the 
stable manifold in \cite{BS1} is done 
by using \emph{retardable cocycles}. 
For the proof of Theorem~\ref{t.IFS} 
we need to recall its definition. 
In this paper, we use the notion of retarded diffeomorphisms 
in slightly different way, given as follows (see also Figure~\ref{f.reta}):
\begin{figure}[h]
\begin{center}
\includegraphics[width=10cm]{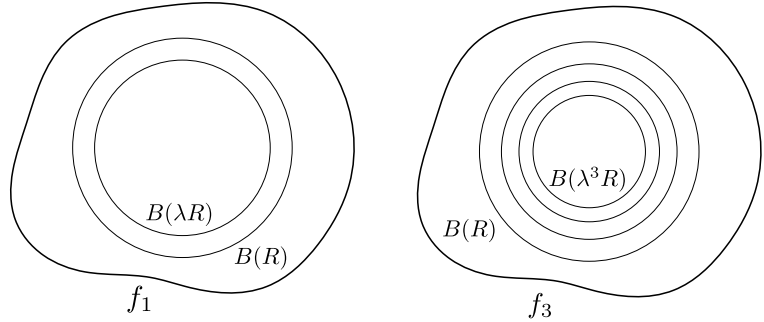}
\caption{An example of retarded family. It is a family 
of diffeomorphisms $\{f_m\}$. Each $f_i$ behaves homothetically
on $B(R) \setminus B(\lambda^{i}R)$.
Roughly speaking, a retarded family is a sequence of
diffeomorphisms obtained by ``inserting'' homothetic regions.}
\label{f.reta}
\end{center}
\end{figure}

\begin{defi}\label{d.reta}
 Let $D$ be a $C^1$-disc in $\mathbb{R}^2$ which contains the 
 origin $\boldsymbol{0}$ in the geometric interior of it.
By $B(R)$ is the closed disc of radius $R$ in $\mathbb{R}^2$ 
centered at $\boldsymbol{0}$.
 A family of diffeomorphisms $\{f_m\}_{m\geq m_0}$ 
(where $m_0$ is some positive integer)
 from $D$ to its image contained in $D$
is called a \emph{retarded family} if it satisfies the following conditions:
\begin{itemize}
 \item There is a radius $R >0 $ such that 
 $B(R) \subset D $ and the maps $f_m$ 
 all coincide on $D \setminus B(R)$,
 \item $\boldsymbol{0}$ is the unique 
 fixed point of $f_m$ for every $m \geq m_0$.
 \item There is $\lambda \in (0,1)$ such that
 for every $m\geq m_0$ 
 the maps $f_m$ coincide with the homothety 
 $H_\lambda=\lambda \mathrm{Id}$ on 
 $B(R)\setminus B(\lambda^{m}R)$.
  The annulus $B(R)\setminus B(\lambda^{m}R)$ is called  \emph{the homothetic region} of $f_m$ and 
 $\lambda$ its \emph{homothetic factor}.
 \item Consider the restriction 
 $f_m|_{ B(\lambda^{m}R)}$. 
 Then, for every $m \geq m_0$, we have
 \[
 f_m|_{ B(\lambda^{m}R)}= 
 H_{\lambda^{m-m_0}}\circ f_{m_0}|_{B(\lambda^{m_0} R)}\circ (H_{\lambda^{m-m_0}})^{-1}.
 \]
 
 The diffeomorphism $f_m|_{B(\lambda^m R)}$ is called \emph{the core dynamics} of $f_m$ 
 and the region
 $B(\lambda^{m}R)$ is called the \emph{core region} of $f_m$.
 
\end{itemize}
\end{defi}

\begin{defi}
A retarded family $\{f_m\}_{m \geq m_0}$ is called \emph{saddle-node}
if 
\begin{itemize}
\item $\boldsymbol{0}$ 
has one positive contracting eigenvalue and the other eigenvalue 
equal to $1$, 
\item there is a neighborhood 
$U_{m_0}$ of $\boldsymbol{0}$ such that 
$f_{m_0}|_{U_{m_0}}$  
 has the form $(x, y) \mapsto (\lambda_0 x, k(y))$ where
 $\lambda_0$ is the contracting eigenvalue and 
 $k(y)$ is a $C^1$ map satisfying $k(0)=0$, $k'(0)=1$
 and topologically attracting in a neighborhood of 
 $0$ (more precisely, there is $\varepsilon>0$ such that 
 $k([-\varepsilon,\varepsilon]) \subset (-\varepsilon, \varepsilon)$ holds).
 \end{itemize} 
 For $f_{m}$, $H_{\lambda^{m-m_0}}(U_{m_0})$ is called the diagonal region of $f_m$.
\end{defi}
%

The arguments of the paper \cite{BS1}
(see Proposition~2.2 and the proof of Theorem~1.1 of \cite{BS1}) 
show that if a diffeomorphism has 
an $\varepsilon$-flexible periodic point
with a large stable manifold, 
then one can produce a saddle-node retarded family 
by giving an $\varepsilon$-small perturbation
along the orbit, keeping the largeness of the 
stable manifold. 
In other words, an $\varepsilon$-flexible point 
can be deformed into a saddle-node point, keeping 
the stable manifold large and inserting homothethic 
fundamental domains as much as one wishes. 
Thus we have the following proposition. 
Recall that for a periodic point $q$ of a Markov IFS, 
$\pi(q)$ denotes its period. 

\begin{prop}\label{p.flexible} Let 
$(\cD=\coprod D_i, F=\{f_{j}\})$ 
be a Markov IFS.
Let $\varepsilon>0$ and a separated 
$\varepsilon$-flexible point with large stable manifold $q$ 
be given.

Given a neighborhood $V$ of $\mathrm{orb}(q)$,
there is an $\varepsilon$-$C^1$-small 
family of perturbations 
$G_m = \{g_{j, m}\}_{m \geq 1}$  along $q$
of $F$ supported in $V$ satisfying the following:
\begin{itemize}
\item For every $m$, $q$ is a periodic point of $G_m$ and 
the orbit is the same as that of $F$.
\item $q$ has a large stable manifold for every
$G_m$ $(m\geq 1)$. 
\item The family of diffeomorphisms 
$\{(G_{m, q})^{\pi(q)}|_{D_q}\}$ (here $(G_{m, q})$
denotes the map of $G_m$ on $\coprod D_{q_i}$, 
see Section~\ref{ss.peri}) 
is a saddle-node retarded family of diffeomorphisms, up to
a coordinate change which is independent of $m$.
\end{itemize}
\end{prop}

The proof of Proposition~\ref{p.flexible} is almost 
immediate from the argument 
of Proposition~2.2 and the proof of Theorem~1.1 of \cite{BS1}.
So we omit the proof.
Based on Proposition~\ref{p.flexible}, we give 
the following definition.

%
%

\begin{defi} 
A family of Markov IFSs  $\{(\cD, F_m)\}_{m\geq m_0}$ is 
said to be \emph{retarded at a periodic point $q$} if the 
following holds:
\begin{itemize}
\item $F_m$ is a perturbation of $F_{m_0}$ along $q$. 
\item $\{(F_{m, q})^{\pi(q)}|_{D_q}\}$ is 
a retarded family on the disc containing $q$.
\item $q$ has a large stable manifold.
\end{itemize}
%
\end{defi}

We prepare one more definition
\begin{defi}\label{d.unibou}
Let $\{(\cD, F_n)\}_{n\geq 1}$ be a family of Markov IFSs. 
We say that $(F_n)$ is uniformly bounded if 
the $C^1$-norm of $(F_n)$ is uniformly bounded. 
\end{defi}

Proposition~\ref{p.flexible} implies the following.
\begin{lemm}\label{l.flexihomoclinic} 
Let $(\cD,F)$ be a Markov
IFS, $\varepsilon>0$ and $q$ be a separated 
$\varepsilon$-flexible 
periodic point with a large stable manifold. 
Assume that there is a $q$-free $u$-homoclinic point 
$Q$ of $q$. 

Then there are $\varepsilon$-$C^1$-small 
perturbations $(F_m)_{m\geq 1}$ of $F$ along $q$ supported in an arbitrarily small neighborhood of the 
periodic orbit of $q$ such that 
\begin{itemize}
 \item $\{(\cD, F_m)\}_{m \geq 1}$ is retarded at $q$
 with the same orbit.
 \item $Q$ is a $q$-free 
 $u$-homoclinic point of $q$ for every $m \geq 1$. 
\end{itemize}
Note that
\begin{itemize} 
\item $(F_m)$ is uniformly bounded since each $(F_m)$
is an $\varepsilon$-$C^1$-small perturbation of 
a single IFS $F$.
\item The $C^0$-distance between $F_m$ and $F$ can be 
chosen arbitrarily small uniformly with respect to $m$,
for the size of the support can be chosen arbitrarily small.
\end{itemize}
\end{lemm}

\subsection{Wells and pre-solutions}
\label{ss.well-preso}

In the setting of Theorem~\ref{t.IFS}, 
we have a separated periodic point $q$ and 
a $q$-free homoclinic point $Q$ of $q$. 
In order to prove the theorem we want 
to reduce the problem to a perturbation problem which 
involves information only 
on the disc $D_q$. Let us formulate it. 

Consider the disc $D_q$ and $\orb(Q)$
(see Section~\ref{sss.per} and \ref{sss.hh} for the definitions). 
Recall that the points $F^{-i}(Q) \in \mathrm{orb}(q)$
for every sufficiently large $i$. We choose the smallest $i$
such that this holds and denote the point by 
$Q_1 \in \mathrm{orb}(q)$. 
Then there is the smallest $\mathfrak{t}>0$ such that 
$F_Q^{\mathfrak{t}}(Q_1) \in \cup D_{q_i}$. 
We set $Q_2 =F_Q^{\mathfrak{t}}(Q_1)$. 
Finally, we choose 
the smallest $\mathfrak{a}\geq0$
such that $F_q^{\mathfrak{a}}(Q_2) \in D_q$.

Then we define the following objects:
\begin{defi}
A \emph{transition well} is a nested sequence of discs 
$(\Xi_n)_{n=1,\ldots, \mkt}$ in $D_q$ and a sequence
of unions of discs 
$\Theta_n \subset \mathrm{Int}( \Xi_n \setminus \Xi_{n+1})$ for $n=1,\ldots, \mkt -1$ 
defined as follows:
\begin{itemize}
\item Definition of $\Xi_n$: Put $\Xi'_n := D_{F^{-n}(Q_2)}$. 
Then
$\Xi_n := 
F_q^{\mathfrak{a}}\circ F_Q^{n}(\Xi'_n)$.
\item Definition of $\Theta_n$:
Put $\Theta'_n := 
(D_{F^{-n}(Q_2)} \cap F(\cD)) $. 
Then 
\[
\Theta_n :=F_q^{\mathfrak{a}}\circ  F_Q^{n}(\Theta'_n) \setminus \Xi_{n+1}.
\]
\end{itemize}
\end{defi}

\begin{defi}\label{d.wells}
A \emph{periodic well} is a nested sequence of discs 
$(T_n)_{n\geq 0}$ in $D_q$ and 
a sequence of unions of discs $(S_n)_{n \geq 0}$
such that $S_n\subset \mathrm{Int}(T_n \setminus T_{n+1})$
defined as follows:
\begin{itemize}
\item Definition of $T_n$:
First put $T'_n := D_{F_q^{-n}(Q_1)}$. Then
\[T_n := F_q^{\mathfrak{a}} \circ F_{Q}^{\mathfrak{t}} \circ F^n_q(T'_n).\]
\item Definition of $S_n$:
First we put $S'_n := D_{F_{q}^{-n}(Q_1)} \cap F(\cD)$. 
Then, put
\[S_n := F_q^{\mathfrak{a}} \circ F_{Q}^{\mathfrak{t}}\circ F^n_q(S'_n) \setminus T_{n+1}.\]
\end{itemize}
\end{defi}
Note that $T_0 = \Xi_{\mkt}$. Thus $(\Xi_i)$ and $(T_i)$ defines a nested 
sequence of discs in the fundamental domain $D_q \setminus F_q^{\pi(q)}(D_q)$. 

The definitions of wells seem complicated, but it can be well
understood as follows: First consider the backward orbit  
$\{F^{-i}(F_q^{\mka}(Q_2))\}_{i \geq 0}$. 
 It initially belongs to periodic discs, then passes 
 transition discs and
 then finally comes back to the periodic discs. 
 The discs 
 $\Xi_n$ and $T_n$ are nothing but the images of these discs 
 in $D_q$ which $F^{-i}(F_q^{\mka}(Q_2))$ passes, 
 and the definition of $\Theta_n$ and $S_n$ are the 
 images of image discs in $\Xi_n$ and $T_n$ excluded 
 $\Xi_{n+1}$ or $T_{n+1}$.

We also define another class of discs in the fundamental domain, see Figure~\ref{f.obst}. 
\begin{defi}
\label{d.objet}
For $i \geq 0$,
let \[
\Lambda_i := F_q^{i}(\partial D_{F^{-i}(q)})
\] 
and 
call it a \emph{stratum}. 
Note that $\Lambda_0 = \partial D_q$
and $\Lambda_{\pi(q)} = \partial (F_q^{\pi(q)}(D_q))$.
We also set 
\[\bar{\Lambda}_i :=  F_q^{i}(D_{F^{-i}(q)})\] 
and 
call it the $i$-th image disc.
In the following, by $i_{\Xi}$ we denote the integer such that 
$\Xi_1 \subset \bar{\Lambda}_{i_{\Xi}} \setminus \bar{\Lambda}_{i_{\Xi}+1}$ holds.

For $0 \leq i \leq \pi(q)-1$, set 
\[
\Delta'_i :=(D_{F^{-i}(q)} \cap F({\cD})) \setminus F_q(D_{F^{-(i+1)}(q)})\]
and put 
\[
\Delta_{i, \ast} = (F_q)^i(\Delta'_i).
\] 
Note
that it is a disjoint union of finitely many discs in the annulus
bounded by $\Lambda_{i}$ and $\Lambda_{i+1}$. Thus we can choose 
a disc $\Delta_i$ contained in the annulus which contains all $\Delta_{i,\ast}$
and disjoint from $\Xi_1$ 
(if the annulus contains $\Xi_1$).  
We fix such $(\Delta_i)$ and
we say that $\Delta_i$ is the \emph{$i$-th obstruction}. 
Also, for $i \geq \pi(q)$ we define $\Delta_i$ 
setting $\Delta_{i+\pi(q)} = F_q^{\pi(q)}(\Delta_{i})$
recursively.

\begin{figure}[h]
\begin{center}
\hspace*{-1.3cm}\includegraphics[width=15cm]{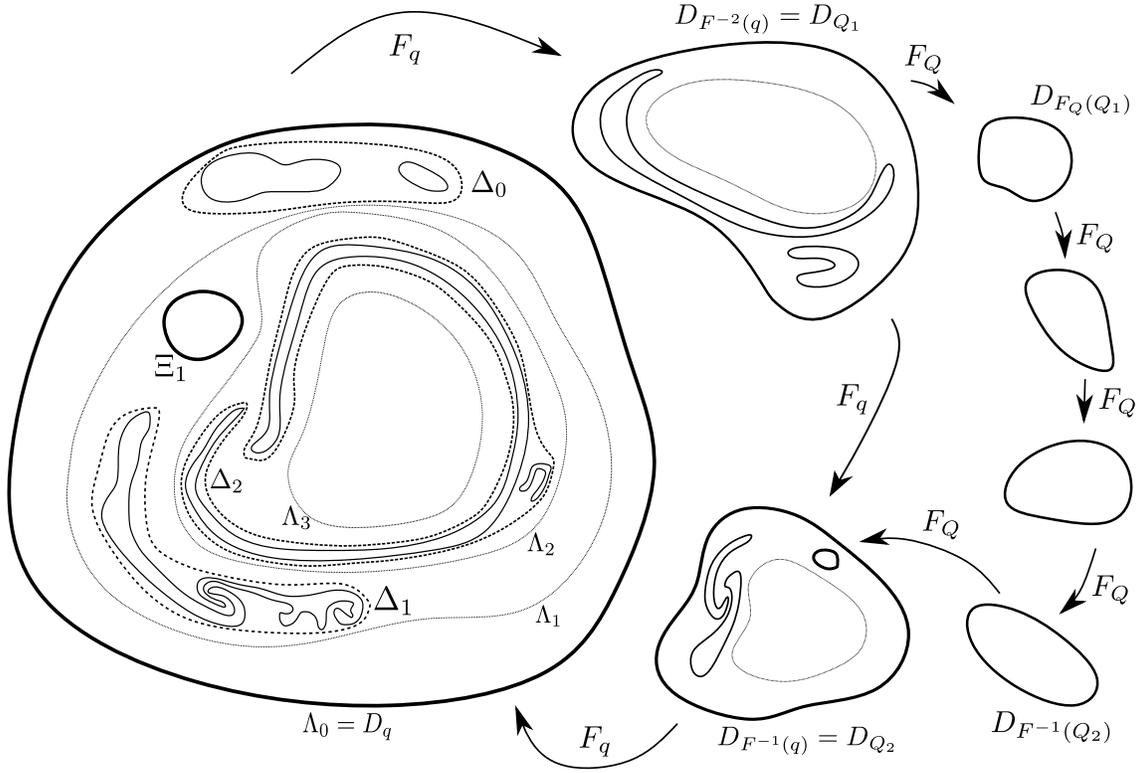}
\caption{A graphical explanation of obstructions and 
wells. $\Delta_i$ is a disc in $\Lambda_i \setminus \Lambda_{i+1}$ which contains all the images of 
image discs except $\Xi_1$. }
\label{f.obst}
\end{center}
\end{figure}

\bigskip

\end{defi}

\begin{rema}\label{r.selfsi}
Let $\mkd \geq 0$ be the smallest integer such that 
$F_{q}^{\mkd}(D_q) \subset D_{Q_1}$ holds. 
Then, recall that the following holds:
\[
T_{\mkd} = F_{q}^\mka \circ F_{Q}^{\mkt} \circ 
F_{q}^{\mkd}(D_q).
\]
Consider the disc $T_{\mkd + j}$ for $j\geq 0$. We have 
\[
T_{\mkd + j}=
F_{q}^\mka \circ F_{Q}^{\mkt} \circ F_{q}^{\mkd}
\big(\bar{\Lambda}_j \big)
\]
Also, $S_{\mkd+k\pi(q)+i}$ is contained in 
\[
\begin{cases}
(F_{q}^\mka \circ F_{Q}^{\mkt} \circ F_{q}^{\mkd})(\Delta_{k\pi(q)+i}), \quad
&\mbox{if $\Xi_1$ is not in 
$\bar{\Lambda}_i \setminus \bar{\Lambda}_{i+1}$},\\
(F_{q}^\mka \circ F_{Q}^{\mkt} \circ F_{q}^{\mkd})(\Delta_{k\pi(q)+i} \cup F_q^{k\pi(q)}(\Xi_1)), \quad
&\mbox{otherwise.} 
\end{cases}
\]
\end{rema}

In order to prove Theorem~\ref{t.IFS}, 
we only need to deal with these information. 
Recall that Theorem~\ref{t.IFS} has two conclusions. 
One is that there is an invariant curve which is univalent in 
some refinement. The other one is that it is $\eta$-weak. 
Let us consider the first part.

\begin{defi}
Let $(\cD, F)$ be an IFS with a separated periodic point $q$ and 
a q-free homoclinic point $Q$ of $q$. 
Let $(\Xi_i)$ and $(T_i)$ be 
the transition well and the periodic well. 
We say that $F$ is a \emph{pre-solution of depth $l$}
if the following holds:
\begin{itemize}
\item[(S1)] $q$ has a strong stable manifold and $Q$ is a $u$-strong homoclinic point of $q$.
\item[(S2)] $W_{\mathrm{loc}}^{ss}(q) \cap \Delta_i = \emptyset$ for $i=0,\ldots, \pi(q)-1$.
\item[(S3)] $W_{\mathrm{loc}}^{ss}(q) \cap \Xi_i$ is a connected $C^1$-curve disjoint from $\Theta_i$
for $i=1,\ldots, \mkt -1$.
\item[(S4)] for $i=0,\ldots, l-1$,
$W_{\mathrm{loc}}^{ss}(q) \cap T_i$ is a connected $C^1$-curve 
disjoint from $S_i$.
\item[(S5)] 
$W_{\mathrm{loc}}^{ss}(q) \cap T_l$ is a connected 
$C^1$-curve satisfying the following:
\[
W_{\mathrm{loc}}^{ss}(q) \cap T_l =  
F_{q}^{\mka} \circ F_Q^{\mkt}(W^{ss}_{\loc}(Q_1) \cap (F_q)^l(D_{F^{-l}(Q_1)})).
\]
\end{itemize}
\end{defi}

\begin{figure}[t]
\begin{center}
\hspace*{-1.3cm}\includegraphics[width=12cm]{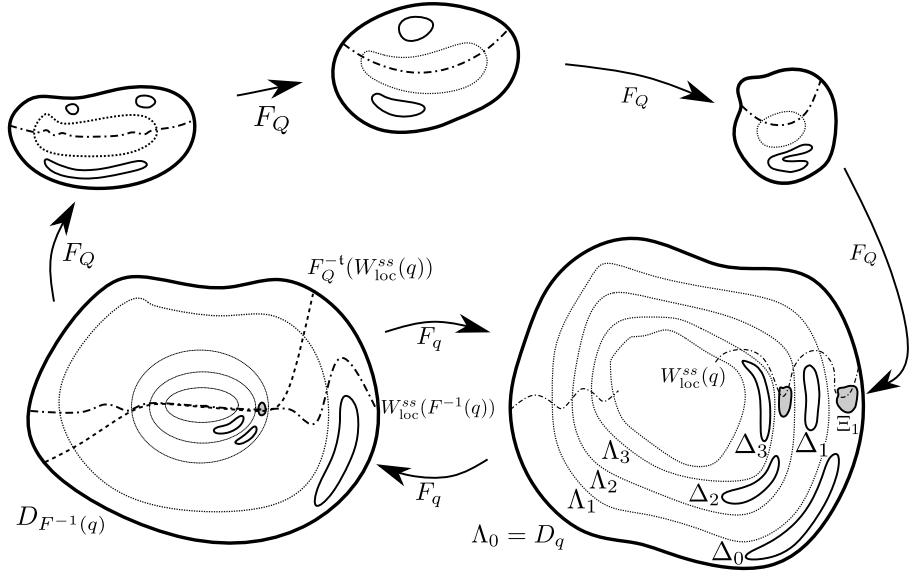}
\caption{A pictorical explanation of pre-solution of depth 
$l$ in $D_q$. In $D_q$, the curve 
$W^{ss}_{\mathrm{loc}}$ avoids every $\Delta_i$ but 
has an non-empty intersection with $\Xi_i$. Thus
in the first fundamental domain 
the backward iteration is well-defined only for 
$W^{ss}_{\mathrm{loc}} (q) \cap \Xi_1$. We can define 
backward images of $W^{ss}_{\mathrm{loc}} (q) \cap \Xi_1$
under $F_Q^{-1}$ until it arrives at the periodic disc.
The backward images avoids all of the image discs 
in the transition discs by condition 
$(S3)$. When it comes back to the periodic discs,
the inverse image of $W^{ss}_{\mathrm{loc}}(q)$ avoids 
all intermediate $\Xi_i$ and $\Delta_i$ but in discs of 
depth $l$ it coincides with the local stable manifold of the 
periodic point (in this picture for the sake of better
visibility this 
coincidence is depicted in $D_{Q_1} = D_{F^{-1}(q)}$). 
Thus by taking the $l$-refinement we can 
take univalent invariant family of curves.
}
\label{f.presol_dq}
\end{center}
\end{figure}

\begin{figure}[h]
\begin{center}
\includegraphics[width=10cm]{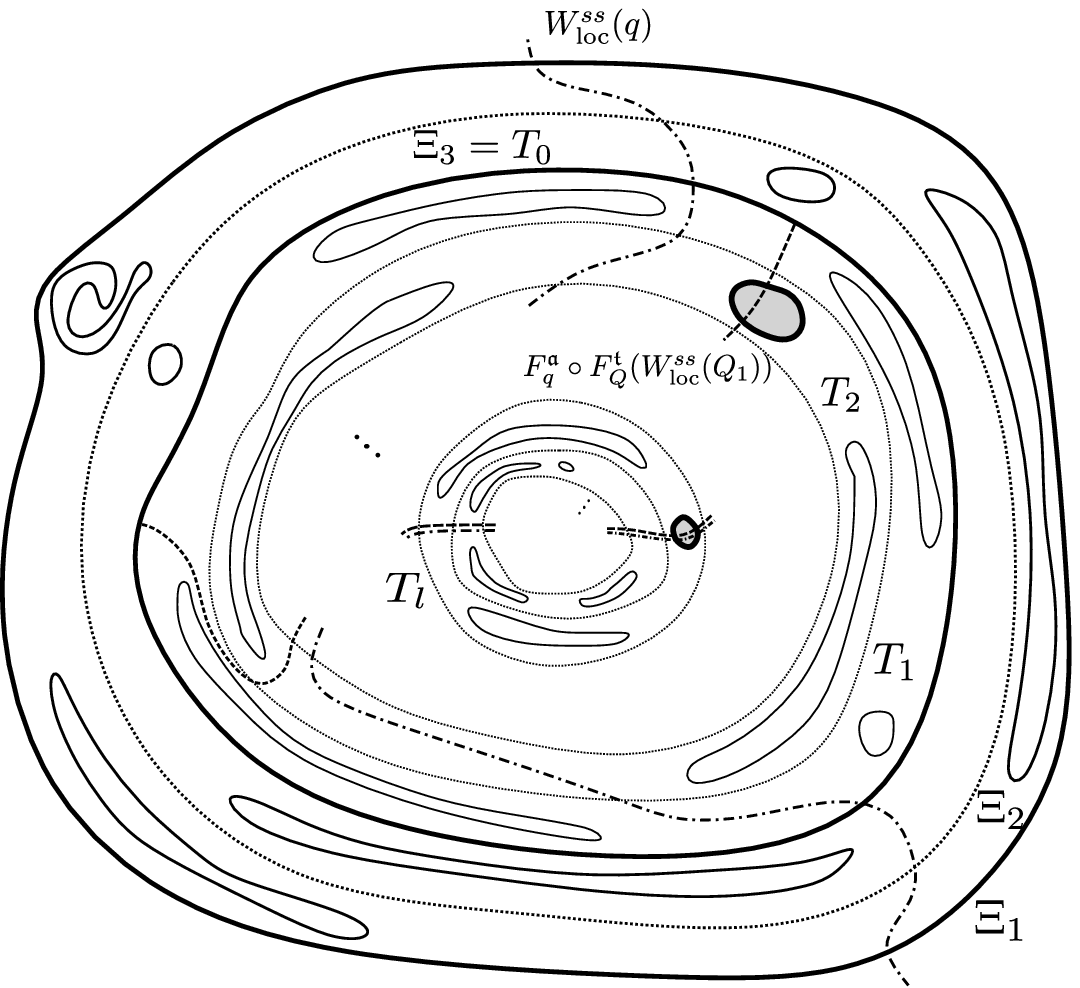}
\caption{A pictorical explanation of pre-solution of depth 
$l$ in $\Xi_1$. There are two curves: 
$W_{\mathrm{loc}}^{ss}(q)$ and  
$F_{q}^{\mka} \circ F_Q^{\mkt}(W^{ss}_{\loc}(Q_1))$.
The former avoids all $\Theta_i$ for $i=1,\ldots, \mkt -1$ 
and $S_i$ for $i<l$ by condition (S3-4).
Inside $T_l$, these two curves coincides. Note that the 
existence of the $u$-homoclinic point implies that 
$F_{q}^{\mka} \circ F_Q^{\mkt}(W^{ss}_{\loc}(Q_1))$ must have non-empty
intersection with $F_{q}^{\mka} \circ F_Q^{\mkt}(\Xi_1)$
which is depicted as a shaded disc.
}
\label{f.presol_obst}
\end{center}
\end{figure}

The following proposition says that
if we have a pre-solution of depth $l$, then 
we can obtain a family of invariant curves which 
is univalent in the $l$-refinement.
\begin{prop}\label{p.preso}
Suppose $(\cD, F)$ has a pre-solution of depth $l$, 
then for the $l$-refinement $(F^l(\cD), \wedge_lF)$ 
there is a family of univalent 
invariant curves $\Gamma_l$ such that it contains $q$
and $Q$ in $(F^l(\cD), \wedge_lF)$. 
\end{prop}

\begin{proof}
First, consider $(F^l(\cD), \wedge_lF)$. 
In the disc $D_q$, 
we have a curve $W^{ss}_{\loc}(q)$. 
Then take its backward images. 

By the conditions $(S1$-$3)$ in the 
definition of the pre-solution, the 
backward images appear on the transition discs 
which contains $\orb(Q)$. 
Note that, if we take the $l$-refinement, 
the corresponding transition discs 
in the refinements are given by the images of transition discs 
or periodic discs in $\cD$ under $(F_q)^l$. 
This, together with $(S4)$, implies that the backward 
images of $W^{ss}_{\loc}(q)$ in the refinement defines 
a connected curve in each corresponding transition discs. 
By the last condition $(S5)$ of the pre-solution 
the collection 
of backward images of $W^{ss}_{\loc}(q_i)$ forms a
family of univalent invariant curves in $(F^l(\cD), \wedge_lF)$. 
\end{proof}

\begin{rema}\label{r.itipre}
The points in $\Gamma_l$ have simple backward 
itineraries. 
If $x \in \Gamma_l$ is in a periodic disc, 
then as $i$ increases the point 
$F^{-i}(x)$ spends some time in the periodic 
discs. Then it arrives at the first fundamental domain of 
$D_q$. 
Then after $\mkt$ backward iterations, passing through
the transition discs, the point comes back to the periodic 
disc. Note that by the definition of pre-solution of depth $l$,
after coming back to the periodic discs this point has at least 
$l$-backward orbit contained in the periodic discs.  
\end{rema}

In the following, we are interested in the following 
special kind of perturbations:
\begin{defi}
Let $\{(\cD, F_n)\}$ be a retarded family with a separated 
periodic orbit $q$ and a separated homoclinic point $Q$.
A perturbation 
$G$ of $F_n$ along the orbit of $q$ is called
\emph{admissible} 
if $G= \phi \circ F_n$ where $\phi$ is a $C^1$-diffeomorphism whose 
support is contained in 
$\bigcup_{i=1}^{n} (F_{n,q}^{\pi(q)})^i(\Xi_1)$. 
\end{defi}

\begin{rema}\label{r.admissible} 
\begin{enumerate}
\item If $(G_m)_{m \geq m_0}$ is a family of admissible perturbation of 
$F_n$, then there exists a neighborhood $W$ of 
$q$
such that $G_{m,q}^{\pi(q)}|_W = F_{n,q}^{\pi(q)}$ for every $m \geq m_0$.
\item If $G$ is an admissible perturbation of $F_n$ then for every $k \geq 0$ we have
$$ 
(F_{n,q}^{\pi(q)})^k(D_q) = G_q^k(D_q).
$$
Also, the shapes of $\Delta_i, \Xi_i, \Theta_i, S_i$ and 
$T_i$ are all 
the same for $F_n$ and $G$.
\end{enumerate}
\end{rema}

\subsection{On the distribution of itineraries for pre-solutions}
\label{ss.weakinv}
For a pre-solution, we have a family 
of univalent invariant curves in some refinements. 
We also want 
to control the normal strength of the invariant curve. 
A priori, there is no information available about 
the information of normal strength. However, if we know 
that the periodic point has a neutral eigenvalue, 
then pre-solutions of the large depths have small strength. 
In Section~\ref{s.solution} we prove that such a construction
is possible for a special type of retarded family 
called prepared family
(see Proposition~\ref{p.solutionprepared}).
We also prove that for every 
saddle-node family, by an arbitrarily small perturbation 
we can make it 
into a prepared one (see Proposition~\ref{p.cleaning}). 
We prove Theorem~\ref{t.IFS} by these propositions. 
see Section~\ref{ss.finfin}.

In this subsection, we prove a result which 
enables us to estimate the distribution of the orbits in 
the invariant curve for a pre-solution of profound depth, 
which will be a fundamental tool for the proof of 
Theorem~\ref{t.IFS}.

\begin{prop}\label{p.frequ}
Let $(\cD, F)$ be a Markov IFS
with a separated periodic point $q$ and 
its $q$-free homoclinic point $Q$ such that $q$ has 
a large stable manifold, where $F$ is a member of some 
retarded family $(F_n)$. 
Given
a neighborhood $W$ of $\orb(q)$ and $r \in (0, 1)$, 
there exists an integer $L_0$ such that the following holds:
Given a pre-solution $G$ of depth $L \geq L_0$  which 
is an admissible perturbation of $F_n$,
consider the Markov IFS $(G^{L}(\cD), \wedge_{L}G)$.  
For every point $x \in \Gamma_L$ 
(see Proposition~\ref{p.preso}),
if $(G)^{-L}(x)$ is defined, then 
one of the following holds:
\begin{itemize}
\item In the interval $[0, L-1]$,
there is a connected interval $H \subset [0, L-1]$ such that
for any $i \in H$,  $G^{-i}(x) \in W \cap G^L(\cD)$ 
and $\# H > rL$.
\item In the interval $[0, L-1]$, 
there are two disjoint connected intervals
$H_1, H_2 \subset [0, L-1]$ such that
for any $i \in H_1 \cup H_2$,  
$G^{-i}(x) \in W \cap G^L(\cD)$ 
and $\# H_1 + \# H_2 > rL$.
\end{itemize}
\end{prop}
\begin{proof}
Let $W$ and $r$ be given.
First, since $q$ has a large stable manifold,
there exists $\ell$ such that  
$(F_q)^{\ell}(\cup D_{q_i}) \subset W$. 
We fix such $\ell$ and   
denote it by $\ell_0$. 
We fix $L_0$ which satisfies 
$(L_0-\ell_0- \mka-\mkt)/L_0 > r$.
Notice that for every $L \geq L_0$, 
$(L-\ell_0-\mka-\mkt)/L > r$ holds.
Let us show that any pre-solution $G$ of depth $L \geq L_0$
which is an admissible perturbation of $F$ 
satisfy the desired condition.

To see this, let us take $x \in \Gamma_L$ such that 
$G^{-L}(x)$ is well-defined. 
We define two sets of integers:
\[
I_1:= \{ 0 \leq i < L \mid 
G^{-i}(x) \in (F_q)^{\ell_0}(\cup D_{q_i}) \}, \quad
I_2 := [0, L-1] \setminus I_1
\]
Let us consider the length of connected intervals 
of $I_1$ and $I_2$.
\begin{itemize}
\item By Remark~\ref{r.itipre} and definition of $\ell_0$, 
the connected intervals of $I_1$ which is bounded by 
the points of $I_2$ has the length at least 
$L-\ell_0+1$. Indeed, let $i_1$ be the first integer of such a connected interval, then by By Remark~\ref{r.itipre}
we know $G^{-i_0}(x) \in (F_q)^{L_0}(\cup D_{q_i})$.
Thus $G^{-i_0-k}(x)$ belongs to $(F_q)^{\ell_0}(\cup D_{q_i})$ for $k=0,\ldots, L-\ell_0$. 
\item By Remark~\ref{r.itipre}, 
the connected intervals of $I_2$ has the 
length no longer than $\ell_0+\mka+\mkt$.  
Note that, together with the definition of $L$, 
this implies that $I_1$ is not empty.
\end{itemize}
Thus we can deduce the following:
\begin{itemize}
\item If $I_2$ is empty, then the conclusion is obvious.
\item If the number of connected intervals 
in $I_2$ is more than one, then there is at least one connected 
interval in $I_1$
which is bounded by the points of $I_2$. 
Let us denote one of them by $H$. Then we have
\[
\#H \geq  L-\ell_0 +1 > Lr.
\]
\item If there is only one connected interval
in $I_2$, then $I_1$ has at most two connected intervals. 
If there is only one, say $H$, 
since the length of the connected interval in $I_2$
is no more than $\ell_0+\mka+\mkt$, we have
\[
\#H\geq  L- (\ell_0+\mka+\mkt) > Lr.
\]
If there are two (and only two) connected components, 
we obtain the conclusion by letting them $H_1$ and $H_2$
and repeating a similar argument.
\end{itemize}
Thus the proof is completed.
\end{proof}

\section{Solution of local problem}
\label{s.solution}
The aim of this section is to complete our construction by 
building  perturbations having the announced  invariant curves. 
Usually, constructions of invariant objects are done by means
of fixed point theorems argument in some infinite dimensional 
setting. Our construction is somehow unconventional: We 
directly propose families of pre-solutions having arbitrarily 
profound depth and we prove that such families are realizable 
by a small perturbation. 

The confirmation of the smallness of the perturbation is the 
main step of the proof. 
For getting our pre-solution by a $C^1$-small perturbation we need enough homothetic fundamental domains, and this changes the deepness of the pre-solution we are looking for. Seemingly, the ``cost'' 
of the perturbation would increase as we require the 
coincidence only on the deeper part: Seemingly 
it would lead us to a vicious circle. 

However, by carefully observing the proof of the 
fragmentation lemma and choosing the curve correctly,
we see that it does not depend on the depth we demand,
if the intermediate dynamics is enough ``clean''. More 
precisely, we will see that, after cleaning, we can choose the 
curves which are graphs with bounded 
derivatives independent of the depth, 
and all these curves have the same cost. 

To follow this strategy, we need to 
examine the geometric information of objects we treat 
and establish a method for cleaning the possible difficulties by 
a small perturbation. 



\vskip 2mm
Let us now explain how our strategy is structured in this section. 
The first step 
is that a saddle-node retarded family can be 
perturbed into a \emph{prepared family} (that is the announced family which is enough ``clean''). 
In Section~\ref{ss.prepare}
we give the precise definition of it and its construction.

Then  we need to show that  any prepared family admits 
a pre-solution of arbitrarily profound depth by small perturbations.
To find a perturbation which makes a given prepared game 
into a pre-solution with large depth, we need to have
an estimate of the $C^1$-size,
which we shall refer as the \emph{cost}
of the perturbation (see Definition~\ref{d.cost}). 
We will prove that there is an upper bound of the cost 
which is independent of the depth.
Our estimation will be obtained as follows: 
First, we prepare a quantitative
version of fragmentation lemma which relates the 
cost of the perturbation and the geometric complexity 
of the curves (Section~\ref{ss.fragme}). 
Then we observe the geometric complexity of the 
curve we need to produce is bounded thanks to the preparedness 
of the family (Section~\ref{ss.choix}). 
The combination of these two techniques enables us 
to conclude Theorem~\ref{t.IFS} (Section~\ref{ss.finale},~\ref{ss.finfin}). 

In this section, by a support of a diffeomorphism
$f:M\to M$ we mean $\mathrm{supp}(f, \mathrm{id})$, 
that is, the closure of the set $\{x \in M \mid x \neq f(x)\}$.

%
%

\subsection{Prepared family}
\label{ss.prepare}
\begin{figure}[h]
\begin{center}
\includegraphics[width=10cm]{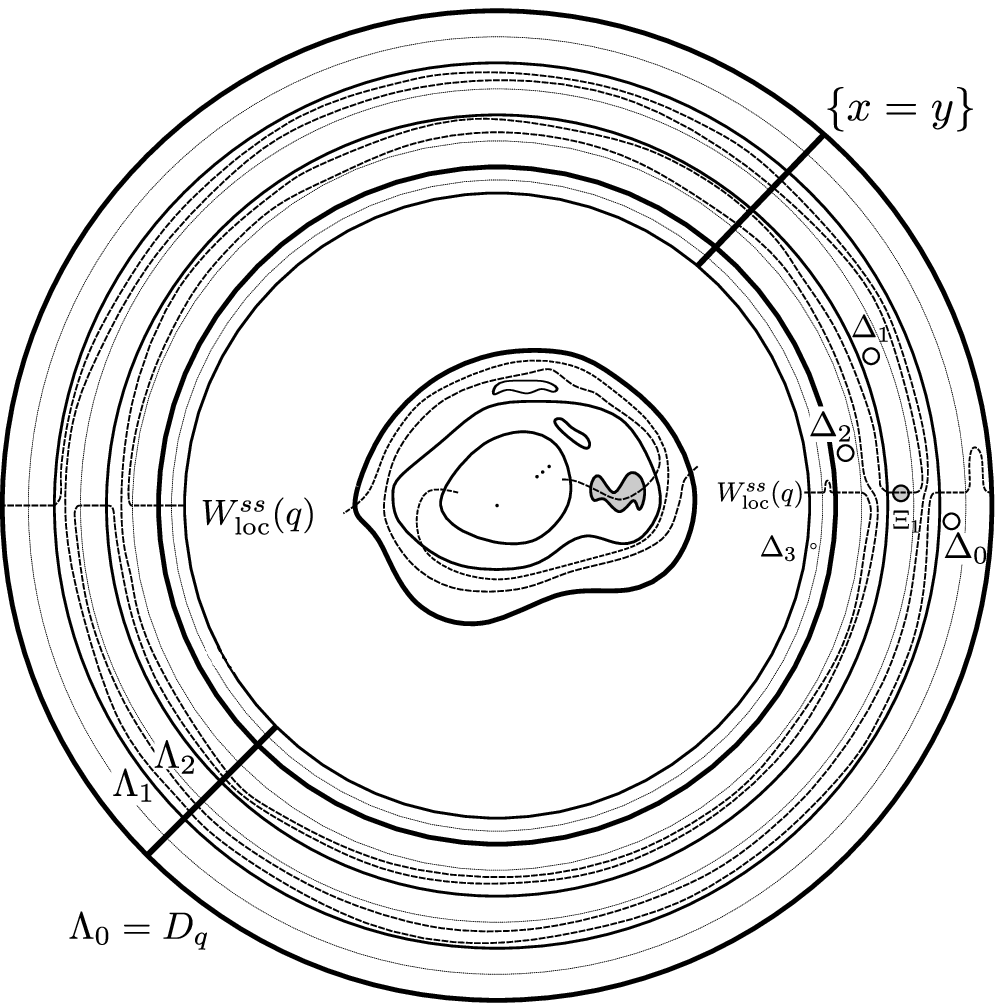}
\caption{A graphical explanation of a prepared family outside
the diagonal region. 
In the a few first fundamental domains the map is a 
contracting homothety. The well $\Xi_1$ and 
the obstructions $\Delta_i$ are all round discs 
(indicated by tiny circles)
and 
$W^{ss}_{\mathrm{loc}}(q)$ coincides with the 
$x$-axis except some region in each annulus 
bounded by $\Lambda_{i}$ and $\Lambda_{i+1}$.}
\label{f.prepared_homo}
\end{center}
\end{figure}

\begin{figure}[h]
\begin{center}
\includegraphics[width=11cm]{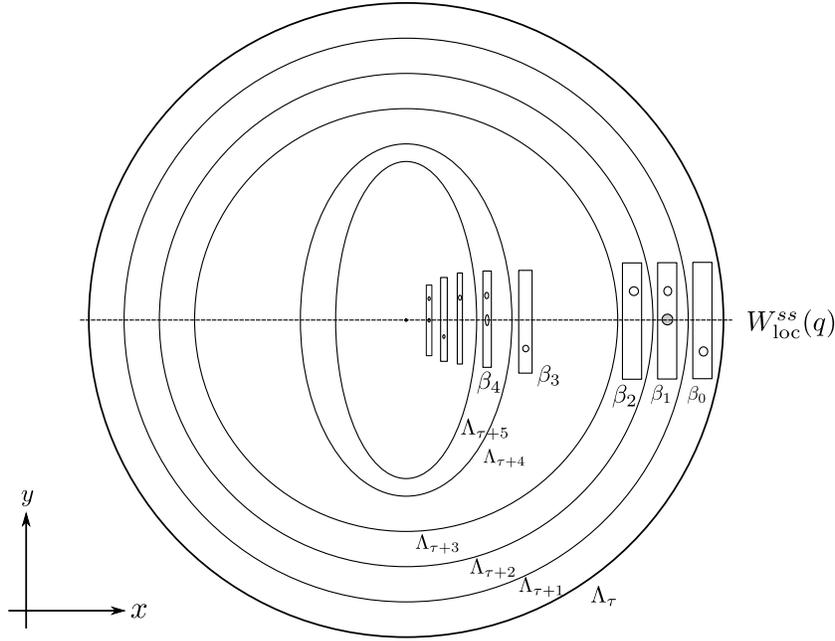}
\caption{A graphical explanation of the prepared 
family near the diagonal region. $\bar{\Lambda}_{\tau +i}$ 
are round discs for $i=0,\ldots, \pi(q)-1$. Each 
annulus bounded by $\bar{\Lambda}_{\tau +i}$ 
and $\bar{\Lambda}_{\tau +i +1}$
contains a rectangle which contains the image 
of the well and the obstruction. 
Note that in the diagonal region 
the $x$-axis coincides with the strong stable manifold 
of $q$.}
\label{f.prepared_center}
\end{center}
\end{figure}

We begin with the definition of the prepared families.
It is a saddle-node retarded family having convenient behavior
of the objects we treat such as the obstructions, images of discs
and the strong stable manifold of the flexible point.
\begin{defi}\label{d.preparation} A saddle-node
family $\{(\cD, F_n)\}_{n \geq 1}$ for a periodic point $q$ and 
a homoclinic point $Q$ is \emph{prepared} 
the following holds
(see Figure~\ref{f.prepared_homo} and~\ref{f.prepared_center}):
\begin{itemize}
\item[(P0)] $D_q$ is a round disc $B(1)$.
 \item[(P1)]  The homothetic region of $(F_{1,q})^{\pi(q)}$ is 
 $B(1) \setminus B(\lambda) =  D_q \setminus B(\lambda)$, 
 where $0<\lambda<1$ is the homothetic factor. Also, 
 $\Lambda_i = B(\lambda_i)$ for $i=0,\ldots, \pi(q)$, where 
 $1=\lambda_0>\cdots>\lambda_{\pi(q)} = \lambda$.
 \item[(P2)] There are $\tau>0$  
 and rectangles $\beta_i \subset D_q$ ($i=0,\ldots,\pi(q)-1$) satisfying the following:
 \begin{itemize} 
\item[(P2-1)] $\bar{\Lambda}_{\tau +i}$ are round 
discs contained in the diagonal region of $F_{1,q}^{\pi(q)}$
whose center is $q$ for $i=0, \ldots, \pi(q)-1$.
\item[(P2-2)] For $i=0, \ldots, \pi(q)-1$,
$\beta_{i}$ is a rectangle in 
the interior of the annulus bounded 
by $\Lambda_{\tau+i}$ and 
$\Lambda_{\tau+i+1}$
such that
its sides are parallel to the coordinate axes (which are eigendirections of $DF_{1,q}^{\pi(q)}|_{q}$) and  
its center is on the positive side of the $x$-axis.
Furthermore, it is disjoint from the line $\{x =y\}$.
 \item[(P2-3)] $\De_{\tau +i}$ is contained 
 in the interior of $\beta_{i}$
 for $i=0, \ldots, \pi(q)-1$.
 \item[(P2-4)] There is an integer $\tau'>0$ such that
 $(F_{1,q}^{\pi(q)})^{\tau'}(\Xi_1)$ is contained in the interior of $\beta_{i}$ for some $i=0, \ldots, \pi(q)-1$..
\end{itemize} 
\item[(P3)] For each $i=0, \ldots, \pi(q)-1$,
there is $\lambda_{i}^{\ast} \in (\lambda_{i+1}, \lambda_i )$
such that the following holds: Set 
$A_i=B(\lambda_{i}) \setminus B(\lambda_{i+1})$
and 
$A'_i=B(\lambda_{i}^{\ast}) \setminus B(\lambda_{i+1})$. 
Then for every $i$ we have
 \begin{itemize}
 \item[(P3-1)]
 The intersection of $W^{ss}(q, F_{1,q}^{\pi(q)})$ with 
 $A_i$ consists of two connected components connecting 
 $\Lambda_i$ and $\Lambda_{i+1}$.
 The intersection of $W^{ss}(q, F_{1,q}^{\pi(q)})$ with 
 $A'_i$ coincides with the $x$-axis. 
 \item[(P3-2)]
 $\Delta_i$ is contained in $A'_i$ and is 
disjoint from the $x$-axis and the line $\{x=y\}$.  
 \item[(P3-3)]  
$\Xi_1$ is contained in some $A'_i$ and is
 a round disc whose center is on the $x$-axis and disjoint from the line $\{x=y\}$.  
\end{itemize}
\end{itemize}
\end{defi}

%

\begin{rema}
\begin{itemize}
\item In condition (P2), $\beta_i$ is defined for $i=0, \ldots, \pi(q)-1$.
For $i\geq \pi(q)$ we define $\beta_i$ setting 
$\beta_{i+\pi(q)} = F_{1,q}^{\pi(q)}(\beta_i)$. Note that they are rectangles 
satisfying similar conditions, but they may touch the 
line $\{x =y\}$. This won't bring any inconvenience 
to our construction.
\item The definition of prepared family is stated as the condition of $F_1$. 
Note that is $F_1$ satisfies the conditions, then $F_n$ satisfies the corresponding 
conditions. More precisely,
\begin{itemize}
\item Condition (P1) holds replacing $B(1) \setminus B(\lambda)$ with
 $B(1) \setminus B(\lambda ^{n})$. 
\item Condition (P2) holds in the corresponding diagonal region of $F_n$ by letting
replacing $\tau$ with $\tau + n-1$.  
\item Condition (P3) holds for $F_n$ as it is. 
\end{itemize}
\end{itemize}
\end{rema}

The following proposition says that
any saddle-node retarded
family admits an arbitrarily small perturbation 
such that 
the perturbed family is 
$C^1$-conjugated to a prepared game. 
\begin{prop}\label{p.cleaning} 
Let  $\{(\cD, F_n)\}_{n \geq 1}$ be a saddle-node family for 
an $\varepsilon$-flexible, 
separated periodic point $q$ and its
$q$-free homoclinic point $Q$.
Then, for every $\delta >0$
there exist a $C^1$-coordinate change
$\varphi$ of $\cD$ (see Definition~\ref{d.coord}) independent 
of $n$
and a saddle-node 
family $(G_n = \{g_{n,j}\})_{n \geq 1}$
such that the following holds:
\begin{itemize}
\item There is $n_1 \geq 1$ such that for every $n \geq 1 $, 
$G_{n}$ is a perturbation of $F_{n+n_1}$ 
along $q$ and the two IFSs
$F_{n+n_1}$ and 
$G_{n}$ are $\delta$-$C^1$ close.
\item $\{G_{n}\}_{n\geq 1}$ 
is a prepared family up to the coordinate
change $\varphi$, that is,
$(\varphi(\cD) ,\{\varphi \circ g_{n,j} \circ \varphi^{-1} \})_{n\geq 1}$ 
is a prepared family for $\varphi(q)$ and $\varphi(Q)$.
\end{itemize}
Note that if $(F_n)$ is uniformly bounded (see 
Definition~\ref{d.unibou}) then the same holds for 
$(G_n)$, since each $G_n$ is a perturbation of one 
of $F_n$ whose size is uniformly bounded.
\end{prop}

The proof of Proposition~\ref{p.cleaning} consists of several steps.
Since the proof is lengthy, we give the outline 
of the proof before going into the detail.

First, we consider a 
deformation of retarded family $(F_n)$ 
with which we can achieve the desired conclusion.
Then, we investigate the ``cost'' of the deformation, 
where the word ``cost'' means the number of diffeomorphisms 
which are $\delta$-close to the identity whose composition 
realizes above perturbations.
Then we use the retardability of the family to obtain 
homothetic regions where we realize above $\delta$-small 
diffeomorphisms. 
These perturbations produce $(G_n)$ for which 
the properties of prepared families hold inside 
homothetic regions. Then, since the number of 
fundamental domains outside the homothetic region
is finite, one can complete the proof just by 
taking a coordinate change outside 
the homothetic region. 

\begin{proof}
Let a saddle-node family $\{(\cD, F_n)\}$ 
for $q$ and $Q$ be given. 
Let $B(R) \setminus B(\lambda^{n}R)$ be 
the homothetic region of $F_{n,q}^{\pi(q)}$.

{\bf Step 1: First preparation.}

Since the stable manifold of $q$ is large, 
$(F_{n, q}^{\pi(q)})^j(\Lambda_i)$ converges 
to $\{q\}$ as $j \to \infty$ for every $n$ and $i=0,\ldots,\pi(q)-1$. 
The same holds for $(F_{n, q}^{\pi(q)})^j(\Delta_i)$ ($i=0,\ldots,\pi(q)-1$) and 
$(F_{n, q}^{\pi(q)})^i(\Xi_1)$. 
Thus, we can take the projections of 
$\Lambda_i$, $\Delta_i$, ($i=0,\ldots, \pi(q)-1$) and $\Xi_1$ to the
orbit space  
(the quotient space obtained by identifying the 
points in the same orbit) of the punctured 
disc $D_q \setminus \{q\}$ which 
is diffeomorphic to the $2$-torus $\TT^2$ 
(for the detail of the orbit space, 
see Section 2 and Section 3 of \cite{BS1}). 
Notice that 
the projections of these objects are independent of the choice 
of $n$.
We project the two branches of 
the strong stable manifold $W^{ss}(q)$, 
$\Lambda_i$, $\Delta_i$  and $\Xi_1$ to $\TT^2$. 
We denote them by $\tilde{\sigma}_{1}$, $\tilde{\sigma}_{2}$, 
$\tilde{\Lambda}_i$, $\tilde{\Delta}_i$ and 
$\tilde{\Xi}_1$ respectively. In the following,  
given a object in $D_q \setminus \{q\}$, we denote its 
projection to $\TT^2$ by putting a tilde. 

\medskip

{\bf Step 2: Preparation from outside.}

We fix constants $1 = \lambda_0>\cdots>\lambda_{\pi(q)-1}>\lambda$
and consider $B(\lambda_i R)$. Also, we fix 
$\{\lambda_{i}^{\ast}\}_{i=0,\ldots, \pi(q)-1}$ 
satisfying $\lambda_{i}>\lambda_{i}^{\ast}>\lambda_{i+1}$.
We denote the boundary of $B(\lambda_i R)$ by $C_i$
and its projection to $\TT^2$ by $\tilde{C}_i$.
Note that $C_i$ are in the same homotopy 
class as $\Lambda_i$ in $D_q \setminus \{q\}$.
Thus we can find an $C^1$-ambient isotopy 
which maps 
$\cup \tilde{\Lambda}_i$ to $\cup \tilde{C}_i$ 
in $\TT^2$ (i.e., a $C^1$-diffeomorphism 
isotopic to the identity such that it
maps 
$\cup \tilde{\Lambda}_i$ to $\cup \tilde{C}_i$).
We denote it by $\tilde{X}_1: \TT^2 \to \TT^2$.

We choose round discs 
$D_i \subset 
\mathrm{Int}(B(\lambda_i R) \setminus B(\lambda_{i+1} R))$
such that
$D_i$ is disjoint from the $x$-axis and 
the line $\{x=y\}$ for $i=0,\ldots, \pi(q)-1$
and $D_{\Xi}$ is a round disc contained in 
$\mathrm{Int}(B(\lambda_{i_\Xi} R) \setminus B(\lambda_{i_\Xi+1} R))$ disjoint from $\{x=y\}$ and 
$D_{i_\Xi}$ such that
its center is on the $x$-axis 
(see Definition~\ref{d.objet} for the definition of $i_{\Xi}$). 
By deforming $\tilde{X}_1$, we can obtain another ambient isotopy 
$\tilde{X}_2$ which satisfies the condition of 
$\tilde{X}_1$, 
$\tilde{X}_2(\tilde{\Delta}_i)=\tilde{D}_i$ and 
$\tilde{X}_2(\tilde{\Xi}_1)=\tilde{D}_\Xi$.

\medskip

{\bf Step 3: Preparation from inside.}

Now we consider the information in the diagonal region.
We fix concentric round circles $E_i$ 
($i=0,\ldots\pi(q)-1$) in $D_q$ contained in the 
diagonal region of $F_{1,q}^{\pi(q)}$ 
such that 
\begin{itemize}
\item the center of $E_i$ is $q$.
\item $E_{i+1}$ is contained in $E_i$
for $i=0,\ldots\pi(q)-1$. 
\item $E_{\pi(q)-1}$ contains
$E_{\pi(q)} :=F_{1,q}^{\pi(q)}(E_0)$ in its interior.
\end{itemize} 
We denote their projections to $\TT^2$ by 
$\tilde{E}_i \subset \TT^2$.
Consider $\tilde{E}_i$ and $\tilde{C}_i$ in $\TT^2$
($i=0,\ldots, \pi(q)-1$).
Since they are in the same homotopy class in 
$D_q \setminus \{q\}$, we see that they are ambient 
isotopic in $\TT^2$. 
We denote such an ambient isotopy 
by $\tilde{Y}_1$.

Let us consider the position of the strong stable manifold 
$\tilde{\sigma}_j$ ($j=1, 2$) in the orbit space.
First, in each annulus $A_i$ ($i=0,\ldots,\pi(q)-1$), 
let $r_{1,i}$ (resp. $r_{2, i}$) 
be the intersection with $A_i$ and 
the half $x$-axis in the positive (resp. negative) side. 
Similarly, using $A'_i$ ($i=0,\ldots,\pi(q)-1$), 
we define $r'_{1,i}$ and $r'_{2, i}$ in a similar way.

Then, we choose two family of curves 
$\gamma_{i, j, k}$ ($j=1, 2$, $k \in \mathbb{Z}$) 
in $\TT^2$ such that 
the following holds:
\begin{itemize}
\item In a small neighborhood of $r'_{i, j}$, 
$\gamma_{i, j, k}$ coincides with $r_{i,j}$.
\item In $A_i \setminus A'_i $, $\gamma_{i, j, k}$ winds
$k$-times in the counter-clockwise direction (if $k$ is negative then it winds $k$-times in the clockwise direction).
\item $\gamma_{i, 1, k}$ and $\gamma_{i, 2, k}$  
are disjoint.
\end{itemize}

Recall that $E_i$ are round circles in the diagonal region
and $E_i$ and $E_{i+1}$ bounds an annulus. 
Consider the intersection of the annulus and 
$\sigma_j$. It is a union of two disjoint curves.
We denote each connected component 
by $\sigma_{j, i}$ ($j=1, 2$). 
Then, in $\TT^2$, 
$\tilde{\sigma}_{i, j}$ is a curve 
which 
connects $\tilde{E}_i$ and $\tilde{E}_{i+1}$. 

Now, we consider a diffeomorphism 
$\tilde{Y}_2$ which is isotopic to $\tilde{Y}_1$ satisfying 
the following:
\begin{itemize}
\item $\tilde{Y}_2$ is still an ambient isotopy 
between $\cup \tilde{C}_i$ and $\cup \tilde{E}_i$:
$\tilde{Y}_2(\cup \tilde{C}_i) = \cup\tilde{E}_i$.
\item $\tilde{\sigma}_{i,j} \cap \tilde{A}'_i$ 
coincides with 
$\tilde{Y}_2(\tilde{r}'_{i,j})$, and 
\item $\tilde{\sigma}_{i,j}$ coincides with $\tilde{Y}_2(\tilde{r}'_{i,j})$ 
in a small neighborhood of $\tilde{Y}_2(\tilde{r}_{i,j})$.
\end{itemize}
Then we see that for each $i$,
$\cup_{j}\tilde{\sigma}_{i,j}$ is 
ambient isotopic to 
$\cup \tilde{Y}_2(\tilde{\gamma}_{i, j,k_i})$
in the annulus bounded 
by $\tilde{E}_i$ and $\tilde{E}_{i+1}$ 
for some $k_i$. Thus, we can take $\tilde{Y}_3$ which 
is isotopic to $\tilde{Y}_2$ satisfying all the conditions 
of $\tilde{Y}_2$ and 
$\tilde{\sigma}_{i,j} = \tilde{Y}_3(\tilde{\gamma}_{i, j,k_i})$ for every $i$ and $j$.

Finally, we choose $\beta_i$.
Consider the annulus bounded by 
$E_i$ and $E_{i+1}$.
We choose a subset $\beta_i$ 
which satisfies the following:
It is a rectangle
whose edges are parallel to two coordinate axes and whose center is on the $x$-axis.
Notice that $\mathrm{Int}(\beta_i) \cup W^{ss}(q)$
coincides with the $x$-axis. 
Now we perform final modification to $\tilde{Y}_3$:
We take $\tilde{Y}_4$ which is isotopic to $\tilde{Y}_3$,
satisfying all the conditions of $\tilde{Y_3}$ and furthermore 
it satisfies
\begin{itemize}
\item $\tilde{Y}_4(\tilde{D}_i)$ is a disc contained in 
$\mathrm{Int}(\tilde{\beta}_i)$.  
\item $\tilde{Y}_4(\tilde{D}_{\Xi})$ is a disc contained in 
$\mathrm{Int}(\tilde{\beta}_i)$ for some $i$. 
\end{itemize}

{\bf Step 4: Estimation of the cost of the deformation.}

We have finished the preparation of the perturbation. Now
we give perturbations to the family $(F_n)$.
In the following, the maps which will be perturbed is just 
$F_{n,F^{-1}(q)}$.

As is in the argument of \cite{BS1} 
(see the proof of Lemma~3.1 of \cite{BS1}), 
we realize $\tilde{X}_2$ and $\tilde{Y}_4$ as follows:
First, by applying fragmentation lemma (see 
Theorem~\ref{t.fragme} in the next section), 
we take $C^1$-diffeomorphisms 
$\{\tilde{\chi}_i\}_{i=1,\ldots, K}$ and 
$\{\tilde{\upsilon}_i\}_{i=1,\ldots, L}$ of 
$\TT^2$
which are $\delta$-$C^1$-close to the identity
and are supported on small discs 
such that the following equalities hold.
\[
\tilde{X}_2 = \tilde{\chi}_K \circ \cdots \circ \tilde{\chi}_1, \quad
\tilde{Y}_4 = \tilde{\upsilon}_L \circ \cdots \circ \tilde{\upsilon}_1
\]
Now, consider the family $\{(\cD, F_n)\}$ where 
$n = 3K+k+3L$ and $k\geq 1$. 
We will prove the theorem setting $n_1 = 3K + 3L$
(for $n_1$ see the statement of Proposition~\ref{p.cleaning}).

{\bf Step 5: Realization of deformations.}

In each three
consecutive fundamental domains in the homothetic region 
$B(\lambda^{3i-3}R) \setminus B(\lambda^{3i} R)$, 
we take the 
lift of $\tilde{\chi}_i$ and denote it by $\chi_i$. 
We take $\tilde{\chi}_i$ in such a way that their 
supports are so small that 
we can take the lift of them supported in 
$B(\lambda^{3i-3}R) \setminus B(\lambda^{3i} R)$. 

Now consider the following diffeomorphism $X$ of $D_q$
(whose support is contained in $B(R) \setminus B(\lambda^{3K} R)$):
\[
X = \chi_K \circ \cdots \circ \chi_1. 
\]
Since the conjugation by a homothetic transformation 
does not affect the $C^1$-distance, $X$ is also 
$\delta$-$C^1$-close to the identity map. 
We define the new IFS $(\cD, F'_n)$ by composing 
$X$ to $F_{n, F^{-1}(q)}$ and keep the other maps intact.

Now, there exist $K$, $K_0$ and $n_1$ such that the following holds
for every $n = n_1 +k$:
\begin{itemize}
\item $(F_{n,q}^{\prime\,\pi(q)})^{K_0}(\Lambda_i) 
= B(\lambda^{K}\lambda_i R)$ for 
$i=0,\ldots, \pi(q)-1$.
\item $F_{n,q}^{\prime\,\pi(q)}$ is 
a homothety of homothetic factor $\lambda$
on $B(\lambda^{3K}R) \setminus B(\lambda^{3K+k+3L}R)$.
\item $(F_{n,q}^{\prime\,\pi(q)})^{K_0}(\Delta_i)$ is 
a round disc contained in 
$(F_{n,q}^{\prime\,\pi(q)})^{K_0}
(\bar{\Lambda}_i \setminus \bar{\Lambda}_{i+1})$
disjoint from the $x$-axis and the line $\{x=y\}$.
\item $(F_{n,q}^{\prime\,\pi(q)})^{K_0}(\Xi_1)$ is 
a round disc contained in 
$(F_{n,q}^{\prime\,\pi(q)})^{K_0}
(\bar{\Lambda}_{i_\Xi} \setminus \bar{\Lambda}_{i_\Xi+1})$
whose center is on the $x$-axis and disjoint from 
the line $\{x=y\}$.
\end{itemize}

\medskip

We perform another perturbation.
For $n=n_1 +k$, take the lift of $\tilde{\upsilon}_i$ and compose them 
on
$B(\lambda^{3K+k}R) \setminus B(\lambda^{3K+k+3L}R)$
and define $Y$ in the similar way. 
That is, first, we take the lift of $\tilde{\upsilon}_i$ on 
$B(\lambda^{3K+k+3(i-1)}R) \setminus 
B(\lambda^{3K+k+3i}R)$
and denote them by $\upsilon_i$ 
(remark that $\upsilon_i$ does depend on $k$, while 
$\chi_i$ does not).
Then consider the following diffeomorphism of $D_q$:
\[
Y = \upsilon_L \circ \cdots \circ \upsilon_1
\]
Then for $F'_n$ we compose $Y$ to
$F'_{n, F^{-1}(q)}$ and keep the other maps intact. 
We denote this IFS by $G_k$.

We can check that it
satisfies, in addition to the previous four
conditions, the following ones:
\begin{itemize}
\item $W^{ss}(q, F^{\pi(q)}_{n,q})$ coincides 
with the $x$-axis on 
$B(\lambda^{3K}\lambda^{\ast}_iR)  \setminus
B(\lambda^{3K}\lambda_{i+1}R)$
for every $i=0,\ldots,\pi(q)-1$. 
Note that the same holds for 
$B(\lambda^{3K+j}\lambda^{\ast}_iR) \setminus
B(\lambda^{3K+j}\lambda_{i+1}R)$
for $j=0,\ldots, k-1$.
\item $(G_{k,q}^{\pi(q)})$ satisfies the condition (P2)
in Definition~\ref{d.preparation}.
\end{itemize}

Thus $(G_{k,q}^{\pi(q)})_{k \geq 1}$ is a retarded family 
for $q$ and $Q$
with a homothetic region 
$B(\lambda^{3K}R) \setminus B(\lambda^{3K+k}R)$
 such that every 
$G_k$ is $\delta$-$C^1$-close to $F_{n_1+k}$.

\medskip

{\bf Final step: Taking conjugacy.}

Now, let us take a coordinate change between $D_q $ 
and $B(1)$ under  
which the family $(G_k)$ satisfies the condition (P1)
(we change coordinates only on $D_q$).
First, let $R_G = \lambda^{3K}R$ be the radius of a ball such that 
$G_{k,q}^{\pi(q)}$ 
is a homothety on $B(R_G) \setminus B(\lambda^{k}R_G)$. Recall that $K_0 \geq 0$ satisfies 
$(G_{k,q}^{\pi(q)})^{K_0}(D_q) = B(R_G)$. 
Then we define a family of diffeomorphisms 
$\{h_{k}$\} defined on the disc 
$B(\lambda^{-K_0}R_g)$ as follows:
\begin{itemize}
\item $h_k = G_{k,q}^{\pi(q)}$ inside 
$(G_{k,q}^{\pi(q)})^{K_0}(D_q) 
=B(R_g)$.
\item Outside $(G_{k,q}^{\pi(q)})^{K_0}(D_q)$, 
$h_k$ is a homothety of homothetic factor $\lambda$. 
\end{itemize}
Now we can define a conjugacy 
$\rho:D_q \to B(\lambda^{-K_0}R_G) $ between 
$G_{k,q}^{\pi(q)}$
and $h_k$ as follows: Set $\rho$ to be the identity map 
inside $(G_{k,q}^{\pi(q)})^{K_0}(D_q)$. 
Outside $(G_{k,q}^{\pi(q)})^{K_0}(D_q)$, 
we extend $\rho$ subject to the formula 
$\rho = h_{k}^{-1}\circ \rho \circ G_{k,q}^{\pi(q)}$.
Notice that $\rho$ can be extended to a conjugacy 
$\rho : D_q \to B(\lambda^{-K_0}R_G)$
since $(G_{k,q}^{\pi(q)})^{K_0}(D_q) = B(R_G)$.

Finally, by conjugating the family $\{h_k\}$ by 
the homothety with homothetic factor $\lambda^{K_0}$, 
we obtain the desired family. Notice that the conjugation by the 
homothety keeps $x$-axis and the roundness of discs.
Thus we see that this conjugated family satisfies (P1). 
Using (P1) and the conditions on $(G_k)$, we can 
also check that the resulted family also satisfies 
(P2) and (P3).  This completes the proof.
\end{proof}

\begin{rema}
If we construct a prepared family by Proposition~\ref{p.cleaning}, then the obtained 
prepared family is automatically bounded, for 
$(F_n)$ are uniformly bounded
each $G_k$ is $\delta$-close to $F_{k+n_1}$.
\end{rema}

In the rest of this section we will prove the following:
\begin{prop}\label{p.solutionprepared} 
Let $\{(\cD, F_n)\}_{n \geq 1}$ be a prepared family for
a periodic point $q$ and its $u$-homoclinic point $Q$ such that $(F_n)$ is uniformly bounded.  
For every $\varepsilon_2 >0$, there exist $n_0, m_0\geq 1$ and 
a neighborhood $W_{n_0}$ of $q$ such that
there is a family of adapted $\varepsilon_2$-perturbations 
$(G_{m})_{m\geq m_0}$ of $F_{n_0}$ such that 
$G_{m}$ is a pre-solution of depth $m$ and
$(G_{m,q})^{\pi(q)}|_{W_{n_0}} 
= (F_{n_0,q})^{\pi(q)}|_{W_{n_0}}$ holds for 
every $m \geq m_0$.
\end{prop}
Notice that if we prove Proposition~\ref{p.solutionprepared},
then with Proposition~\ref{p.cleaning} we can conclude 
Theorem~\ref{t.IFS}.
In the end of this section we will prove Proposition~\ref{p.solutionprepared}
using Proposition~\ref{p.monoprix}, 
which will be proved in the next two sections. 

\subsection{Fragmentation lemma and the cost of a curve}
\label{ss.fragme}

In this subsection, we give an important ingredient of the 
proof of Proposition~\ref{p.solutionprepared}. 
First, let us recall the statement of classical fragmentation lemma.

\begin{theo}[Fragmentation lemma]\label{t.fragme}
Given any smooth closed Riemannian manifold $M$, 
any diffeomorphism $f\colon M\to M$ isotopic to the 
identity map and any $\varepsilon>0$, there is 
a sequence $\{ \varphi_i \}_{i=1,\ldots,k}$ of diffeomorphisms of $M$ with the following properties:
\begin{enumerate}
\item For every $i$, the $C^1$-distance between $\varphi_i$ and 
the identity map $\mathrm{Id}_{M}$ is less than $\varepsilon$.
\item For every $i$, $\varphi_i$ coincides with the identity map outside a disc of radius $\varepsilon$.
\item $f=\varphi_k\circ\cdots \circ \varphi_1$. 
\end{enumerate}
\end{theo}

Theorem~\ref{t.fragme} enables us to decompose a given
diffeomorphism into a composition of diffeomorphisms 
whose $C^1$-distance from the identity is arbitrarily small,
while it does not give any information about the number of 
diffeomorphisms needed. 
In order to prove Proposition~\ref{p.solutionprepared}, 
we need to establish the upper bound of it. To clarify the 
meaning  of the upper bound, we introduce a definition.

\begin{defi}\label{d.cost} Let $\gamma_1, \gamma_2$ 
be $C^1$-curves in a disc
$D \subset \mathbb{R}^2$ transverse to $\partial D$.
We assume that $\gamma_1$ and $\gamma_2$ 
coincide near $\partial D$.
Let $\eta>0$. The \emph{$\eta$-cost from 
$\gamma_1$ to $\gamma_2$}, 
denoted by $c_\eta(\gamma_1, \gamma_2)\in\NN$, 
is the minimum integer $n$ such that 
there are diffeomorphisms $\varphi_1,\dots,\varphi_n$ of 
$D$
satisfying the following:
\begin{itemize} 
\item $\varphi_i$ is supported in a disc of radius $\eta$ contained in 
$D \setminus \partial D$ for every $i$. 
\item $\varphi_i$ is $\eta$-$C^1$-close to the identity map 
for every $i$, and 
\item $\varphi_n\circ\cdots\circ\varphi_1(\gamma_1)=\gamma_2$.
 \end{itemize}
\end{defi}

We will prove the following:
\begin{prop}\label{p.monoprix} 
Given a real number $\eta>0$ and a prepared family 
$\{(\cD, F_n)\}_{n\geq 1}$ for 
a periodic point $q$ and its $u$-homoclinic point $Q$,
there exist $c = c_{\eta}\in \NN$, $m_2 = m_{2, c}\in \NN$ and 
a curve $\gamma_0 \subset \Xi_1$ such that
for every $n\geq 1$ and $m \geq m_2 + (n-1)\pi(q)$ there is a 
curve $\gamma_{n, m}\subset \Xi_1$ 
which coincides with $\gamma_0$ near 
 $\partial \Xi_1$ and satisfying 
the following (see Definition~\ref{d.wells} for the 
definitions of $T_i$, $S_i$ and $\Theta_i$):
\begin{itemize}
 \item $\gamma_0 = \{x\mbox{-\emph{axis}}\} 
 \cap  \Xi_1
 =W^{ss}_{\mathrm{loc}}(q, F_{n,q}^{\pi(q)}) \cap \Xi_1.$
 \item $\gamma_{n,m}\cap T_m 
 = F_{n,q}^{\mka} \circ F_{n,Q}^{\mkt} \circ F_{n,q}^{m}\left(W_{\mathrm{loc}}^{ss}(F^{-m}(Q_1))\right)$ 
 (see Section~\ref{ss.well-preso} for the definitions of $Q_1$,
 $\mkt$ and $\mka$).
 \item $\gamma_{n,m} \cap \Xi_i$ is a connected $C^1$-curve and 
 $\gamma_{n,m}\cap \Theta_i=\emptyset$
 for $i=1,\ldots,\mkt-1$.
 \item $\gamma_{n,m} \cap T_i$ is a connected $C^1$-curve and
 $\gamma_{n,m}\cap S_i=\emptyset$
 for $i=0,\ldots,m-1$.
 \item 
$c_\eta(\gamma_{n, m}, \gamma_{0})\leq c$ for every
$n\geq 1$ and $m\geq m_2+(n-1)\pi(q)$.
\end{itemize}
\end{prop}

Next lemma will be one of our main tool: 
We give a bound of the cost of curves in a simple situation.
In the following by $\DD$ we denote the unit disc 
$B(1) \subset \mathbb{R}^2$. 

\begin{lemm} \label{l.boufra}
Let $0<\delta <1$ be given. 
Given $\alpha>0$ and $\varepsilon >0$, there 
is $N\in \NN$ satisfying the following property: 
Suppose $f\colon [-1, 1] \to\RR$ is a $C^1$-map satisfying
\begin{itemize}
\item $f(t)=0$ if $t \in [-1, -1+\delta] \cup [1-\delta, 1]$,
\item $|f'(t)| < \alpha$ for every $t$.
\item The graph of $f$, 
that is, $\{(t, f(t)) \in \RR^2 \mid t \in [-1, 1]\}$ 
is contained in $\DD$.
\end{itemize}
Then there is a sequence of 
diffeomorphisms $\{\varphi_i \}_{i = 1,\ldots, K}$ 
where $K\leq N$ and satisfying the following:
\begin{itemize}
\item For every $i=1,\ldots, K$, 
the support of $\varphi_i$ is contained in 
$\DD \setminus \partial \DD$ and has diameter less than $\varepsilon$, 
\item for every $i$, the $C^1$-distance 
between $\varphi_i$ and the identity map is 
smaller than $\varepsilon$, and
\item$\varphi_K\circ\cdots\circ \varphi_1([-1, 1]\times \{0\})$ is 
equal to the graph of $f$.
\end{itemize}
\end{lemm}

\begin{proof} Let $\theta\colon \mathbb{R}\to [0,1]$ be a 
smooth bump function satisfying the following:
\begin{itemize}
\item $\theta$ is equal to zero in 
$(-\infty, -2+\delta_0] \cup [2-\delta_0, +\infty)$
for some small $\delta_0>0$. 
\item $\theta$ is identically $1$ on $[-1,1]$.
\end{itemize}
For every sufficiently large $n\in\NN$ and every $i \in\{0,\dots,n^3-1\}$, we denote by $\psi_{i,n}$ 
the diffeomorphism of $\DD$ defined as the time one map of the vector field
\[X_{i,n}(x,y)=  \theta(n(y-\frac i{n^3}f(x)))\frac{f(x)}{n^3}\frac\partial{\partial y}.
\]

Note that $\psi_{i,n}(x, \frac i{n^3} f(x))=(x, \frac {i+1}{n^3} f(x))$, 
so we have $\psi_{n^3-1,n}\circ\cdots\circ \psi_{0,n} \big( (x,0) \big)= (x,f(x))$ for every $x \in [-1, 1]$.
Furthermore, the support of $\psi_{i,n}$ is contained in 
the $2n^{-3}$-neighborhood of $\{(x,\frac i{n^2} f(x)) \mid x\in[-1+\delta,1-\delta]\}$. 
A simple calculation shows that the $C^1$-distance 
of $\psi_{i,n}$ from the identity map 
is bounded by $K \alpha n^{-2}$ (where $\alpha$ is the constant 
in the statement of this lemma and $K>0$ is some constant independent of 
$f$), 
hence tends to $0$ when $n\to \infty$. 

However, the diameter of the support of $\psi_{i,n}$ does not tends to $0$. In order to obtain 
this property, we take finer factorization 
in products of diffeomorphisms with smaller support. 

We put
\[\theta_{j,n}(x)=  \frac{\theta(n(x-\frac jn))}{\sum_{k=-\infty}^{+\infty} \theta(n(x-\frac kn))}.\]
The family $\{\theta_{j,n}\}$ is a partition of 
unity whose differential is proportional to $n$. 
Then, we define diffeomorphisms $\{ \varphi_{i,j,n} \}$ 
as  time one map of the vector field 

$$X_{i,j,n}(x,y)= \theta_{j,n}(x) X_{i,n}(x,y).$$

Given a pair $(i,n)$, each $\varphi_{i,j,n}$ commutes
for restricted to each line $\{x=x_0\}$ each 
$\varphi_{i,j,n}$ is a flow generated by proportional
vector fields. 
Let $\rho_{i,n}$ be the product of  
$\{\varphi_{i,j,n}\}$ for $j\in\{-2,\dots, n +2\}$. 
Then we have 
$\rho_{n^3-1, n}\circ\cdots\circ \rho_{0, n} \big( (x,0) \big)= (x,f(x))$
for every $x \in [-1, 1]$,
the $C^1$ distance between $\varphi_{i,j,n}$ and 
the identity map is bounded by a constant proportional 
to $K'\alpha n^{-1}$ (where $K' >0$ is some constant independent of 
$f$) and 
the diameter of their support is bounded by $4n^{-1}$.
Thus the proof is completed.
\end{proof}

\subsection{Choice of the curves: Proof of Proposition~\ref{p.monoprix}}
\label{ss.choix}

Using Lemma~\ref{l.boufra}, we can complete the proof of 
Proposition~\ref{p.monoprix}.

\begin{proof}[Proof of Proposition~\ref{p.monoprix}]
Let $\{(\cD, F_n)\}_{ n\geq 1}$ be 
a prepared family for a periodic point $q$ and its 
homoclinic point $Q$.
Let $\lambda$ be the homothetic 
factor of the retarded family $(F_{n,q}^{\pi(q)})$.
Also, let $\eta >0$ be given.

In order to construct the family $\{\gamma_{n, m}\}$,
given $\eta'>0$
we only need to construct a family of curves 
$\{\alpha_{n, m}\}_{n \geq 1, m \geq m_3+\pi(q)(n-1)}$ 
(where $m_3$ is some non-negative integer)
in $D_q$ such that the following holds:
\begin{itemize}
\item Let $\alpha_0$ be the diameter $D_q \cap \{x=y\}$.
Every $\alpha_{n, m}$ coincides with $\alpha_0$ near 
$\partial D_q$.
\item $\alpha_{n, m} \cap \bar{\Lambda}_i$
is a connected curve for every $i=0,\ldots,m$.
\item $\alpha_{n, m} 
\cap \Delta_i = \emptyset$
for $i=0,\ldots, m-1$.
\item $\alpha_{n, m} 
\cap (F_{n, q}^{\pi(q)})^i(\Xi_1) = \emptyset$
for $i=0,\ldots, [m/\pi(q)]-1$.
\item $\alpha_{n, m} \cap \bar{\Lambda}_m$ 
coincides with the $x$-axis. 
\item There exists an integer $\eta'>0$ such that 
\[
c_{\eta'}(\alpha_{n,m}, \alpha_0) \leq c_{\eta'}
\]
for every $n \geq 1$ and $m \geq m_3 +\pi(q)(n-1)$.
\end{itemize} 
If we have constructed such a family, then we can 
obtain the conclusion. To see this, recall that there 
is an integer $\mkd$ 
such that $T_{\mkd} = F_{n,q}^{\mka}\circ 
F_{n,Q}^{\mkt} \circ F_{n,q}^{\mkd}(D_q)$ holds
(see Remark~\ref{r.selfsi}). If we take the image of
the family $\{\alpha_{n, m}\}$ under 
$F_{n,q}^{\mka}\circ F_{n,Q}^{\mkt} \circ F_{n,q}^{\mkd}$,
it gives a family of curves in $T_{\mkd}$ such that 
for each $n$, the $K_n\eta'$-cost
between $F_{n,q}^{\mka}\circ F_{n,Q}^{\mkt} \circ F_{n,q}^{\mkd}(\alpha_0)$ and 
$F_{n,q}^{\mka}\circ F_{n,Q}^{\mkt} \circ F_{n,q}^{\mkd}(\alpha_{n, m})$ are 
uniformly bounded, 
where $K_n>0$ is some constant determined by 
$F_{n,q}^{\mka}\circ F_{n,Q}^{\mkt} \circ F_{n,q}^{\mkd}$.
Note that $K_n$ depends on $n$ but $F_{n,Q}^{\mkt}$
are independent of $n$ and the $C^1$-norm of 
$F_{n,q}$ is uniformly bounded. 
Thus the sequence $(K_n)$ is also uniformly bounded, 
for $(F_n)$ is uniformly bounded 
by assumption. As a result, 
the $K\eta'$-cost of the family 
$F_{n,q}^{\mka}\circ F_{n,Q}^{\mkt} \circ F_{n,q}^{\mkd}(\alpha_{n, m})$ is uniformly bounded by 
$c_{\eta'}$, where $K$ is the uniform bound of $(K_n)$.

Let us explain how to extend these curves to obtain $\{\gamma_{n,m}\}$ and
$\gamma_0$.
We extend the curve $F_{n,q}^{\mka}\circ F_{n,Q}^{\mkt} \circ F_{n,q}^{\mkd}(\alpha_0)$  to 
a $C^1$-curve $\gamma'_0$ in  
$\Xi_1$ which 
coincides with the $x$-axis near the boundary of $D_q$. 
Since the objects in $\Xi_1$ outside $T_{\mkd}$ 
(such as $\Theta_i$ for $i=1,\ldots, \mkt -1$ and 
$S_i$ for $i=1,\ldots, \mkd -1$) 
are the same for every $n$, we can extend
each image $F_{n,q}^{\mka}\circ F_{n,Q}^{\mkt} \circ F_{n,q}^{\mkd}(\alpha_{n, m})$ to $\Xi_1$ 
such that the following holds:
\begin{itemize}
\item Each extension coincides with $\gamma_0$ near the boundary 
of $\Xi_1$,
\item The $\eta$-cost between each extension and 
$\gamma'_0$ is uniform bounded. Note that the same is 
true for the cost between each extension and $\gamma_0$.
\item The intersections between each extension and 
$\Xi_i$ ($i=0,\ldots,\mkt -1$) and $T_i$
($i=0,\ldots,\mkd$) are connected.
\end{itemize}
Thus the extension gives the desired family 
$\{\gamma_{n,m}\}$.

\bigskip

Now, given $\eta'>0$ 
let us construct the family of curves 
$\{\alpha_{n, m}\}$ in $D_q$ which has uniform 
$\eta'$-cost from $\alpha_0$.
We only need to construct a family $\{\alpha_{1, m}\}_{m\geq m_3}$ 
since for general $\{\alpha_{n, m}\}$ we only need to extend the homothetic 
image of $\{\alpha_{1, m}\}_{l\geq m_3}$ 
by a straight line in the homothetic region.
Note that by increasing $n$ by $1$ the number of the 
homothetic region increases by $1$, thus the number 
of intermediate strata increase by $\pi(q)$.
Thus $\{\alpha_{n,m}\}$ is defined only for 
$m \geq m_3 + \pi(q)(n-1)$.

Let $\tau$ be an integer in (P2) (see Definition~\ref{d.preparation}).
Recall that we have rectangles
$\beta_k \subset \Lambda_{\tau +k} \setminus 
\Lambda_{\tau +k+1}$ for $k \geq 0$. 
%
Now, we construct 
a family of curves $\{\zeta_i\}_{i \geq 0}$
in $\bar{\Lambda}_{\tau}$ satisfying the following conditions 
(see Figure~\ref{f.presol-curves}):
\begin{itemize}
\item $\beta_{k} \cap \zeta_i = \emptyset$ if 
$k < i$.
\item $\beta_{k} \cap \zeta_i$ coincides with the $x$-axis 
if $k \geq i$.
\item For each $i$, the two endpoints of $\zeta_i$ are  
$\Lambda_{\tau} \cap \{x=y\}$.
\item For every $i \geq 0$ and  $j \geq \tau$, 
$\zeta_i \cap \bar{\Lambda}_j$ is a connected 
curve.
\end{itemize}

\begin{figure}[h]
\begin{center}
\includegraphics[width=10cm]{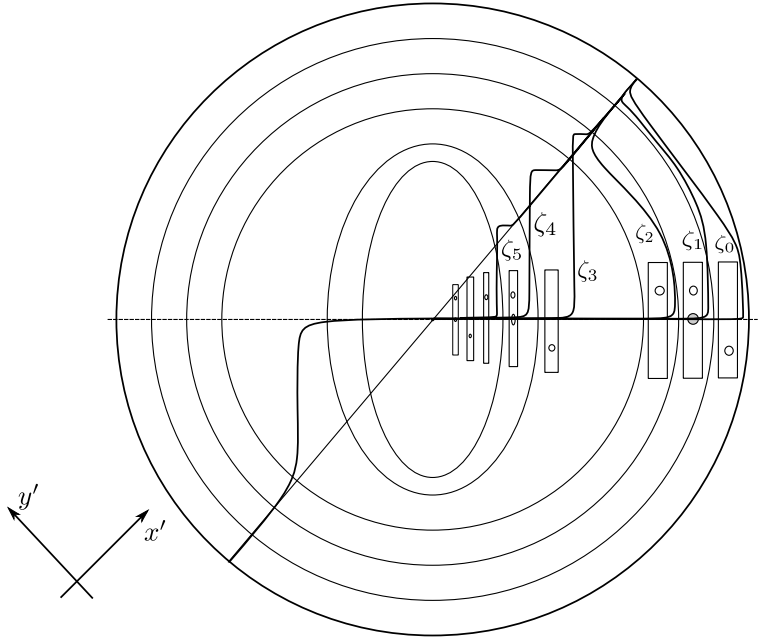}
\caption{A graphical explanation the family $\{\zeta_i\}$.
They are chosen in such a way that they have bounded 
Lipschitz constants with respect to the orthogonal 
$x'y'$-coordinates
where $\{x=y\}$ corresponds to the $x'$-axis, except 
some finitely many $\zeta_i$. In the picture, 
$\zeta_0, \zeta_1, \zeta_2$ are the exceptions.
They need to turn back to reach the line $\{x=y\}$ and 
may fail to be graphs in $x'y'$-coordinates.}
\label{f.presol-curves}
\end{center}
\end{figure}

Let us prove the following:
\begin{clai}
$\{ \zeta_i\} \subset \Lambda_{\tau}$ can be chosen such that the $\eta'$-cost from 
$\zeta_{-\infty} := \bar{\Lambda}_{\tau} \cap \{x=y\}$
is uniformly bounded. 
\end{clai}
\begin{proof}
We need to achieve two properties: 
$\zeta_i$
avoids the intersection 
with $\beta_k$ for $k < i$ and they must have intersection 
for $i \geq k$. 
The other condition is that the intersection with $\zeta_i$
and $\bar{\Lambda}_j$ is connected. 

To obtain both, 
we choose the family of curves $\{\zeta_i\}$ described 
as in Figure~\ref{f.presol-curves}. Namely, $\zeta_i$ is a
curve such that:
\begin{itemize}
\item In the positive $x$-half plane,
first it follows the $x$-axis.
\item Then it makes an almost vertical turn to
the $y$-direction
between $\beta_k$ and $\beta_{k-1}$.
\item After the curve reaches higher than all the 
rectangles $\{\beta_i\}$,
it makes another almost vertical turn to the $x$-direction.
Recall that we require that $\beta_i$ ($i=0,\ldots \pi(q)-1$)
is disjoint from the line $\{x=y\}$.
\item When it arrives at the line $\{x=y\}$, 
the curve follows it to reach $\Lambda_{\tau}$.
\item On the negative $x$-half plane 
which contains no $\beta_k$, 
we just take 
some extension, paying attention to keep the condition 
about the connectedness. 
\end{itemize}
We can construct such $\zeta_i$ for sufficiently large $i$.
If we introduce orthogonal $x'y'$-coordinate system 
on $\Lambda_{\tau}$
which sends $\zeta_{-\infty}$ 
to the $x'$-axis, then it is not difficult 
to see that we can realize $\zeta_i$ as the graph
of functions whose derivatives are uniformly bounded. 
Note that in the diagonal region the map $F_{1, q}^{\pi(q)}$
has the diagonal form. It guarantees that if the turns of
$\zeta_i$ are enough vertical then the connectedness of 
the curves $\bar{\Lambda}_j \cap \zeta_i$ are guaranteed.

Now, Lemma~\ref{l.boufra} implies the uniform boundedness 
of the cost from $\zeta_{-\infty}$ to $\zeta_i$.
There may be some curves where this construction does not 
hold, but there are at most finitely many such curves and their contributions
are irrelevant to the boundedness of the cost.
\end{proof}

Then we need to connect $\{\zeta_i\}$ to 
$\Lambda_0$. The shape of $(\Delta_i)$ and 
$(F_{1,q}^{\pi(q)})^j(\Xi_1)$ may be complicated in
$\bar{\Lambda}_0 \setminus \bar{\Lambda}_\tau$ 
but the size of this region is bounded. Thus, 
we can extend $\{\zeta_i\}$ to $\{\alpha_{1, m}\}$
keeping the uniform boundedness of the $\eta'$-cost. 
\end{proof}

\subsection{Proof of Proposition~\ref{p.solutionprepared}}
\label{ss.finale}
In this subsection we finish the proof of 
Proposition~\ref{p.solutionprepared}.
\begin{proof}
Let $\{(\cD, F_n)\}_{n \geq 1}$ be a prepared family and 
fix $\varepsilon_2>0$. 
We take an integer $c>0$ and the curves $\gamma_0$ and 
$\{\gamma_{n,m}\}_{n\geq 1, m \geq  m_2 + \pi(q)(n-1)}$ in 
$\Xi_1$ 
by Proposition~\ref{p.monoprix}. We have 
$c_{\varepsilon_2}(\gamma_{m, n}, \gamma_0 )< c$ 
for every $n \geq 1$ and $m \geq m_2+\pi(q)(n-1)$. 

By definition of the $\varepsilon_2$-cost, for each $m \geq m_2+\pi(q)(n-1)$ there are 
diffeomorphisms $\varphi_{1, m},\dots,\varphi_{c, m}$ 
supported in the interior of 
$\Xi_1$, $\varepsilon_2$-$C^1$-close to the identity such that 
$$\varphi_{c, m}\circ\cdots\circ\varphi_{1,m}(\gamma_{n, m})=
\gamma_0.$$

Fix now an integer $n_0$ greater than $c+1$. We compose 
$\varphi_{i, m}$
to $F_{n_0,F^{-1}(q)}$ 
as in the proof of Proposition~\ref{p.cleaning}.
Namely, for $i\in\{1,\dots,c\}$ let $\psi_i$ be a diffeomorphism 
of $D_q$ which satisfies
$\psi_{i, m}=(F_{n_0,q}^{\pi(q)})^i\circ \varphi_{i, m}\circ (F_{n_0,q}^{\pi(q)})^{-i}$ on 
$(F_{n_0,q}^{\pi(q)})^i(D_q)$
and equal to the identity map outside. 
Note that $\psi_{i, m}$ 
is supported on $(F_{n_0,q}^{\pi(q)})^i(\Xi_1)$, 
which is contained in the annulus 
$(F_{n_0,q}^{\pi(q)})^i(D_q)
\setminus (F_{n_0,q}^{\pi(q)})^{i+1}(D_q)$.
Remark that 
this is contained in the homothetic region of $F_{n_0}$. 

Notice that the $C^1$-distance between $\psi_{i, m}$ and the identity map is bounded by $\varepsilon_2$, as $\psi_{i, m}$ is conjugated 
to $\varphi_{i, m}$ by a contracting homothety. 
The maps $\{\psi_{i, m}\}_{i=1,\ldots, c}$ have disjoint support so they commute.  Let $\psi_m$ be the product $\psi_m=\psi_{c, m}\circ\cdots\circ\psi_{1,m}$ and 
let $G_m$ be a Markov IFS obtained by composing
$\psi_m$ to $F_{n_0,F^{-1}(q)}$ and keep the other maps
intact.

The map $G_m$ is an admissible $\varepsilon_2$-$C^1$-small perturbation of $F_{n_0}$. 
Let us check the following:
\begin{clai}
 $G_m$ is a pre-solution of depth $m$.
\end{clai}
\begin{proof}
Let us see that 
\[ W^{ss}(q, G_{m, q}^{\pi(q)}) \cap \Xi_1 = \gamma_{n_0, m}. \] 
Then, the fact that $G_m$ is an admissible perturbation of 
$F_{n_0}$ and the definition of $\gamma_{n_0, m}$ immediately
implies the conclusion.

Consider a point $\bar x\in\Xi_1$ and $k>c$. Then 
\begin{align*}
(G_{m, q}^{\pi(q)})^k(\bar x)
&=(F_{n_0,q}^{\pi(q)})^{k-c} \circ \left[\prod_{i=1}^c  
[(F_{n_0,q}^{\pi(q)})^{i}\circ \varphi_{i,m}\circ (F_{n_0,q}^{\pi(q)})^{-i} ] \right]\circ F_{n_0,q}^{\pi(q)}(\bar x) \\
&=(F_{n_0,q}^{\pi(q)})^k\circ\varphi_{c,m}\circ\cdots\circ \varphi_{1,m}(\bar x).
\end{align*}
Therefore, 
for $\bar x\in\Xi_1$, 
we have $\bar x\in  W^{ss}(q,G_{m, q}^{\pi(q)})$ if and only if 
$\varphi_{c,m}\circ\cdots\circ \varphi_{1,m}(\bar x)$ belongs to
$W^{ss}(q,F_{n_0,q}^{\pi(q)}) \cap \Xi_1
=\gamma_0$, 
that is, $\bar x\in \gamma_{n_0,m}$.
\end{proof}
Thus the proof of Proposition~\ref{p.solutionprepared} is completed.
\end{proof}

\subsection{Weakness of the invariant curves}
\label{ss.finfin}
Now we are ready to finish the 
proof of Theorem~\ref{t.IFS}. 
Let us complete it.

\begin{proof}
Let a Markov IFS $(\cD, F)$ with an $\varepsilon$-flexible
point $q$ having a large stable manifold and its $u$-homoclinic
point $Q$ be given. Also, let $\eta >0$ and $\varepsilon_0>0$
be given.

First, we apply Lemma~\ref{l.flexihomoclinic} and Proposition~\ref{p.cleaning} successively. Then 
we obtain a bounded prepared family
$\{(\cD, F_n)\}$ such that each $F_n$ is a $C^1$-$\varepsilon$-perturbation of $F$ and the $C^0$-distance between 
$F_n$ and $F$ is less than $\varepsilon_0$.


Then we apply Proposition~\ref{p.solutionprepared} to 
$(F_n)$ letting $\varepsilon_2>0$ small.
Then we obtain $n_0$, $m_0$, $W_{n_0}$ and a
family of pre-solutions
$(G_{m})_{m \geq m_0}$ of depth $m$ which 
are $\varepsilon_2$-admissible perturbations of $F_{n_0}$.
Each $G_{m}$ has a family of univalent invariant 
curves $\Gamma_m$ in $(0, m)$-refinement 
by Proposition~\ref{p.preso}. 
Note that $(G_m)$
is bounded, too.
Since $\varepsilon_2$ can be chosen arbitrarily small, 
we see that $(G_m)$ is $C^1$-$\varepsilon$-close to $F$
and $C^0$-$\varepsilon_0$-close to $F$.
Note that by shrinking $W_{n_0}$, 
we may assume that 
the differential of 
$F_{n, q}^{\pi(q)}|_{D_q\cap W_{n_0}}$ on the $x$-axis 
is equal to one.

In the following, we will choose convenient $m$ so that 
in the resulted dynamics this curve has two kinds of 
hyperbolicity in the definition of contracting
invariant curves, see Definition~\ref{d.nhiv}, \ref{d.ntan} and \ref{d.nnor}.
We only explain how to establish the weakness
of the curves. The choice of 
$m$ for the contraction in the tangential 
direction is left to the reader. 

We choose 
\begin{align*}
M_{+} &= \sup_m \{ 
\max_{x \in \Gamma_m} \|D(G_{m})^{-1}(x)|_{T\cD/T\Gamma_m}\|, 
\max_{x \in \Gamma_m}\|D(G_{m})^{-1}|_{T\cD/T\Gamma_m}\|, 1\},\\
M_{- } &= \inf_m \{\min_{x \in \Gamma_m} \|D(G_{m})^{-1}(x)|_{T\cD/T\Gamma_m}\|, \min_{x \in \Gamma_m} \|D(G_{m})^{-1}(x)|_{T\cD/T\Gamma_m}\|, 1 \}.
\end{align*}
Notice that, even though 
there are infinitely many maps $G_{m}$ in the 
argument of the supremum and the infimum, 
 $M_{+}$ and $M_{-}$ are positive finite values, for 
$G_{m}$ are bounded. 
Then we apply Proposition~\ref{p.frequ} to 
$F_{n_0}$, $W$ and some $r$. Then for every 
$G_{m}$ the conclusion of Proposition~\ref{p.frequ} holds.
In the following, we consider the case where the 
first conclusion of Proposition~\ref{p.frequ} holds,
that is, there is one connected interval $H$
satisfying the conditions. 
The proof for the second case is similar, so we omit the 
argument of that case.

Consider a point $x \in \Gamma_m$ for which 
$(G_{m})^{-m}(x)$ is defined. 
Since the normal expansion of $\Gamma_m$
is $1$ in $W$,
$x \in \Gamma_m$ and $x$ goes around $\orb(q)$ at least
$([mr/\omega(q)] -2)$-times (where $[x]$ denotes the 
integer part of $x$), we know 
\begin{align*}
&M_{-}^{m-([mr/\omega(q)] -2)\omega(q)} 
\leq \|D(G_{m})^{-m}|_{TD/T\Gamma_m}(x)\| 
\leq 
M_{+}^{m-([mr/\omega(q)] -2)\omega(q)}. 
\end{align*}
Thus, by taking $r$ arbitrarily close to $1$, 
we see that the contribution of the derivatives outside 
$W$ will be negligible. Thus we have that 
the average normal derivative of $\Gamma_m$ tends 
to $1$ by letting $r$ close to one.
\end{proof}

\section{On the proof of general cases}
\label{s.expla}
In this section, we discuss the proof of Theorem~\ref{t.main}
based on the argument of the proof of Theorem~\ref{t.IFS}. 
Indeed, the proof can be done in a parallel way, adding some 
modification to avoid the interference. Let us briefly 
see it.

Let $\{q_i\}$ and $\{Q_j\}$ be given. 
\begin{itemize}
\item First, 
we apply Lemma~\ref{l.flexihomoclinic}
to each $q_i$. It gives us a retarded family for each $q_i$.
Note that there is no interference between two different 
periodic points, due to the assumption of 
mutually separatedness. 
\item Then, for each $q_i$, we consider the homo/heteroclinic
points which arrives at $\mathrm{orb}(q_i)$. For each homo/heteroclinic 
point we define the periodic 
well and the transition well in the first fundamental domain. 
\item Note that two transition wells of different homo/heteroclinic orbits may share some discs, 
since two different homo/heteroclinic orbits may share
their itinerary.
\item However, they cannot have totally the same itinerary
by assumption.
Thus at least we know that their periodic wells are disjoint. 
\item Now we give small perturbation to each retarded family 
which makes it to be prepared. Then for each periodic well 
we find curves $\{\gamma_{n, m}\}$ with bounded cost.
\item Since the periodic wells are disjoint, we have 
no interference when we perform the perturbation in 
Section~\ref{ss.finale}. Thus we can obtain the pre solution 
of arbitrarily profound depths for each homo/heteroclinic orbits.
\end{itemize}

\bigskip
\bigskip
\bigskip

\begin{itemize}
\item[]  \emph{Christian BONATTI (christian.bonatti@u-bourgogne.fr)}
\begin{itemize}
\item[] Institut de Math\'ematiques de Bourgogne CNRS - URM 5584
\item[] Universit\'e de Bourgogne Dijon 21004, France
\end{itemize}
\item[] \emph{Katsutoshi SHINOHARA (ka.shinohara@r.hit-u.ac.jp)}
\begin{itemize}
\item[]  Graduate School of Business Administration, 
\item[]  Hitotsubashi University,
\item[]  2-1 Naka, Kunitachi, Tokyo 186-8601, Japan
\end{itemize}
\end{itemize}

\end{document}